\documentclass[11pt]{article}%
\usepackage{graphicx}
\usepackage{amsmath}
\usepackage{fullpage}
\usepackage{amsfonts}
\usepackage{amssymb}%
\providecommand{\U}[1]{\protect\rule{.1in}{.1in}}
\newtheorem{theorem}{Theorem}

\newtheorem{corollary}[theorem]{Corollary}

\newtheorem{definition}[theorem]{Definition}
\newtheorem{example}[theorem]{Example}

\newtheorem{lemma}[theorem]{Lemma}
\newtheorem{notation}[theorem]{Notation}
\newtheorem{problem}[theorem]{Problem}
\newtheorem{proposition}[theorem]{Proposition}
\newtheorem{remark}[theorem]{Remark}

\newenvironment{proof}[1][Proof]{\textbf{#1.} }{\ \rule{0.5em}{0.5em}}
\numberwithin{theorem}{section} \numberwithin{equation}{section}
\begin{document}

\title{Relation and radical approach to the theory of C*-algebras}
\author{Edward Kissin, Victor S. Shulman and Yurii V. Turovskii}
\maketitle

\begin{abstract}
In this paper we pursue three aims. The first one is to apply Amitsur's
relations and radicals theory to the study of the lattices Id$_{A}$ of closed
two-sided ideals of C*-algebras $A$. We show that many new and many well-known
results about C*-algebras follow naturally from this approach.

To use "relation-radical" approach, we consider various subclasses of the
class $\mathfrak{A}$ of all C*-algebras, which we call C*-\textit{properties}%
,\textit{ }as they often linked to some properties of C*-algebras. We consider
C*-properties $P$ consisting of $CCR$- and of $GCR$-algebras; of C*-algebras
with continuous trace; of real rank zero, AF, nuclear C*-algebras, etc. Each
$P$ defines reflexive relations $\ll_{_{P}}$ in all lattices Id$_{A}.$ Our
second aim is to determine the hierarchy and interconnection between
properties in $\mathfrak{A}.$ 

Our third aim is to study the link between the radicals of relations
$\ll_{_{P}}$in the lattices Id$_{A}$ and the topological radicals on
$\mathfrak{A}.$

\end{abstract}

\section{Introduction and preliminaries}

In his research of the radical theory of algebras and rings, Amitsur
\cite{Am2}, \cite{Am3} discovered that a significant part of the results can
be formulated and proved in terms of the general theory of lattices. In
\cite{Am} he developed the theory of radicals for relations in lattices that
was used in various areas of algebra: group theory, non-associative rings, Lie
algebras, universal algebras, etc. This theory was further developed in the
work of Kurosh \cite{K}.

Later Dixon \cite{Di} initiated the radical approach to some problems of
functional analysis and laid the basis of the theory of topological radicals
of Banach algebras. This theory was further developed and applied to the
theory of invariant subspaces of operator algebras and to classification of
Banach and operator Lie algebras in \cite{ST}, \cite{TR}, \cite{KST1}.

In this paper we pursue three aims. The first one is to apply Amitsur's
relations and radicals theory to the study of the lattices Id$_{A}$ of closed
two-sided ideals of C*-algebras $A$. We show that many new and many well-known
results about C*-algebras follow naturally from this approach. For example, it
is well known that each $A$ has the largest $GCR$ ideal $I_{A}$ and an
ascending transfinite chain $\{I_{\lambda}\}$ of ideals such that all
$I_{\lambda+1}/I_{\lambda}$ are $CCR$ algebras and $\overline{\cup I_{\lambda
}}=I_{A}.$ From the "relation-radical" approach this result follows
immediately if we consider the relation $\ll_{_{CCR}}$ in Id$_{A}$; the ideal
$I_{A}$ coincides with the radical $\mathfrak{r}_{_{CC\mathsf{R}(P)}%
}^{\triangleright}$ of the relation $\ll_{_{CCR}}^{\triangleright}$
constructed from $\ll_{_{CCR}}.$ In addition, we obtain that if $I_{A}%
\nsubseteq J\in$ Id$_{A}$ and $J\neq A$ then there is $K\in$ Id$_{A}$ such
that $J\subsetneqq K$ and $K/J$ is a $CCR$ algebra. All these results follow
from much more general results (Theorem \ref{T3.1}).

To use "relation-radical" approach, we consider various subclasses of the
class $\mathfrak{A}$ of all C*-algebras, which we call C*-\textit{properties}%
,\textit{ }as they often linked to some properties of C*-algebras. We consider
C*-properties" consisting of $CCR$- and of $GCR$-algebras; of C*-algebras with
continuous trace; of real rank zero, AF, nuclear C*-algebras, etc. Each
C*-property $P$ defines reflexive relations $\ll_{_{P}}$ in all lattices
Id$_{A}$ by (\ref{4.1})$.$ The collection $\mathcal{P}$ of all properties in
$\mathfrak{A}$ is also a lattices. Our second aim is to determine the
hierarchy and interconnection between properties in $\mathcal{P}.$ We
introduce and investigate three maps (closure operators) in $\mathcal{P}$:
$\mathsf{G}$, $\mathsf{dG}$ and $\mathsf{R}(P).$ For example, for the
C*-property $C$ of all C*-algebras isomorphic to $C(\mathcal{H})$ on separable
Hilbert spaces $\mathcal{H}$, the C*-property $C\mathsf{R}(P)$ coincides with
the class of all $CCR$ algebras and the C*-property $\mathsf{G}(C\mathsf{R}%
(P))$ coincides with the class of all $GCR$ algebras. In some cases the
properties $P$ and $\mathsf{G}P$ coincide. For example, they coincide if $P$
is the C*-property of all AF, or the C*-property of all nuclear algebras.

Our third aim is to study the link between the radicals of relations generated
by properties $P$ in the lattices Id$_{A}$ and the topological radicals on
$\mathfrak{A}.$

Amitsur's theory is based on the investigation of special relations in
complete lattices $(Q,\leq)$: $\mathbf{H}$-relations and dual $\mathbf{H}%
$-relations (see (\ref{1.2}) and (\ref{dh1})). In general, they do not have
radicals. Amitsur developed a certain procedure for constructing an
$\mathbf{R}$-order $\ll^{\triangleright}$ from an $\mathbf{H}$-relation $\ll$
and proved that $\ll^{\triangleright}$ has a unique $\ll^{\triangleright}%
$\textit{-radical }$\mathfrak{r}^{\triangleright}$ in $Q$ (see (\ref{2.11})).
Similarly, a dual $\mathbf{R}$-order $\ll^{\triangleleft}$ can be constructed
from a dual $\mathbf{H}$-relation $\ll$ and $\ll^{\triangleleft}$ has a unique
dual $\ll^{\triangleleft}$\textit{-radical }$\mathfrak{p}^{\triangleleft}$
(see (\ref{2.11})). This procedure was further refined in \cite{KST2} and we
describe it briefly in Section 2.

Simultaneously $\mathbf{H}$- and dual $\mathbf{H}$-relations ($\mathbf{HH}%
$-relations) play an important role in the theory of operator algebras $A$ on
Banach spaces$.$ For example, the relation $\ll_{_{\infty}}$ defined in the
lattice Lat $A$ of all invariant subspaces by $L\ll_{_{\infty}}M$ if
$L\subseteq M$ in Lat $A$ and $\dim(M/L)<\infty,$ is an $\mathbf{HH}$-relation
(\cite{KST1}, \cite{KST2}). The lattice Lat $A$ has the $\ll_{_{\infty}%
}^{\triangleright}$-radical $\mathfrak{r}_{_{\infty}}^{\triangleright}$ and an
ascending transfinite sequence $(\mathfrak{r}_{\lambda})$ of subspaces such
that $\mathfrak{r}_{_{\infty}}^{\triangleright}=\overline{\cup\text{
}\mathfrak{r}_{\lambda}}$ and $\dim\mathfrak{r}_{\lambda+1}/\mathfrak{r}%
_{\lambda}<\infty.$ It also has the dual $\ll_{_{\infty}}^{\triangleleft}%
$-radical $\mathfrak{p}_{_{\infty}}^{\triangleleft}$ and a descending sequence
$(\mathfrak{p}_{\lambda})$ of subspaces such that $\mathfrak{p}_{_{\infty}%
}^{\triangleleft}=\cap$ $\mathfrak{p}_{\lambda}$ and $\dim\mathfrak{p}%
_{\lambda}/\mathfrak{p}_{\lambda+1}<\infty.$ The radicals $\mathfrak{r}%
_{_{\infty}}^{\triangleright}$ and $\mathfrak{p}_{_{\infty}}^{\triangleleft}$
are invariant for the operator Lie algebra of all derivations of $A$. A
description of $\mathbf{HH}$-relations in the lattices of the projections in
W*-algebras was obtained in \cite{Ki}.

For $A\in\mathfrak{A},$ the set Id$_{A}$ of all closed two-sided ideals of $A$
(we call them ideals) is a complete lattice\textit{ }with\textit{ }relation
$\leq$ being the inclusion of ideals $\subseteq$.

\begin{definition}
\label{D1.1}\emph{(i) }We say that a subclass $P$ of $\mathfrak{A}$ is a
C*-\textbf{property}$\mathfrak{,}$ if%
\begin{equation}
\{0\}\in P\text{ and\emph{,} if }A\in P\text{ then all }B\in\mathfrak{A}\text{
isomorphic to }A\text{ }(B\approx A)\text{ belong to }P. \label{1.7}%
\end{equation}
If $A\in P,$ we say that $A$ has property $P,$ or that it is a $P$%
-\textit{algebra}. In the same way we use terms $P$-ideal, or $P$%
-quotient.\textbf{ }\smallskip

\emph{(ii) }For a C*-property $P$ in $\mathfrak{A}$ and $A\in\mathfrak{A},$
define the reflexive relation $\ll_{_{P}}$ in \emph{Id}$_{A}$ as follows%
\begin{equation}
I\ll_{_{P}}J\text{ if }I\subseteq J\text{ in \emph{Id}}_{A}\text{ and
}J/I\text{ is a }P\text{-algebra},\text{ i.e., }J/I\in P. \label{4.1}%
\end{equation}

\end{definition}

Note that when a C*-property is defined, $\{0\}$ is automatically included.
For example, if $P$ is the C*-propery of all one-dimensional algebras, then
$\{0\}\in P$ although $\{0\}$ is not one-dimensional.

We will consider a wide variety of properties $P$ in $\mathfrak{A}$. For each
of them$,$ we treat Id$_{A}$ as a lattice with two relations: $\subseteq$ and
$\ll_{_{P}}$ and, using methods of the lattice theory, analyze its structure.
We show that many results in the theory of C*-algebras can be obtained by this approach.

\begin{definition}
\label{D4.0}\emph{(i)\ }We say that a C*-property $P$ in $\mathfrak{A}$ is
\textbf{lower stable} if $A\in P$ implies \emph{Id}$_{A}\subset P.\smallskip$

\emph{(ii) }We say that $P$ is \textbf{upper stable} if $A\in P$ implies that
all quotients $A/I\in P$ for $I\in$ \emph{Id}$_{A}.$
\end{definition}

For example, the classes of all $CCR$- and of all $GCR$-algebras are lower and
upper stable properties. We will see later that many properties studied in the
theory of C*-algebras are lower, or upper stable, or both lower and upper
stable. If $P$ is an upper stable C*-property then $\ll_{_{P}}$ is an
$\mathbf{H}$-relation, if $P$ is lower stable then $\ll_{_{P}}$ is a dual
$\mathbf{H}$-relation (Theorem \ref{T4.1}). So the relations $\ll_{_{P}%
}^{\triangleright}$ and $\ll_{_{P}}^{\triangleleft}$ constructed from them are
$\mathbf{R}$-, or dual $\mathbf{R}$-orders and Id$_{A}$ has the $\ll_{_{P}%
}^{\triangleright}$-radical $\mathfrak{r}_{_{P}}^{\triangleright},$ or the
dual $\ll_{_{P}}^{\triangleleft}$-radical $\mathfrak{p}_{_{P}}^{\triangleleft
}.$ These radicals provide important information about the structure of $A.$
In particular, they are invariant for all automorphisms of $A$ (Corollary
\ref{C9.1}).

For each C*-property $P,$ we introduce a wider \textit{generalized} $P$
\textit{C*-property}: $\mathsf{G}P=\{A\in\mathfrak{A}$: $\mathfrak{r}_{_{P}%
}^{\triangleright}=A\}$. We show that the relation $\ll_{_{\mathsf{G}P}}$ is
an $\mathbf{R}$-order in Id$_{A}$, $\ll_{_{P}}^{\triangleright}=$
$\ll_{_{\mathsf{G}P}}$ and the radical $\mathfrak{r}_{_{P}}^{\triangleright
}=\mathfrak{r}_{_{\mathsf{G}P}}$ is the largest $\mathsf{G}P$-ideal in $A$ and
the smallest ideal such that $A/\mathfrak{r}_{_{P}}^{\triangleright}$ is an
$\mathsf{NG}P$-algebra, i.e., has no $P$-ideals. The radical $\mathfrak{r}%
_{_{P}}^{\triangleright}$ can be "approached" by an ascending transfinite
$\ll_{_{P}}$-series $(\mathfrak{r}_{\lambda})$ of ideals: $\mathfrak{r}_{_{P}%
}^{\triangleright}=\overline{\cup\text{ }\mathfrak{r}_{\lambda}}$ and
$\mathfrak{r}_{\lambda+1}/\mathfrak{r}_{\lambda}\in P.$ This is a natural
generalization of the well-known result for $CCR$ and $GCR$ algebras.
Moreover, if $\mathfrak{r}_{_{P}}^{\triangleright}\subsetneqq J\in$ Id$_{A}$
then $J$ is contained in $K\in$ Id$_{A}$ such that $K/J\in P.$

Dual results are also considered in Section 3, where we introduce \textit{the
dual} $\mathsf{G}P$ \textit{C*-property}: \textsf{d}$\mathsf{G}P=\{A\in
\mathfrak{A}$: $\mathfrak{p}_{_{P}}^{\triangleleft}=\{0\}\}$. For example,
$CCR$ algebras are \textsf{d}$GCR$ algebras. All residually finite-dimensional
algebras and, in particular, the group C*-algebras generated by the free
groups $\mathbb{F}_{n}$ on $n$ generators are \textsf{d}$GCR$ algebras
(Example \ref{E5}).

In Propositions \ref{L3.1} and \ref{P3.7} we show that the maps $\mathsf{G}$:
$P\mapsto\mathsf{G}P$ and \textsf{d}$\mathsf{G}$: $P\mapsto$ \textsf{d}%
$\mathsf{G}P$ are closure operators (order preserving, increasing idempotant
maps (see (\ref{2.6}))) on the subclasses of all upper stable and of all lower
stable properties in $\mathfrak{A,}$ respectively.

Different properties $P_{1},P_{2}$ may have equal generalized
properties:\ $\mathsf{G}P_{1}=\mathsf{G}P_{2}.$ In Section 5 we show that the
C*-property $P_{c.t.}$ of all continuous trace C*-algebras is upper and lower
stable. So the corresponding relation $\ll_{_{P_{c.t.}}}$ is an $\mathbf{HH}%
$-relation in all Id$_{A}.$ Although the properties $P_{c.t.}$ and $CCR$ are
different, $\mathsf{G}P_{c.t.}=GCR$ and their radicals coincide:
$\mathfrak{r}_{_{CCR}}^{\triangleright}=\mathfrak{r}_{_{P_{c.t.}}%
}^{\triangleright}$ and $\mathfrak{p}_{_{CCR}}^{\triangleleft}=\mathfrak{p}%
_{_{P_{c.t.}}}^{\triangleleft}$.

We also consider properties $P$ of all real rank zero C*-algebras, exact
C*-algebras, $AF$- and nuclear C*-algebras. We show that they are lower and
upper stable and, therefore, generate $\mathbf{HH}$-relations $\ll_{_{P}}$ in
each Id$_{A}.$ The corresponding $\ll_{_{P}}^{\triangleright}$-radicals
$\mathfrak{r}_{_{P}}^{\triangleright}$ in Id$_{A}$ contain all ideals of $A$
belonging to $P$ and there are ascending transfinite series $(\mathfrak{r}%
_{\lambda})$ of ideals such that $\mathfrak{r}_{_{P}}^{\triangleright
}=\overline{\cup\text{ }\mathfrak{r}_{\lambda}}$ and $\mathfrak{r}_{\lambda
+1}/\mathfrak{r}_{\lambda}\in P$.

Moreover, for the properties $P$ of all $AF$- and of all nuclear algebras,
$P=\mathsf{G}P.$ So $\ll_{_{P}}$ is an $\mathbf{R}$-order and the radical
$\mathfrak{r}_{_{P}}^{\triangleright}$ is a $P$-algebra (an $AF$- or a nuclear
C*-algebra, respectively).

In Section 4 we establish a link between topological radicals on C*-algebras
and the radicals of the relations $\ll_{_{P}}$ in the lattices Id$_{A}$. We
show that a map $R$: $A\mapsto R\left(  A\right)  \in$ Id$_{A}$ is a
topological radical on $\mathfrak{A}$ if and only if there is an upper stable
C*-property $P_{1}$ and a lower stable C*-property $P_{2}$ such that
$R(A)=\mathfrak{r}_{_{P_{1}}}^{\triangleright}(A)=\mathfrak{p}_{_{P_{2}}%
}^{\triangleleft}(A)$ for all $A.$ Moreover, $P_{1}=\mathbf{Rad}\left(
R\right)  =\left\{  A\in\mathfrak{M}\text{: }R\left(  A\right)  =A\right\}  $
and $P_{2}=\mathbf{Sem}\left(  R\right)  =\left\{  A\in\mathfrak{M}\text{:
}R\left(  A\right)  =0\right\}  .$

Each C*-property $P$ generates a relation-valued function $f_{_{P}}$:
$A\in\mathfrak{A\mapsto}$ $\ll_{_{P}}$ in Id$_{A}$ (see (\ref{4.1})). Let $f$:
$A\in\mathfrak{A\mapsto}\ll^{A}$ be a relation-valued function on
$\mathfrak{A},$ where $\ll^{A}$ is a reflexive relation in Id$_{A}%
\mathfrak{.}$ In Section 4 we give necessary and sufficient conditions for $f$
to be generated by a C*-property: $f=f_{_{P}}.$ If, in addition, each relation
$f(A)=$ $\ll^{A}$ has a unique $\ll^{A}$-radical $\mathfrak{r}(A)$ in Id$_{A}$
(see (\ref{2.11})), the question arises as to whether the map $\mathfrak{r}$:
$A\in\mathfrak{A}\mapsto\mathfrak{r}\left(  A\right)  $ is a topological
radical. In Theorem \ref{T4.5} we prove that $\mathfrak{r}$ is a topological
radical if and only if $f=f_{_{P_{f}}},$ where $P_{f}=\{A\in\mathfrak{A}$:
$\{0\}\ll^{A}A\}$\emph{ }is an upper stable C*-property and $P_{f}%
=\mathsf{G}P_{f}.$

In Section 6 we have exactly the above situation. Unlike the relations
constructed in the previous sections with the help of various properties in
$\mathfrak{A}$, we consider relations defined in all Id$_{A}$ by the property
of ideals vis-\`{a}-vis the algebras $A.$ Namely, an ideal $I$ of $A$ is
\textit{small} if $I+K\neq A$ for all $K\in$ Id$_{A},$ $A\neq K.$ For example,
the ideals $\{0\}$ and $C(H)$ are small in $B(H).$ The reflexive relation
$\ll_{\text{sm}}^{A}$ in Id$_{A}$ is defined by the condition: $I\ll
_{\text{sm}}^{A}J$ if $I\subseteq J$ and $J/I$ is a small ideal in $A/I.$

We show that $\ll_{\text{sm}}^{A}$ is a transitive $\mathbf{H}$-relation in
each Id$_{A},$ so that the relation $(\ll_{\text{sm}}^{A})^{\triangleright}$
is an $\mathbf{R}$-order and Id$_{A}$ has the $(\ll_{\text{sm}}^{A}%
)^{\triangleright}$-radical $\mathfrak{r}_{\text{sm}}^{\triangleright}\left(
A\right)  .$ However, unlike the maps $\mathfrak{r}_{_{P}}^{\triangleright}$:
$A\mapsto\mathfrak{r}_{_{P}}^{\triangleright}\left(  A\right)  $ generated by
properties $P$, the map $\mathfrak{r}_{\text{sm}}^{\triangleright}$:\emph{
}$A\mapsto\mathfrak{r}_{\text{sm}}^{\triangleright}\left(  A\right)  $\emph{
}is not a topological radical on\emph{ }$\mathfrak{A}$ (Remark \ref{R3}). If
$A$ is unital then $\ll_{\text{sm}}^{A}$ is an $\mathbf{R}$-order and
$\mathfrak{r}_{\text{sm}}(A)=\mathfrak{r}_{\text{sm}}^{\triangleright}(A)$ is
the largest small ideal of $A.$ We also establish that $\mathfrak{r}%
_{\text{sm}}(A)=$ rad$_{_{K}}$($A$), where rad$_{_{K}}$($A)$ is the radical
introduced by Kasch \cite{Kas} for rings and algebras. Aristov \cite{A}
extended Kasch's results to small submodules and morphisms in algebras.

\section{Radicals generated by H- and dual H-relations}

A partially ordered set $(Q,\leq)$ with a reflexive, anti-symmetric,
transitive relation $\leq$ is a \textit{lattice} if all $a,b\in Q$ have a
least upper bound $a\vee b$ and a greatest lower bound $a\wedge b.$ It is
$\vee$-\textit{complete}, if each subset $G\subseteq Q$ has a least upper
bound $\vee G;$ it is $\wedge$-\textit{complete }if each $G$ has a greatest
lower bound $\wedge G.$ A lattice is \textit{complete}, if it is\textit{
}$\vee$\textit{- and }$\wedge$\textit{-complete. }{In this case} we write
$\mathbf{0}=\wedge Q$ and $\mathbf{1}=\vee Q.$

Let $\ll$ be a \textit{reflexive} relation in a complete lattice $(Q,\leq)$
stronger than $\leq$: $a\ll b$ implies $a\leq b$ (in what follows we only
consider such relations). Amitsur \cite{Am} (see \cite{Gr}) defined
$\mathbf{H}$-relations and dual $\mathbf{H}$-relations in $Q$ as follows:%
\begin{align}
&  \ll\text{ is an }\mathbf{H}\text{-\textit{relation }if }a\ll b\text{ and
}a\leq c\text{ imply }c\ll b\vee c,\label{1.2}\\
&  \ll\text{ is a \textit{dual }}\mathbf{H}\text{-\textit{relation }if }a\ll
b\text{ and }c\leq b\text{ imply}\emph{\ }a\wedge c\ll c\text{,} \label{dh1}%
\end{align}
for $a,b,c\in Q.$ He proved that (\ref{1.2}) and (\ref{dh1}) are equivalent
respectively to the following conditions:%
\begin{align}
&  \ll\text{ is an }\mathbf{H}\text{-relation if }a\ll b\text{ implies }a\vee
x\ll b\vee x\text{ for each }x\in Q;\label{1.3}\\
&  \ll\text{ is a dual }\mathbf{H}\text{-relation if\ }a\ll b\text{ implies
}a\wedge x\ll b\wedge x\text{ for each }x\in Q.\nonumber
\end{align}
By the Duality Principle \cite[Theorem $1.3^{\prime}$]{Sk}, the results for
dual $\mathbf{H}$-relations follow from the corresponding results for
$\mathbf{H}$-relations and vice versa.

For $a,b\in Q$ and a reflexive relation $\ll,$ set%
\begin{equation}
\left[  a,\ll\right]  =\left\{  x\in Q\text{\textbf{: }}a\ll x\right\}
,\text{ }\left[  \ll,b\right]  =\left\{  x\in Q\text{: }x\ll b\right\}  \text{
and }[a,b]=\{z\in Q:a\leq z\leq b\}. \label{1}%
\end{equation}
A reflexive, transitive relation $\ll$ in $Q$ is%
\begin{align}
\text{an }\mathbf{R}\text{-\textit{order }if it is an }\mathbf{H}%
\text{-relation and }[a,  &  \ll]\text{ is }\vee\text{-complete for each }a\in
Q;\label{1.5}\\
\text{a \textit{dual }}\mathbf{R}\text{-\textit{order }if it is a dual
}\mathbf{H}\text{-relation and }[  &  \ll,a]\text{ is }\wedge\text{-complete
for each }a\in Q. \label{1.6}%
\end{align}

Following \cite{Am}, for a relation\emph{ }$\ll,$\emph{ }we define
\textit{lower }and \textit{upper complement }relations\textit{ }%
$\overrightarrow{\ll}$ and $\overleftarrow{\ll}$ by%
\begin{equation}
a\text{ }\overleftarrow{\ll}\text{ }b\text{ if }\left[  a,\ll\right]
\cap\,[a,b]=\left\{  a\right\}  ;\text{ and }a\text{ }\overrightarrow{\ll
}\text{ }b\text{ if }\left[  \ll,b\right]  \,\cap\,[a,b]=\left\{  b\right\}
\text{ for }a\leq b. \label{A1}%
\end{equation}
An element $\mathfrak{r}\in Q$ is called a $\ll$\textit{-radical} and,
respectively, $\mathfrak{p}\in Q$ is called a \textit{dual }$\ll
$-\textit{radical},\textit{ }if%
\begin{equation}
\mathbf{0}\ll\mathfrak{r}\text{ }\overleftarrow{\ll}\text{ }\mathbf{1}\text{
and, respectively, }\mathbf{0}\text{\ }\overrightarrow{\ll}\text{
}\mathfrak{p}\ll\mathbf{1}. \label{2.11}%
\end{equation}

Amitsur proved that, for an $\mathbf{R}$-order $\ll,$%
\begin{equation}
\mathfrak{r}=\vee\left[  \mathbf{0},\ll\right]  \text{ is a unique }%
\ll\text{-radical in }Q\text{ and }[\ll,\mathfrak{r}]=[\mathbf{0}%
,\mathfrak{r}]; \label{1.1}%
\end{equation}
for a dual $\mathbf{R}$-order $\ll,$ $\mathfrak{p}=\wedge\left[
\ll,\mathbf{1}\right]  $ is a unique dual $\ll$-radical in $Q$ and$\mathit{\ }%
[\mathfrak{p},\ll]=[\mathfrak{p},\mathbf{1}].$ Thus $\mathfrak{r}$ (resp.
$\mathfrak{p})$ is the largest (resp. smallest) element of $Q$ $\ll$-related
to $\mathbf{0}$ (resp. $\mathbf{1}$).

The set of $\ll$-radicals in $Q$ may be empty or have many elements. Even if
$\ll$ is an $\mathbf{H}$- or a dual $\mathbf{H}$-relation, it does not,
necessarily, has a radical. However, if it has a radical then it is unique.

Numerous natural relations in lattices are often $\mathbf{H}$- or dual
$\mathbf{H}$-relations, but seldom $\mathbf{R}$- or dual $\mathbf{R}$-orders.
Amitsur \cite{Am} introduced a procedure for construction of $\mathbf{R}$- and
dual $\mathbf{R}$-orders from $\mathbf{H}$- and dual $\mathbf{H}$-relations.
This procedure was refined in \cite{KST2}. We sketch it below. If%
\begin{equation}
x\ll y\neq x,\text{ }x\text{ is called a}\ll\text{-\textit{predecessor} of
}y\text{ and }y\text{ is a }\ll\text{-\textit{successor} of }x. \label{2.4}%
\end{equation}
A subset $G$ in $Q$ is a \textit{lower }$\ll$-set if each $x\in G\diagdown
\{\wedge G\}$ has a $\ll$-predecessor in $G$; it is an \textit{upper }$\ll
$-set if each $x\in G\diagdown\{\vee G\}$ has a $\ll$-successor in $G$.

Consider now the following relations in $Q$:%
\begin{align}
a  &  \ll^{\text{lo}}b\text{ if }[a,b]\text{ is a lower }\ll\text{-set,
}\nonumber\\
a  &  \ll^{\text{up}}b\text{ if }[a,b]\text{ is an upper }\ll\text{-set for
}a\leq b\text{ in }Q. \label{2.1}%
\end{align}
It was proved in \cite{KST2} that $\ll^{\text{lo}}$ and $\ll^{\text{up}}$ are
naturally linked to the relations $\overleftarrow{\ll}$ and
$\overrightarrow{\ll}$
\begin{equation}
\overleftarrow{\left(  \overrightarrow{\ll}\right)  }=\text{ }\ll^{\text{lo}%
}\text{ and }\overrightarrow{\left(  \overleftarrow{\ll}\right)  }=\text{ }%
\ll^{\text{up}}. \label{p11}%
\end{equation}

A chain\textit{ }$C$ in $Q$ is a linearly ordered set: $x\leq y,$ or $y\leq x$
for $x,y\in C.$ We say that a chain $C$ is a \textit{descending transfinite}
$\ll$-series from $b$ to $a$ if there is a transfinite number $\gamma$ such
that $C=\left(  x_{\lambda}\right)  _{1\leq\lambda\leq\gamma}$, where
$x_{1}=\vee C=b,$ $x_{\gamma}=\wedge C=a,$ $x_{\lambda+1}\ll x_{\lambda}$ for
all $\lambda<\gamma,$ and $x_{\beta}=\wedge_{\lambda<\beta}(x_{\lambda})$ for
limit ordinals $\beta.$

Similarly, $C$ is an \textit{ascending transfinite} $\ll$-series from $a$ to
$b$, if there is a transfinite number $\gamma$ such that $C=\left(
x_{\lambda}\right)  _{1\leq\lambda\leq\gamma}$, where $x_{1}=\wedge C=a,$
$x_{\gamma}=\vee C=b,$ $x_{\lambda}\ll x_{\lambda+1}$ for all $\lambda
<\gamma,$ and $x_{\beta}=\vee_{\lambda<\beta}(x_{\lambda})$ for limit ordinals
$\beta.$

Finally, consider two more relations generated by a relation $\ll$. We write%
\begin{align}
a  &  \ll^{\triangleleft}b\text{ if there is a descending transfinite}%
\ll\text{-series from }b\text{ to }a;\nonumber\\
a  &  \ll^{\triangleright}b\text{ if there is an ascending transfinite}%
\ll\text{-series from }a\text{ to }b. \label{2.2}%
\end{align}
Clearly,
\begin{equation}
a\ll b\text{ implies }a\ll^{\triangleleft}b\text{ and }a\ll^{\triangleright}b.
\label{2.3}%
\end{equation}

\begin{theorem}
\label{inf}\emph{(\cite{KST2}) (i) }Let $\ll$ be an $\mathbf{H}$-relation in
$Q.$ Then $\ll^{\text{\emph{up}}}$ $=$ $\ll^{\triangleright}$ \textit{is an
}$\mathbf{R}$-\textit{order}\emph{,} $\overleftarrow{\ll}$ $=$
$\overleftarrow{\ll^{\triangleright}}$ \textit{is a dual }$\mathbf{R}%
$-\textit{order and }the $\ll^{\triangleright}$-radical coincides with the
dual $\overleftarrow{\ll}$-radical$.$

An $\mathbf{H}$-relation $\ll$ is an $\mathbf{R}$-order if and only if $\ll$
$=$ $\ll^{\triangleright}.\smallskip$

\emph{(ii) }Let $\ll$ be a dual $\mathbf{H}$-relation in $Q.$ Then
$\ll^{\text{\emph{lo}}}$ $=$ $\ll^{\triangleleft}$ is a dual $\mathbf{R}%
$-order\emph{,} $\overrightarrow{\ll}$ $=$ $\overrightarrow{\ll^{\triangleleft
}}$ is an $\mathbf{R}$-order and the dual $\ll^{\triangleleft}$-radical
coincides with the $\overrightarrow{\ll}$-radical$.$

A dual $\mathbf{H}$-relation $\ll$ is a dual $\mathbf{R}$-order if and only if
$\ll$ $=$ $\ll^{\triangleleft}.$
\end{theorem}

We summarize below some results obtained in \cite{Am} and \cite{KST2} about
$\mathbf{H}$- and dual $\mathbf{H}$-relations.

\begin{theorem}
\label{T3.5}\emph{(i) }Let $\ll$ be an $\mathbf{H}$-relation. Then
$\ll^{\triangleright}$ is an $\mathbf{R}$-order\emph{, }$(\ll^{\triangleright
})^{\triangleright}=$ $\ll^{\triangleright}$ and\smallskip

$1)$ $\mathfrak{r}^{\triangleright}=\vee\lbrack\mathbf{0},\ll^{\triangleright
}]$ is the $\ll^{\triangleright}$-radical in $Q,$ i.e.$,$ $\mathbf{0}%
\ll^{\triangleright}\mathfrak{r}^{\triangleright}\,\overleftarrow{\ll}$
$\mathbf{1};$\smallskip

$2)$ $\mathfrak{r}^{\triangleright}$ is the largest $\ll^{\triangleright}%
$-successor of $\mathbf{0}$\emph{ }and the smallest $\,\overleftarrow{\ll}%
$-predecessor of $\mathbf{1}$\emph{;}$\mathbf{\smallskip}$

$3)$ for each $z\in\lbrack\mathbf{0},\mathfrak{r}^{\triangleright}],$ there is
an ascending transfinite $\ll$-series from $z$ to $\mathfrak{r}%
^{\triangleright},$ i.e.\emph{, }$z\ll^{\triangleright}\mathfrak{r}%
^{\triangleright};\smallskip$

$4)$ each $z\in Q\diagdown\lbrack\mathfrak{r}^{\triangleright},\mathbf{1}]$
has a $\ll$-successor and $\mathfrak{r}^{\triangleright}$ has no $\ll
$-successor\emph{.\smallskip}

\emph{(ii) }Let $\ll$ be a dual $\mathbf{H}$-relation. Then $\ll
^{\triangleleft}$ is a dual $\mathbf{R}$-order\emph{, }$(\ll^{\triangleleft
})^{\triangleleft}=$ $\ll^{\triangleleft}$ and\smallskip

$1)$ $\mathfrak{p}^{\triangleleft}=\wedge\lbrack\ll^{\triangleleft}%
,\mathbf{1}]$ is the dual $\ll^{\triangleleft}$-radical in $Q,$ i.e.$,$
$\mathbf{0}\,\overrightarrow{\ll}\,\mathfrak{p}^{\triangleleft}\,\ll
^{\triangleleft}\mathbf{1};$\smallskip

$2)$ $\mathfrak{p}^{\triangleleft}$ is the smallest $\ll^{\triangleleft}%
$-predecessor of $\mathbf{1}$\emph{ }and the smallest $\,\overleftarrow{\ll}%
$-successor of $\mathbf{1}$\emph{;}$\mathbf{\smallskip}$

$3)$ for each $z\in\lbrack\mathfrak{p}^{\triangleleft},\mathbf{1}],$ there is
a descending transfinite $\ll$-series from $z$ to $\mathfrak{p}^{\triangleleft
},$ i.e.\emph{, }$\mathfrak{p}^{\triangleleft}\ll^{\triangleleft}z;\smallskip$

$4)$ each $z\in Q\diagdown\lbrack\mathbf{0},\mathfrak{p}^{\triangleleft}]$ has
a $\ll$-predecessor and $\mathfrak{p}^{\triangleleft}$ has no $\ll
$-predecessor\emph{.}
\end{theorem}

A one-to-one map $g$: $Q\rightarrow Q$ is a lattice automorphism if $x\leq
y\Leftrightarrow g\left(  x\right)  \leq g\left(  y\right)  .$

\begin{theorem}
\label{T9.1}\emph{\cite{KST2} }Let $g$ be a lattice automorphism of $Q$ that
preserves a relation $\ll$ in $Q$\emph{: }$x\ll y$\emph{ }$\Leftrightarrow
$\emph{ }$g(x)\ll g(y)$\emph{.}$\smallskip$

\emph{(i) \ }If $\ll$ is an $\mathbf{H}$-relation and $\mathfrak{r}%
^{\triangleright}$ is the $\ll^{\triangleright}$-radical$,$ then
$g(\mathfrak{r}^{\triangleright})=\mathfrak{r}^{\triangleright}.\smallskip$

\emph{(ii)} If $\ll$ is a dual $\mathbf{H}$-relation and $\mathfrak{p}%
^{\triangleleft}$ is the dual $\ll^{\triangleleft}$-radical$,$ then
$g(\mathfrak{p}^{\triangleleft})=\mathfrak{p}^{\triangleleft}.$
\end{theorem}

\section{$H$-relations in Id$_{A}$ generated by properties}

\subsection{Basic definitions and constructions}

Recall that $\mathfrak{A}$ denotes the class of all C*-algebras and $P$
denotes a C*-property in $\mathfrak{A}$ (Definition \ref{D1.1}). If $A\in P,$
we say that $A$ has property $P,$ or that it is a $P$-\textit{algebra}.

For $A\in\mathfrak{A},$ Id$_{A}$ is the lattice of all closed two-sided ideals
(called just ideals) of $A$. The relation $\leq$ in Id$_{A}$ is the inclusion
of ideals $\subseteq$. For each $G\subseteq$ Id$_{A},$%
\[
\vee G=\overline{%
{\displaystyle\sum}
\{I\in G\}}\text{ \ and }\wedge G=\cap\{I\in G\}.
\]

We often use the fact that if $J\in$ Id$_{A}$ and $I\in$ Id$_{J}$ then $I\in$
Id$_{A}$. It is well known also that the sum of two closed ideals is closed:
$I+J\in$ Id$_{A}$, for $I,J\in$ Id$_{A}$, and%
\begin{equation}
(J+I)/I\approx J/(I\cap J)\text{ and }K/J\approx(K/I)/(J/I)\text{ if }I\subset
J\subset K. \label{4.01}%
\end{equation}
The proof of these properties can be found in Section 1.8 of \cite{D}.\textbf{
}Furthermore, for each $I,J\in$ Id$_{A}$, $I\cap J=IJ$\textbf{ }- the closed
linear span of $\{ab$: $a\in I,b\in J\}$ (see, for example, \cite{Murphy},
Section 3.1). It follows easily that the lattice Id$_{A}$ is distributive:
\begin{equation}
(I+J)\cap K=I\cap K+J\cap K. \label{distrib}%
\end{equation}

\begin{notation}
\label{N1}Let $P$ be a \emph{C*}-property in $\mathfrak{A}.$ For
$A\in\mathfrak{A,}$ let $\ll_{_{P}}$ be the relation in \emph{Id}$_{A}$
defined in \emph{(\ref{4.1}).\smallskip}

\emph{(i) \ }If $\ll_{_{P}}$ is an $\mathbf{R}$-order\emph{,} the $\ll_{_{P}}%
$-radical in \emph{Id}$_{A}$ is denoted by $\mathfrak{r}_{_{P}}(A).$

\ \ \ \ \ If $\ll_{_{P}}$ is an $\mathbf{H}$-relation then the $\ll_{_{P}%
}^{\triangleright}$-radical is denoted by $\mathfrak{r}_{_{P}}^{\triangleright
}(A).\smallskip$

\emph{(ii) }If $\ll_{_{P}}$ is a dual $\mathbf{R}$-order\emph{,} the dual
$\ll_{_{P}}$-radical in \emph{Id}$_{A}$ is denoted by $\mathfrak{p}_{_{P}%
}(A).$

\ \ \ \ If $\ll_{_{P}}$ is a dual $\mathbf{H}$-relation then the dual
$\ll_{_{P}}^{\triangleleft}$-radical is denoted by $\mathfrak{p}_{_{P}%
}^{\triangleleft}(A).$
\end{notation}

We write $\mathfrak{r}_{_{P}}^{\triangleright},$ $\mathfrak{p}_{_{P}%
}^{\triangleleft}$ instead of $\mathfrak{r}_{_{P}}^{\triangleright}(A),$
$\mathfrak{p}_{_{P}}^{\triangleleft}(A)$ if it is clear what algebra $A$ we
are dealing with.

The following result establishes the conditions on $P$ for $\ll_{_{P}}$ to be
an $\mathbf{H}$- or a dual $\mathbf{H}$-relation.

\begin{theorem}
\label{T4.1}\emph{(i) }A \emph{C*}-property\emph{\ }$P$ is upper stable if and
only if $\ll_{_{P}}$ is an $\mathbf{H}$-relation in \emph{Id}$_{A}$ for each
$A\in\mathfrak{A}.$ In this case the relation $\ll_{_{P}}^{\triangleright}$ is
an $\mathbf{R}$-order in $\emph{Id}_{A}\emph{\ }$and%
\begin{equation}
\{0\}\ll_{_{P}}^{\triangleright}\mathfrak{r}_{_{P}}^{\triangleright}\text{
}\overleftarrow{\ll_{_{P}}^{\triangleright}}\text{ }A\mathbf{.} \label{3.6}%
\end{equation}

\emph{(ii)\ }A \emph{C*}-property\emph{\ }$P$ is lower stable if and only
if\emph{ }$\ll_{_{P}}$ is a dual $\mathbf{H}$-relation in \emph{Id}$_{A}$ for
each $A\in\mathfrak{A}.$\emph{ }In this case the relation\emph{ }$\ll_{_{P}%
}^{\triangleleft}$ is a dual $\mathbf{R}$-order in $\emph{Id}_{A}\emph{\ }$and%
\begin{equation}
\{0\}\text{ }\overrightarrow{\ll_{_{P}}^{\triangleleft}}\text{ }%
\mathfrak{p}_{_{P}}^{\triangleleft}\ll_{_{P}}^{\triangleleft}A\mathbf{.}
\label{3.4}%
\end{equation}

\end{theorem}

\begin{proof}
Let $A\in\mathfrak{A}$ and $I\ll_{_{P}}J$ in Id$_{A},$ i.e., $J/I\in P$ (see
(\ref{4.1})).

(i) Let $I\subseteq K\in$ Id$_{A}.$ As $I\subseteq K\cap J,$ we have that
$(K\cap J)/I$ is an ideal of $J/I.$ As $J/I\in P$ and since $P$ is upper
stable, the quotient $(J/I)/((K\cap J)/I)\in P.$ By (\ref{4.01}),%
\[
(J+K)/K\approx J/(K\cap J)\approx(J/I)/((K\cap J)/I)\in P,\text{ so that
}(J+K)/K\in P.
\]
Thus $K\ll_{_{P}}(J+K).$ So, by (\ref{1.2}), $\ll_{_{P}}$ is an $\mathbf{H}$-relation.

Conversely, if $A\in P$ and $I\in$ Id$_{A}$ then $\{0\}\ll_{_{P}}A$ by
(\ref{4.1}). As $\ll_{_{P}}$ is an $\mathbf{H}$-relation in Id$_{A}$,
$I\ll_{_{P}}(A+I)=A$ by (\ref{1.2}). Thus $A/I\in P$ by (\ref{4.1}). So $P$ is
upper stable. The rest follows from Theorem \ref{T3.5}.

(ii) Let $K\in$ Id$_{A}$ and $K\subseteq J.$ As $K+I$ is an ideal of $J,$
$(K+I)/I$ is an ideal of $J/I.$ As $J/I\in P$ and since $P$ is lower stable,
$(K+I)/I\in P.$ By (\ref{4.01}), $K/(I\cap K)\approx(K+I)/I.$ Hence $K/(I\cap
K)\in P.$ Thus $I\cap K\ll_{_{P}}K.$ By (\ref{dh1}), $\ll_{_{P}}$ is a dual
$\mathbf{H}$-relation.

Conversely, let $A\in P$ and $I\in$ Id$_{A}.$ Then $\{0\}\ll_{_{P}}A$ by
(\ref{4.1}). Since $\ll_{_{P}}$ is a dual $\mathbf{H}$-relation,
$\{0\}\ll_{_{P}}I=A\cap I$ by (\ref{dh1}). Thus $I\in P.$ So $P$ is lower
stable. The rest follows from the dual of Theorem \ref{T3.5}.\bigskip
\end{proof}

Stability of many interesting C*-properties was actively studied in the theory
of C*-algebras.

In agreement with general definition a chain $\left(  I_{\lambda}\right)
_{1\leq\lambda\leq\gamma}$ of ideals is a descending transfinite $\ll_{_{P}}%
$-series from $J$ to $I,$ if $I_{1}=J,$ $I_{\gamma}=I,$ $I_{\lambda}\subseteq
I_{\mu}$ for $\mu\leq\lambda\leq\gamma,$%
\begin{equation}
I_{\lambda}/I_{\lambda+1}\in P\text{ for all }\lambda<\gamma\text{, and
}I_{\beta}=\cap_{\lambda<\beta}I_{\lambda}\text{ for all limit ordinals }%
\beta. \label{d}%
\end{equation}
It is an ascending transfinite $\ll_{_{P}}$-series of ideals from $I$ to $J$,
if $I_{1}=I,$ $I_{\gamma}=J,$ $I_{\mu}\subseteq I_{\lambda}$ for $\mu
\leq\lambda\leq\gamma,$%
\begin{equation}
I_{\lambda+1}/I_{\lambda}\in P\text{ for all }\lambda<\gamma\text{,
and}I_{\beta}=\overline{\cup_{\lambda<\beta}I_{\lambda}}\text{ for all limit
ordinals }\beta. \label{a}%
\end{equation}

Combining Theorem \ref{T3.5} and its dual (for dual $\mathbf{H}$-relations)
with Theorem \ref{T4.1} yields

\begin{corollary}
\label{C4.2n} Let $P$ be a \emph{C*}-property\emph{,} $A\in\mathfrak{A}$ and
$I\in$ \emph{Id}$_{A}$.\smallskip

\emph{(i)} If $P$ is upper stable then $(\ll_{_{P}}^{\triangleright})_{_{P}%
}^{\triangleright}=$ $\ll_{_{P}}^{\triangleright}.$ Moreover\emph{,\smallskip}

$\qquad1)$ $\mathfrak{r}_{_{P}}^{\triangleright}=\overline{\sum\{J\in\text{
\emph{Id}}_{A}\text{\emph{: }}\{0\}\ll_{_{P}}^{\triangleright}J\}}$
\emph{(}see \emph{(\ref{3.6}))}$.\mathbf{\smallskip}$

$\qquad2)$ $\mathfrak{r}_{_{P}}^{\triangleright}$ is the largest $\ll_{_{P}%
}^{\triangleright}$-successor of $\{0\}$\emph{ }and the smallest
$\,\overleftarrow{\ll_{_{P}}^{\triangleright}}$-predecessor of $A$%
\emph{.}$\mathbf{\smallskip}$

$\qquad3)$ If $I\subseteq\mathfrak{r}_{_{P}}^{\triangleright},$ there is an
ascending transfinite $\ll_{_{P}}$-series of ideals from $I$ to $\mathfrak{r}%
_{_{P}}^{\triangleright}$\emph{: }$I\ll_{_{P}}^{\triangleright}\mathfrak{r}%
_{_{P}}^{\triangleright}.\smallskip$

$\qquad4)$ If $\mathfrak{r}_{_{P}}^{\triangleright}\nsubseteq I\neq A$ then
$I$ has a $\ll_{_{P}}$-successor\emph{: }there is $J\in(I,A]$ such that
$J/I\in P.\smallskip$

\qquad$5)$ $\mathfrak{r}_{_{P}}^{\triangleright}$ has no $\ll_{_{P}}%
$-successor$.\smallskip$

$\qquad6)$ If there is an ascending transfinite $\ll_{_{P}}$-series of ideals
from $\{0\}$ to $I$ then $I\subseteq\mathfrak{r}_{_{P}}^{\triangleright
}.\smallskip$

\emph{(ii)} If $P$ is lower stable then $(\ll_{_{P}}^{\triangleleft})_{_{P}%
}^{\triangleleft}=$ $\ll_{_{P}}^{\triangleleft}.$ Moreover\emph{,\smallskip}

$\qquad1)$ $\mathfrak{p}_{_{P}}^{\triangleleft}=\cap\{J\in$ \emph{Id}$_{A}$:
$J\ll_{_{P}}^{\triangleleft}A\}$ \emph{(}see \emph{(\ref{3.4}))}%
$\mathbf{.\smallskip}$

$\qquad2)$ $\mathfrak{p}_{_{P}}^{\triangleleft}$ is the smallest $\ll_{_{P}%
}^{\triangleleft}$-predecessor of $A$ and the largest $\overrightarrow{\ll
_{_{P}}^{\triangleleft}}$-successor of $\{0\}$\emph{.}$\mathbf{\smallskip}$

$\qquad3)$ If $\mathfrak{p}_{_{P}}^{\triangleleft}\subseteq I,$ there is a
descending transfinite $\ll_{_{P}}$-series of ideals from $I$ to
$\mathfrak{p}_{_{P}}^{\triangleleft}$\emph{:} $\mathfrak{p}_{_{P}%
}^{\triangleleft}\ll_{_{P}}^{\triangleleft}I.\smallskip$

$\qquad4)$ If $\{0\}\neq I\nsubseteq\mathfrak{p}_{_{P}}^{\triangleleft}$ then
$I$ has a $\ll_{_{P}}$-predecessor\emph{: }there is $J\in\lbrack\{0\},I)$ such
that $I/J\in P.\smallskip$

$\qquad5)$ $\mathfrak{p}_{_{P}}^{\triangleleft}$ has no $\ll_{_{P}}%
$-predecessor\emph{.}$\smallskip$

$\qquad6)$ If there is a descending transfinite $\ll_{_{P}}$-series of
ideals\emph{ }from $A$ to $I,$ then $\mathfrak{p}_{_{P}}^{\triangleleft
}\subseteq I.$
\end{corollary}

\begin{definition}
\label{D3.3}We say that a \emph{C*}-property $P$ is \textbf{extension stable
}if\emph{,} for each $A\in\mathfrak{A,}$%
\begin{equation}
\text{the conditions }I\in P\text{ and }A/I\in P\text{ for some }I\in
\emph{Id}_{A},\text{ imply }A\in P. \label{3.7}%
\end{equation}

\end{definition}

\begin{proposition}
\label{P3.3}\emph{ }The relation $\ll_{_{P}}$ is transitive on each
\emph{Id}$_{A}$ if and only if $P$ is extension stable.
\end{proposition}

\begin{proof}
Let $P$ be extension stable and $I\ll_{_{P}}J\ll_{_{P}}K$ in Id$_{A}.$ Since
$J/I\in P$ and $J/I$ is isomorphic to an ideal $R$ of $K/I$ and since $K/J\in
P$ and $K/J$ is isomorphic to $(K/I)/R$, we have $K/I\in P.$ So $I\ll_{_{P}%
}K.$ Thus $\ll_{_{P}}$ is transitive. The converse is evident.
\end{proof}

\begin{example}
\label{E3.1}The C\emph{*}-property\textbf{ }$P_{un}$ of all unital C*-algebras
is extension stable.

\emph{Let }$A\in\mathfrak{A}$\emph{ and }$I\in$\emph{ Id}$_{A}.$\emph{ If }%
$I$\emph{ is unital, its identity }$e$\emph{ belongs to the centre of }%
$A.$\emph{ Indeed, }$ex,xe\in I$\emph{ for }$x\in A.$\emph{ So }%
$ex=(ex)e=e(xe)=xe.$\emph{ Then }$K=\{x-ex$\emph{: }$x\in A\}\in$\emph{
Id}$_{A}$\emph{ and }$A=I\dotplus K.$\emph{ If }$A/I$\emph{ is unital,
}$K\approx A/I$\emph{ has the identity }$f.$\emph{ So }$e+f=\mathbf{1}_{A}%
$\emph{. Thus if }$I,A/I\in P_{un}$\emph{ then }$A\in P_{un}.$
\end{example}

It will be proved later that $\ll_{_{P}}$ is an $\mathbf{R}$-order in Id$_{A}$
for each $A\in\mathfrak{A},$ if and only if $P$ is upper and extension stable
and the closure of the sums of families of $P$-ideals in $A$ are $P$-ideals.
Furthermore, $\ll_{_{P}}$ is a dual $\mathbf{R}$-order in Id$_{A}$ for each
$A\in\mathfrak{A,}$ if and only if it is lower and extension stable and, for
any family $\{I_{\lambda}$: ${I}_{\lambda}\in$ Id$_{A}$, $A/I_{\lambda}\in
P\},$ the algebra $A/\cap_{\lambda}I_{\lambda}\in P$.

We consider now some link between the radicals in algebras and in their ideals.

\begin{proposition}
\label{P3.1}Let $P$ be a \emph{C*}-property in $\mathfrak{A,}$ let
$A\in\mathfrak{A}$ and $I\in$ \emph{Id}$_{A}.\smallskip$

\emph{(i) }If $P$ is lower stable and $\mathfrak{p}_{_{P}}^{\triangleleft
}(A/I)=\{0\}$ then $I\ll_{_{P}}^{\triangleleft}A$ and $\mathfrak{p}_{_{P}%
}^{\triangleleft}(A)=\mathfrak{p}_{_{P}}^{\triangleleft}(I).\smallskip$

\emph{(ii) }Let $P$ be upper stable. Then $\mathfrak{r}_{_{P}}^{\triangleright
}(I)\subseteq\mathfrak{r}_{_{P}}^{\triangleright}(A).$\smallskip

\qquad$1)$ If $I+J=A$ for some $J\in$ \emph{Id}$_{A},$ and $J\in P$ then
$A/I\in P.\smallskip$

\qquad$2)$ If $\mathfrak{r}_{_{P}}^{\triangleright}(I)\neq\mathfrak{r}_{_{P}%
}^{\triangleright}(A)$ then $\mathfrak{r}_{_{P}}^{\triangleright}(I)\ll_{_{P}%
}J\ll_{_{P}}^{\triangleright}\mathfrak{r}_{_{P}}^{\triangleright}(A)$ for some
$\mathfrak{r}_{_{P}}^{\triangleright}(I)\neq J\in\emph{Id}_{A}.\smallskip$

\qquad$3)$ If $P$ is also lower stable and $\mathfrak{r}_{_{P}}%
^{\triangleright}(I)\ll_{_{P}}J,$ then $I\cap J=\mathfrak{r}_{_{P}%
}^{\triangleright}(I).$
\end{proposition}

\begin{proof}
(i) By Theorem \ref{T4.1}, $\ll_{_{P}}^{\triangleleft}$ is a dual $\mathbf{R}%
$-order in Id$_{A/I}\mathfrak{.}$ As $\mathfrak{p}_{_{P}}^{\triangleleft
}(A/I)=\{0\},$ there is a descending transfinite $\ll_{_{P}}$-series
$\{\widetilde{I}_{\lambda}\}_{1\leq\lambda\leq\gamma}$ of ideals of $A/I$ from
$A/I$ to $\{0\}$ by Corollary \ref{C4.2n}(ii) 3). Let $I_{\lambda}\in$
Id$_{A}$ be such that $I\subseteq I_{\lambda}$ and $\widetilde{I}_{\lambda
}\approx I_{\lambda}/I.$ As $I_{\lambda}/I_{\lambda+1}\approx\widetilde{I}%
_{\lambda}/\widetilde{I}_{\lambda+1}\in P$ by (\ref{d}), $\{I_{\lambda
}\}_{1\leq\lambda\leq\gamma}$ is a descending transfinite $\ll_{_{P}}$-series
of ideals from $A$ to $I.$ So, by (\ref{2.2}), $I\ll_{_{P}}^{\triangleleft}A.$

By (\ref{3.4}), $\mathfrak{p}_{_{P}}^{\triangleleft}(I)$ $\ll_{_{P}%
}^{\triangleleft}I.$ As $\ll_{_{P}}^{\triangleleft}$ is is transitive,
$\mathfrak{p}_{_{P}}^{\triangleleft}(I)$ $\ll_{_{P}}^{\triangleleft}A.$ By
Corollary \ref{C4.2n}(ii) 1), $\mathfrak{p}_{_{P}}^{\triangleleft}%
(A)\subseteq\mathfrak{p}_{_{P}}^{\triangleleft}(I).$ Let $\mathfrak{p}_{_{P}%
}^{\triangleleft}(A)\neq\mathfrak{p}_{_{P}}^{\triangleleft}(I).$ By Corollary
\ref{C4.2n}(ii) 5), $\mathfrak{p}_{_{P}}^{\triangleleft}(I)$ has no $\ll
_{_{P}}$-predecessor. On the other hand, as $\mathfrak{p}_{_{P}}%
^{\triangleleft}(A)\subset\mathfrak{p}_{_{P}}^{\triangleleft}(I),$ it has a
$\ll_{_{P}}$-predecessor by Corollary \ref{C4.2n}(ii) 4). This contradiction
shows that $\mathfrak{p}_{_{P}}^{\triangleleft}(A)=\mathfrak{p}_{_{P}%
}^{\triangleleft}(I).$

(ii) If $P$ is upper stable, $\ll_{_{P}}^{\triangleright}$ is an $\mathbf{R}%
$-order in Id$_{B}$ for each $B\in\mathfrak{A,}$ and $\ \{0\}\ll_{_{P}%
}^{\triangleright}\mathfrak{r}_{_{P}}^{\triangleright}(I)$ by Theorem
\ref{T4.1}. So, by Corollary \ref{C4.2n} (i) 1), $\mathfrak{r}_{_{P}%
}^{\triangleright}(I)\subseteq\mathfrak{r}_{_{P}}^{\triangleright}(A).$

1) By (\ref{4.01}), $A/I=(I+J)/I\approx J/(I\cap J).$ As $J\in P$, $J/(I\cap
J)\in P$. So $A/I\in P$.

2) If $\mathfrak{r}_{_{P}}^{\triangleright}(I)\neq\mathfrak{r}_{_{P}%
}^{\triangleright}(A),$ the result follows from Corollary \ref{C4.2n} (i).

3) Let $\mathfrak{r}_{_{P}}^{\triangleright}(I)\ll_{_{P}}J.$ Set $J_{0}=I\cap
J.$ Then $\mathfrak{r}_{_{P}}^{\triangleright}(I)\subseteq J_{0}.$ As
$J/\mathfrak{r}_{_{P}}^{\triangleright}(I)\in P$ and $P$ is lower stable, the
ideal $J_{0}/\mathfrak{r}_{_{P}}^{\triangleright}(I)$ in $J/\mathfrak{r}%
_{_{P}}^{\triangleright}(I)$ belongs to $P.$ So $\mathfrak{r}_{_{P}%
}^{\triangleright}(I)\ll_{_{P}}J_{0}.$ By Corollary \ref{C4.2n} (i) 5), this
is only possible if $\mathfrak{r}_{_{P}}^{\triangleright}(I)=J_{0}.$ Thus
$\mathfrak{r}_{_{P}}^{\triangleright}(I)=I\cap J.$
\end{proof}

\begin{corollary}
\label{C3.1}Let $P$ be a \emph{C*}-property\emph{,} $A\in\mathfrak{A}$ be
non-unital and $\widehat{A}=A+\mathbb{C}\mathbf{1}.$\smallskip

\emph{(i) }If $P$ is lower stable and contains the \emph{C*-}algebra
$\mathbb{C}\mathbf{1}$\emph{, }then $\mathfrak{p}_{_{P}}^{\triangleleft
}(A)=\mathfrak{p}_{_{P}}^{\triangleleft}(\widehat{A}).\smallskip$

\emph{(ii) }If $P$ is upper stable then $\mathfrak{r}_{_{P}}^{\triangleright
}(A)\subseteq\mathfrak{r}_{_{P}}^{\triangleright}(\widehat{A}).$ If it is also
lower stable and $A/\mathfrak{r}_{_{P}}^{\triangleright}(A)$ is
non-unital\emph{,} then $\mathfrak{r}_{_{P}}^{\triangleright}(A)=\mathfrak{r}%
_{_{P}}^{\triangleright}(\widehat{A}).$\smallskip

\emph{(iii) }Let $P$ be lower\emph{,} upper and extension stable and let
$\mathbb{C}\mathbf{1}\in P.$ If $A/\mathfrak{r}_{_{P}}^{\triangleright}(A)$ is
unital and the class $e+\mathfrak{r}_{_{P}}^{\triangleright}(A),$ $e\in A,$ is
the identity in $A/\mathfrak{r}_{_{P}}^{\triangleright}(A),$ then
$\mathfrak{r}_{_{P}}^{\triangleright}(\widehat{A})=\mathbb{C}(\mathbf{1}%
-e)+\mathfrak{r}_{_{P}}^{\triangleright}(A)\mathbf{.}$
\end{corollary}

\begin{proof}
Replacing $I$ by $A$ and $A$ by $\widehat{A}$ in Proposition \ref{P3.1}, we
get the proof of (i).

(ii) By Proposition \ref{P3.1}(ii), $\mathfrak{r}_{_{P}}^{\triangleright
}(A)\subseteq\mathfrak{r}_{_{P}}^{\triangleright}(\widehat{A})$ and, if
$\mathfrak{r}_{_{P}}^{\triangleright}(A)\neq\mathfrak{r}_{_{P}}%
^{\triangleright}(\widehat{A}),$ then%
\[
\mathfrak{r}_{_{P}}^{\triangleright}(A)\ll_{_{P}}J\ll_{_{P}}^{\triangleright
}\mathfrak{r}_{_{P}}^{\triangleright}(\widehat{A})\text{ and }\mathfrak{r}%
_{_{P}}^{\triangleright}(A)=A\cap J\text{ for some }J\in\text{ Id}%
_{\widehat{A}}\text{ }\backslash\text{ Id}_{A}.
\]
Hence, as $\dim(\widehat{A}/A)=1,$ we have $\dim(J/\mathfrak{r}_{_{P}%
}^{\triangleright}(A))=1.$ So $J=\mathbb{C}(\mathbf{1}-e)\dotplus
\mathfrak{r}_{_{P}}^{\triangleright}(A)$ for some $e\in A.$ Hence
($\mathbf{1}-e)b=b-eb\in A\cap J=\mathfrak{r}_{_{P}}^{\triangleright}(A)$ for
all $b\in A.$ Thus $b=eb+s(b)$ for some $s(b)\in\mathfrak{r}_{_{P}%
}^{\triangleright}(A),$ whence the class $e+\mathfrak{r}_{_{P}}%
^{\triangleright}(A)$ in $A/\mathfrak{r}_{_{P}}^{\triangleright}(A)$ is the
identity. So if $A/\mathfrak{r}_{_{P}}^{\triangleright}(A)$ is non-unital then
$\mathfrak{r}_{_{P}}^{\triangleright}(A)=\mathfrak{r}_{_{P}}^{\triangleright
}(\widehat{A})$.

(iii) As $e+\mathfrak{r}_{_{P}}^{\triangleright}(A)$ is the identity in
$A/\mathfrak{r}_{_{P}}^{\triangleright}(A)$, then $J=\mathbb{C}(\mathbf{1}%
-e)\dotplus\mathfrak{r}_{_{P}}^{\triangleright}(A)$ is an ideal in
$\widehat{A}.$ If $\mathbb{C}\mathbf{1}\in P$ then $\mathfrak{r}_{_{P}%
}^{\triangleright}(A)\ll_{_{P}}J,$ as $J/\mathfrak{r}_{_{P}}^{\triangleright
}(A)\approx\mathbb{C}\mathbf{1.}$ Suppose that $\mathfrak{r}_{_{P}%
}^{\triangleright}(\widehat{A})\neq J.$ Then, by Corollary \ref{C4.2n}(i),
$J\ll_{_{P}}J_{1}$ for some $J\neq J_{1}\in$ Id$_{\widehat{A}}.$ Since $P$ is
extension stable, the relation $\ll_{_{P}}$ is transitive. So $\mathfrak{r}%
_{_{P}}^{\triangleright}(A)\ll_{_{P}}J_{1}.$ Then, by Proposition \ref{P3.1}
(ii) 3), $A\cap J_{1}=\mathfrak{r}_{_{P}}^{\triangleright}(A).$ As in 1), we
get $J_{1}=\mathbb{C}(\mathbf{1}-e_{1})\dotplus\mathfrak{r}_{_{P}%
}^{\triangleright}(A)$ for some $e_{1}\in A.$ Since $J\subseteq J_{1},$ we
have $\mathbf{1}-e=\lambda(\mathbf{1}-e_{1})+a$ for some $\lambda\in
\mathbb{C}$ and $a\in\mathfrak{r}_{_{P}}^{\triangleright}(A).$ Then
$\lambda=1$ and $e_{1}=e+a.$ So $J_{1}=J.$ This contradiction shows that
$\mathfrak{r}_{_{P}}^{\triangleright}(\widehat{A})=J=\mathbb{C}(\mathbf{1}%
-e)\dotplus\mathfrak{r}_{_{P}}^{\triangleright}(A).\bigskip$
\end{proof}

Each automorphism $\phi$ of $A$ generates a lattice automorphism of Id$_{A}$.
Let $P$ be a C*-property in $\mathfrak{A}$ and $I\ll_{_{P}}J$ in Id$_{A}$:
$J/I\in P.$ Then the map $\widehat{\phi}$: $J/I\rightarrow\phi(J)/\phi(I)$
defined by $\widehat{\phi}(x+I)=\phi(x)+\phi(I)$ for $x\in J,$ is an
isomorphism. Thus $J/I\approx\phi(J)/\phi(I),$ so that, by (\ref{1.7}),%
\begin{equation}
\phi(J)/\phi(I)\in P,\text{ i.e., }\phi(I)\ll_{_{P}}\phi(J). \label{3,3}%
\end{equation}
Similarly, as $\phi^{-1}$ is an automorphism of $A,$ $\phi(I)\ll_{_{P}}%
\phi(J)$ implies $I\ll_{_{P}}J.$ Thus $\phi$ preserves $\ll_{_{P}}.$ Hence
Theorem \ref{T9.1} yields

\begin{corollary}
\label{C9.1}Let $P$ be a \emph{C*}-property in $\mathfrak{A}$ and let
$\ll_{_{P}}$ be the corresponding relation on \emph{Id}$_{A}.\smallskip$

If $\ll_{_{P}}$ is an $\mathbf{H}$-relation then $\phi(\mathfrak{r}_{_{P}%
}^{\triangleright})=\mathfrak{r}_{_{P}}^{\triangleright}$ for all
automorphisms $\phi$ of $A.\smallskip$

If $\ll_{_{P}}$ is a dual $\mathbf{H}$-relation then $\phi(\mathfrak{p}_{_{P}%
}^{\triangleleft})=\mathfrak{p}_{_{P}}^{\triangleleft}$ for all automorphisms
$\phi$ of $A.$
\end{corollary}

\subsection{Some $\ll_{_{P}}^{\triangleright}$- and dual $\ll_{_{P}%
}^{\triangleleft}$-radicals in commutative C*-algebras.}

Here we will consider the case of commutative C*-algebras $A.$ If $A$ is
unital, $A=C(X)$ for a compact $X$ and there is a one-to-one correspondence
between closed subsets $Y$ of $X$ and the ideals $I_{Y}=\{f\in A$:\emph{
}$f|_{Y}=0\}$ in Id$_{A}.$ For each $f\in A,$ set $p(f)=f|_{Y}.$ Then $p$:
$A\mapsto C(Y)$ and $\ker p=I_{Y}.$ By Urysohn's Theorem, for each $g\in
C(Y),$ there is $f\in A$ such that $g=p(f)=f|_{Y}.$ Thus
\begin{equation}
C(Y)\approx C(X)/I_{Y}. \label{3.0}%
\end{equation}

If $A$ is non-unital, then $A\approx C_{0}(X)=\{f(x)\in C(X)$: $f(x_{0})=0\}$
for some $x_{0}$ in a compact $X$. There is a one-to-one correspondence
between closed subsets $Y$ of $X$ containing $x_{0}$ and the ideals
\[
I_{Y}=\{f\in C(X)\text{: }f|_{Y}=0\}=\{f\in A\text{: }f|_{Y\backslash
\{x_{0}\}}=0\}\text{ in Id}_{A}.
\]

Denote by $G_{X}$ the group of all continuous automorphisms of a compact $X$.

\begin{corollary}
\label{C3.8}Let $A=C(X).$

\emph{(i) }If $I_{Y}$ is $\mathfrak{r}_{_{P}}^{\triangleright}$\emph{, }or
$\mathfrak{p}_{_{P}}^{\triangleleft}$ in \emph{Id}$_{A}$ for some
\emph{C*}-property $P,$ then $Y$ is invariant for all $g\in G_{X}.\smallskip$

\emph{(ii)} Let $X$ have no closed $G_{X}$-invariant subsets$.$ Then\emph{,}
for any upper \emph{(}resp.\emph{, }lower\emph{)} stable \emph{C*}-property
$P,$ the $\ll_{_{P}}^{\triangleright}$-radical \emph{(}resp.\emph{, }the dual
$\ll_{_{P}}^{\triangleleft}$-radical\emph{)} in $A$ is either $\{0\}$ or $A.$
\end{corollary}

\begin{proof}
(i) Each $g\in G_{X}$ generates an automorphism $\theta_{g}$ of $A$:
$\theta_{g}f(x)=f(gx).$ If $I_{Y}$ is a radical then, by Corollary \ref{C9.1},
$f\in I_{Y}$ implies $\theta_{g}f\in I_{Y}.$ So $f(gx)=0$ for all $x\in Y$ and
$g\in G_{X}.$ If $gx\notin Y$ for some $x\in Y$ and $g\in G_{X},$ then there
exists $f\in I_{Y}$ such that $f(gx)\neq0.$ This contradiction proves (i).
Part (ii) follows from (i).\bigskip
\end{proof}

Recall that a closed subset $E$ of a topological space is \textit{perfect }if
it has no isolated points:%
\begin{equation}
E\cap U_{y}\neq\{y\}\text{ for each }y\in E\text{ and each open neighbourhood
}U_{y}\text{ of }y. \label{3.10}%
\end{equation}
Let $\mathcal{E}_{X}$ be the set of all perfect subsets of $X.$ If $E\subset
X$ is not closed but (\ref{3.10}) holds, then $\overline{E}$ is perfect. Thus
if $\mathcal{E}_{X}\neq\varnothing$ then the set $E_{X}=\overline
{\cup\{E\text{: }E\in\mathcal{E}_{X}\}}$ is the largest perfect subset of $X,$
i.e., it contains all perfect subsets of $X.$

If $N$ is the set of all isolated points in $X$ and $\overline{N}$ its
closure, it is easy to show that $X\diagdown\overline{N}\subseteq E_{X}.$

\begin{theorem}
\label{T3.9}Let $A=C(X)$ and a \emph{C*}-property $P$ contain a
finite-dimensional commutative\textbf{ }\emph{C*-}algebra.\smallskip

\emph{(i) }If $P$ is lower stable then $\mathfrak{p}_{_{P}}^{\triangleleft
}=\{0\}$ in $A.\smallskip$

\emph{(ii) }Let $P$ be upper stable. Then%
\begin{equation}
I_{E_{X}}\subseteq\mathfrak{r}_{_{P}}^{\triangleright}. \label{3.11}%
\end{equation}

If $E_{X}=\varnothing$ then $\mathfrak{r}_{_{P}}^{\triangleright}=A.$

If\emph{ }$E_{X}\neq\varnothing$ and $P$ contains no infinite dimensional
commutative algebras\emph{, }then $I_{E_{X}}=\mathfrak{r}_{_{P}}%
^{\triangleright}\neq A.$
\end{theorem}

\begin{proof}
If $B$ is a finite-dimensional, commutative C*-algebra, it is the direct sum
of one-dimensional ideals. If $P$ is lower, or upper stable, it contains the
one-dimensional C*-algebra $\mathbb{C}\mathbf{1}$.

(i) By Theorem \ref{T4.1}, $\ll_{_{P}}$ is a dual $\mathbf{H}$-relation in
Id$_{A}.$ For $x\in X,$ the ideal $I_{\{x\}}$ is maximal and $\dim
A/I_{\{x\}}=1.$ Hence $A/I_{\{x\}}\in P,$ so that $I_{\{x\}}\ll_{_{P}}A.$
Then, by Corollary \ref{C4.2n}(ii),
\[
\mathfrak{p}_{_{P}}^{\triangleleft}=\cap\{I\in\text{ Id}_{A}\text{: }%
I\ll_{_{P}}^{\triangleleft}A\}\overset{(\ref{2.3})}{\subseteq}\cap\{I\in\text{
Id}_{A}\text{: }I\ll_{_{P}}A\}\subseteq\cap_{x\in X}I_{\{x\}}=\{0\}
\]
which completes the proof of (i).

(ii) Let $J=I_{R}\in$ Id$_{A}$ for a closed subset $R$ in $X.$ By Theorem
\ref{T4.1}, $\ll_{_{P}}$ is an $\mathbf{H}$-relation.

Let $E_{X}\neq X$ and $I_{R}\subsetneqq I_{E_{X}}.$ Then $E_{X}\subsetneqq R.$
Thus $R\ $is not perfect. So, by (\ref{3.10}), there is $y\in R\diagdown
E_{X}$ and an open neighbourhood $U_{y}$ of $y$ in $X$ such that $R\cap
U_{y}=\{y\}.$

As $R\diagdown\{y\}=R\cap(X\diagdown U_{y}),$ it is closed. So $I_{R\diagdown
\{y\}}\in$ Id$_{A},$ $I_{R}\subseteq I_{R\diagdown\{y\}}$ and there is $g\in
C(X)$ such that $g|_{R\diagdown\{y\}}=0$ and $g(y)=1.$ Hence $g\in
I_{R\diagdown\{y\}},$ $g\notin I_{R}$ and $f-f(y)g\in I_{R}$ for each $f\in
I_{R\diagdown\{y\}}.$ Thus $\dim(I_{R\diagdown\{y\}}/I_{R})=1.$ Since $P$
contains one-dimensional C*-algebras, $I_{R\diagdown\{y\}}/I_{R}\in P.$ So
$J=I_{R}\ll_{_{P}}I_{R\diagdown\{y\}},$ i.e., each ideal $J\neq I_{E_{X}}$ in
the segment $[\{0\},I_{E_{X}}]$ in Id$_{A}$ has a $\ll_{_{P}}$-successor (see
(\ref{2.4})). Hence $[\{0\},I_{E_{X}}]$ is an upper $\ll_{_{P}}$-set. So
$\{0\}\ll_{_{P}}^{\text{up}}I_{E_{X}}$ (see (\ref{2.1})). By Theorem \ref{inf}
(i), $\{0\}\ll_{_{P}}^{\triangleright}I_{E_{X}}.$ Hence, by Corollary
\ref{C4.2n} (i) 1), $I_{E_{X}}\subseteq\mathfrak{r}_{_{P}}^{\triangleright},$
so that (\ref{3.11}) holds.

If $E_{X}=X$ then $I_{E_{X}}=\{0\}\subseteq\mathfrak{r}_{_{P}}^{\triangleright
}.$ So (\ref{3.11}) holds for all commutative C*-algebras.

If $E_{X}=\varnothing$ then $I_{E_{X}}=A$ so that $\mathfrak{r}_{_{P}%
}^{\triangleright}=A$ by (\ref{3.11}).

Let $E_{X}\neq\varnothing$. Then $I_{E_{X}}\neq A.$ If $I_{E_{X}}\subsetneqq
I_{R}\in$ Id$_{A}$ then $R\subsetneqq E_{X}.$ Let $y\in E_{X}\diagdown R$ and
$U_{y}$ be an open neighbourhood of $y$ such that $R\cap\overline{U_{y}%
}=\varnothing.$ As $E_{X}$ is perfect, we have from (\ref{3.10}) that
$E_{X}\cap U_{y}\neq\{y\}$ and, moreover, that card($E_{X}\cap U_{y})=\infty.$
So card($E_{X}\cap\overline{U_{y}})=\infty.$

Fix $n\in\mathbb{N}$ and choose some $\{z_{i}\}_{i=1}^{n}$ in $E_{X}%
\cap\overline{U_{y}}.$ As $R\cap(E_{X}\cap\overline{U_{y}})=\varnothing,$
there are $g_{i}\in C(X)$ such that $g_{i}(z_{i})=1,$ $g_{i}(z_{j})=0$ for
$i\neq j,$ and $g_{i}|_{R}=0.$ Hence all $g_{i}\in I_{R}$ are linearly
independent and $g_{i}\notin I_{E_{X}}$. So $\dim(I_{R}/I_{E_{X}})\geq n.$ As
$n$ is arbitrary, $\dim(I_{R}/I_{E_{X}})=\infty.$ If $P$ contains no infinite
dimensional commutative algebras, $I_{E_{X}}\not \ll _{_{P}}I_{R}.$ Thus%
\begin{equation}
I_{E_{X}}\not \ll _{_{P}}J\text{ for each }I_{E_{X}}\subsetneqq J\in\text{
Id}_{A}. \label{3.12}%
\end{equation}

If $I_{E_{X}}\neq\mathfrak{r}_{_{P_{_{\text{fin}}}}}^{\triangleright},$ it
follows from Corollary \ref{C4.2n} (i) 4) and (\ref{3.11}) that $I_{E_{X}}$
has a $\ll_{_{P}}$-successor which contradicts (\ref{3.12}). Hence $I_{E_{X}%
}=\mathfrak{r}_{_{P}}^{\triangleright}\neq A.\bigskip$
\end{proof}

\section{Basic closure operations on C*-properties}

An important part of the Amitzur approach to radicals in lattices was the
possibility to improve an $\mathbf{H}$-relation to obtain an $\mathbf{R}%
$-relation, or to improve a dual $\mathbf{H}$-relation to obtain a dual
$\mathbf{R}$-relation. We consider corresponding C*-algebraic constructions
that allow one to improve C*-properties. These constructions can be regarded
as closure operation on classes of C*-properties. Recall that a map $f$ on a
partially ordered set $(G,\leq)$ is called a \textit{closure operator} (see
Definition I.3.26 \cite{G}) if
\begin{equation}
1)\text{ }x\leq f(x)=f(f(x)),\text{ and }2)\text{ }x\leq y\text{ implies
}f(x)\leq f(y)\text{ for }x,y\in G. \label{2.6}%
\end{equation}

\begin{lemma}
\label{L2.1}Let $f$ be a closure operator on $(G,\leq).$\textbf{ }If\textbf{
}$x\leq y\leq f(x)$\textbf{ }for $x,y\in G,$\textbf{ }then\textbf{
}$f(x)=f(y).$
\end{lemma}

\begin{proof}
As\textbf{ }$x\leq y\leq f(x),$\textbf{ }we have $f(x)\leq f(y)\leq
f(f(x))=f(x)$\textbf{ }by (\ref{2.6}).\textbf{\ }So\textbf{ }%
$f(x)=f(y).\bigskip$
\end{proof}

We study several important closure operators; more closure operators can be
obtained as their compositions.

\subsection{The closure operator $P\mapsto\mathsf{G}P$ on upper stable
properties}

The first closure operator $P\mapsto\mathsf{G}P$ we deal reminds the
transition from the class of CCR-algebras to the class of GCR-algebras.

Recall that a C*-algebra $A$ is called a $CCR$\textit{-algebra} if its
irreducible representations map $A$ to algebras of compact operators, and a
$GCR$\textit{-algebra} if the images of its irreducible representations
contain non-zero compact operators. The "opposite" class - $NGCR$%
\textit{-algebras} - consists of algebras that have no $GCR$-ideals. These
important classes of C*-algebras were intensively studied and can be
characterized by many other conditions. In paricular, $GCR$-algebras can be
charactrized as C*-algebras whose non-zero quotients contain $CCR$-ideals.

In this section we extend this construction to all C*-properties. Namely, for
each C*-property $P$ in $\mathfrak{A},$ we construct a larger C*-property
$\mathsf{G}P$ and complementary C*-property $\mathsf{NG}P$ as follows.

\begin{definition}
\label{D3.1}Let $P$ be a \emph{C*}-property in $\mathfrak{A}.$ We call a
\emph{C}$^{\ast}$-algebra $A\smallskip$

\emph{(i)\ \ }a $\mathsf{G}P$\textbf{-algebra} \emph{(}\textbf{generalized
}$P$\textbf{-algebra}$)$ if either $A=\{0\},$ or each non-zero quotient $A/I$
has a non-zero $P$-ideal $\smallskip$

\emph{(ii) }an $\mathsf{NG}P$\textbf{-algebra} if it does not have non-zero
$P$-ideals$.$
\end{definition}

Clearly, the classes $\mathsf{G}P$ of all $\mathsf{G}P$-algebras and
$\mathsf{NG}P$ of all $\mathsf{NG}P$-algebras are also C*-properties, and
these terms match the terms of $GCR$- and $NGCR$-algebras for $P=CR=CCR.$

\begin{lemma}
\label{L3.6}For each \emph{C*}-property $P,$ the \emph{C*}-property
$\mathsf{G}P$ is upper stable, while $\mathsf{NG}P$ is lower stable.
\end{lemma}

\begin{proof}
Let $A\in$\textsf{G}$P$ and $I\in$ Id$_{A}.$ Since each non-zero quotient of
$A/I$ is isomorphic to a non-zero quotient of $A$ and since all non-zero
quotients of $A$ have non-zero $P$-ideals, $A/I$ is a \textsf{G}$P$-algebra.
Thus the C*-property \textsf{G}$P$ is upper stable (Definition \ref{D4.0}).

Since an ideal of an ideal of $A$ is an ideal of $A$, all ideals of an
$\mathsf{NG}P$-algebra are $\mathsf{NG}P$-algebras. So $\mathsf{NG}P$ is a
lower stable C*-property.
\end{proof}

\begin{example}
\label{exG} \emph{1) If $P=CCR$, then \textsf{G}$P$ = $GCR$}.\smallskip

\emph{2) }Let $P=C$ consist of all \emph{C*-}algebras isomorphic to $C(H)$ for
different $H$. Then \emph{\textsf{G}}$C=Sc$ is the \emph{C*}-property of all
scattered C*-algebras.\emph{ (A Banach algebra is called \textbf{scattered} if
spectra of all its elements are countable or finite}.)

\emph{Indeed, it is known (see \cite[Section 8]{ST}, where there is a long
list of conditions equivalent to $A\in Sc$) that each $A\in Sc$ contains an
ideal isomorphic to $C(H)$. Also if} $a\in A$ \emph{and }$p$\emph{:} $A\mapsto
A/I$ \emph{for} $I\in$ \emph{Id}$_{A},$ \emph{then Sp(}$p(a))\subseteq$
\emph{Sp}$(a).$ \emph{So }$A/I$ \emph{is} \emph{scattered if }$A$ \emph{is
scattered. Thus $Sc\subseteq\mathsf{G}C$}.

\emph{Conversely if $A\in\mathsf{G}C$ then each quotient $B=A/J$ of $A$
contains an ideal $I\approx C(H)$. So $I$ has a minimal projection $p.$ As
$pBp=p(pBp)p\subset pIp$, we have $\dim(pBp)=1$. Thus $p$ is a minimal
projection in $B$. So each quotient of $A$ has a minimal projection; this
condition is equivalent to $A\in Sc$ \cite[Section 8]{ST}. Thus $\mathsf{G}%
C=Sc$}.\smallskip

\emph{3) }Let $P=Comm$ be the \emph{C*}-property of all commutative algebras.
Then $GP=P$.

\emph{Note firstly that $P$ is upper stable: }$A/I\in P$\emph{, if} $A\in
P.$\emph{ So $P\subseteq\mathsf{G}P$ by Definition \ref{D3.1}. Thus we only
have to prove the converse inclusion}.

\emph{Let $A\in\mathsf{G}P$. If $\pi$ is an irreducible representation of $A$,
then $\pi(A)\approx A/\ker(\pi)$ has a commutative ideal $J\neq\{0\}$. By
Lemma 2.11.3(i) \cite{D}, }$H_{\pi}$ \emph{has no }$J$\emph{-invariant
subspaces. So the identity representation of }$J$ \emph{on }$H_{\pi}$\emph{ is
irreducible. As }$J$\emph{ is commutative,} $\dim H_{\pi}=1.$ \emph{Hence}
\emph{$\pi(A)$ is commutative}.

\emph{Suppose that $xy\neq yx$, for some $x,y\in A$. Then there is an
irreducible representation $\pi$ with $\pi(xy-yx)\neq0$. Thus $\pi(x)$ and
$\pi(y)$ do not commute, a contradiction}. \emph{Thus }$A\in P.$\emph{ So
$\mathsf{G}P=P$}.\smallskip

\emph{4) }Let $P=\mathbb{C}\mathbf{1}$ be the \emph{C*}-property of all
one-dimensional algebras. Then $\mathsf{G}P$ is the \emph{C*}-property of all
commutative algebras with dispersed space of maximal ideals.\emph{ (A
topological space is called \textbf{dispersed} if it has no perfect subsets).}

\emph{Indeed, since $P\subset Comm$, we have $GP\subset\mathsf{G}Comm=Comm$ by
Proposition \ref{L3.1} and 3). By Definition \ref{D3.1},} \emph{$A=C(X)$
belongs to $\mathsf{G}P$ if and only if the quotient }$C(X)/I_{Y}$ \emph{has a
one-dimensional ideal for each closed subset }$Y$\emph{ of }$X$\emph{. By
(\ref{3.0}),} $C(X)/I_{Y}\approx C(Y).$ \emph{Thus }$A\in$ \emph{$\mathsf{G}P$
if and only if }$C(Y)$ \emph{has a one-dimensional ideal for each closed
subset }$Y$ \emph{of }$X,$\emph{ i.e.,} \emph{each }$Y$ \emph{has an isolated
point. So }$X$ \emph{has no perfect subsets, i.e., }$X$\emph{ is dispersed.}
\end{example}

As in (\ref{4.1}), we consider the relations $\ll_{_{\mathsf{G}P}}$ and
$\ll_{_{\mathsf{NG}P}}$ in Id$_{A}$: $I\ll_{_{\mathsf{G}P}}J$ if $J/I\in
$\textsf{G}$P,$ and $I\ll_{_{\mathsf{NG}P}}J$ if $J/I\in\mathsf{NG}P.$ It
follows from Definition \ref{D3.1} that
\begin{equation}
I\ll_{_{\mathsf{G}P}}J\Longleftrightarrow\text{ for each }K\in\lbrack
I,J),\text{ there is }L\in(K,J]\text{ such that }L/K\in P. \label{3.21}%
\end{equation}

\begin{theorem}
\label{T3.0}Let $P$ be an upper stable \emph{C*}-property in $\mathfrak{A}$
and $A\in\mathfrak{A}$\emph{. Then}$\smallskip$

\emph{(i) \ \ }$\ll_{_{\mathsf{G}P}}=$ $\ll_{_{P}}^{\triangleright}$ in
\emph{Id}$_{A},$ so that $\ll_{_{\mathsf{G}P}}$ is an $\mathbf{R}$-order in
\emph{Id}$_{A}$ $($see \emph{(\ref{1.5}));}\smallskip

\emph{(ii) \ }$\ll_{_{\mathsf{NG}P}}=$ $\overleftarrow{\ll_{_{P}}}$ in
\emph{Id}$_{A},$ so that $\ll_{_{\mathsf{NG}P}}$ is a dual $\mathbf{R}$-order
in \emph{Id}$_{A}$ $($see \emph{(\ref{1.6}));\smallskip}

\emph{(iii) }$\mathfrak{r}_{_{P}}^{\triangleright}=\mathfrak{r}_{_{\mathsf{G}%
P}}=\mathfrak{p}_{_{\mathsf{NG}P}}.$
\end{theorem}

\begin{proof}
As $P$ is upper stable, $\ll_{_{P}}$ is an $\mathbf{H}$-relation in Id$_{A}$
by Theorem \ref{T4.1}.

(i) Let $I\ll_{_{\mathsf{G}P}}J,$ $I\neq J,$ in Id$_{A}.$ It follows from
(\ref{2.4}) and (\ref{3.21}) that $I\ll_{_{\mathsf{G}P}}J$ if and only if each
$K\in\lbrack I,J)$ has a $\ll_{_{P}}$-successor$,$ i.e., $[I,J]$ is an upper
$\ll_{_{P}}$-set. Thus we have from (\ref{2.1}) that $I\ll_{_{\mathsf{G}P}%
}J\Longleftrightarrow I\ll_{_{P}}^{\text{up}}J$. So $\ll_{_{\mathsf{G}P}}=$
$\ll_{_{P}}^{\text{up}}.$ As $\ll_{_{P}}$ is an $\mathbf{H}$-relation,
$\ll_{_{P}}^{\text{up}}$ $=$ $\ll_{_{P}}^{\triangleright}$ by Theorem
\ref{inf}(i). Thus $\ll_{_{\mathsf{G}P}}=$ $\ll_{_{P}}^{\triangleright}$ is an
$\mathbf{R}$-order in Id$_{A}$ by Theorem \ref{T4.1}(i).\smallskip

(ii) Let $I\ll_{_{\mathsf{NG}P}}J$ and $I\neq J.$ Then $J/I$ is an
$\mathsf{NG}P$-algebra, i.e., it contains no non-zero $P$-ideals. Hence, for
each $K\in(I,J],$ $K/I$ is not a $P$-ideal in $J/I$, i.e., $I\not \ll _{_{P}%
}K.$ Thus $I\ll_{_{\mathsf{NG}P}}J$ if and only if $[I,\ll_{_{P}}]\cap\lbrack
I,J]=I.$ So, by (\ref{A1}), $I\ll_{_{\mathsf{NG}P}}J\Longleftrightarrow I$
$\overleftarrow{\ll_{_{P}}}$ $J.$ Hence $\ll_{_{\mathsf{NG}P}}$ $=$
$\overleftarrow{\ll_{_{P}}}$ is a dual $\mathbf{R}$-order in Id$_{A}$ by
Theorem \ref{inf}(i).\smallskip

(iii) Since $\ll_{_{P}}^{\triangleright}=$ $\ll_{_{\mathsf{G}P}},$ we have
$\mathfrak{r}_{_{P}}^{\triangleright}=\mathfrak{r}_{_{\mathsf{G}P}}.$ As
$\ll_{_{\mathsf{NG}P}}$ $=$ $\overleftarrow{\ll_{_{P}}},$ $\mathfrak{p}%
_{_{\mathsf{NG}P}}$ coincides with the dual $\overleftarrow{\ll_{_{P}}}%
$-radical which, in turn, coincides with $\mathfrak{r}_{_{P}}^{\triangleright
}$ by Theorem \ref{inf}.\bigskip
\end{proof}

The following theorem generalizes various results that hold for $GCR$- and
$NGCR$-algebras.

\begin{theorem}
\label{T3.1}Let $P$ be an upper stable \emph{C*}-property in $\mathfrak{A}$
and $A\in\mathfrak{A}$. Then\smallskip

\emph{(i)} $\mathfrak{r}_{_{P}}^{\triangleright}$ is the largest $\mathsf{G}%
P$-ideal of $A$ and the smallest ideal with $\mathsf{NG}P$-quotient$.$ There
is an ascending transfinite $\ll_{_{P}}$-series of ideals from $\{0\}$ to
$\mathfrak{r}_{_{P}}^{\triangleright}$ $($see $(\ref{a}))$.\smallskip

\emph{(ii) }A \emph{C*}-algebra $A$ is a $\mathsf{G}P$-algebra if and only if
$A=\mathfrak{r}_{_{P}}^{\triangleright};$ a \emph{C*}-algebra $A$ is an
$\mathsf{NG}P$-algebra if and only if $\mathfrak{r}_{_{P}}^{\triangleright
}=\{0\}.$
\end{theorem}

\begin{proof}
(i) By (\ref{3.6}), $\{0\}\ll_{_{P}}^{\triangleright}\mathfrak{r}_{_{P}%
}^{\triangleright}.$ As $\ll_{_{\mathsf{G}P}}=$ $\ll_{_{P}}^{\triangleright}$
by Theorem \ref{T3.0}, we have $\{0\}\ll_{_{\mathsf{G}P}}\mathfrak{r}_{_{P}%
}^{\triangleright}.$ So $\mathfrak{r}_{_{P}}^{\triangleright}$ is a
\textsf{G}$P$-algebra. If $I\in$ Id$_{A}$ is a \textsf{G}$P$-ideal then
$\{0\}\ll_{_{\mathsf{G}P}}I.$ As $\ll_{_{\mathsf{G}P}}=$ $\ll_{_{P}%
}^{\triangleright},$ we have $\{0\}\ll_{_{P}}^{\triangleright}I.$ Hence
$I\subseteq\mathfrak{r}_{_{P}}^{\triangleright}$ by Corollary \ref{C4.2n}$.$
Thus $\mathfrak{r}_{_{P}}^{\triangleright}$ contains all \textsf{G}$P$-ideals
of $A.$

By (\ref{3.6}), $\mathfrak{r}_{_{P}}^{\triangleright}$ $\overleftarrow{\ll
_{_{P}}^{\triangleright}}$ $A.$ By Theorem \ref{inf} (i), $\overleftarrow{\ll
_{_{P}}^{\triangleright}}$ $=$ $\overleftarrow{\ll_{_{P}}}.$ So $\mathfrak{r}%
_{_{P}}^{\triangleright}$ $\overleftarrow{\ll_{_{P}}}$ $A.$ Hence, by Theorem
\ref{T3.0} (ii), $\mathfrak{r}_{_{P}}^{\triangleright}$ $\ll_{_{\mathsf{NG}P}%
}$ $A.$ Thus $A/\mathfrak{r}_{_{P}}^{\triangleright}$ is an $\mathsf{NG}P$-algebra.

Let $A/I$ be an $\mathsf{NG}P$-algebra for some $I\in$ Id$_{A}$. If
$\mathfrak{r}_{_{P}}^{\triangleright}\nsubseteq I$ then, by Corollary
\ref{C4.2n}, $I\ll_{_{P}}J\neq I$ for some $J\in$ Id$_{A}$. So $\{0\}\neq J/I$
is a $P$-ideal of $A/I$ -- a contradiction. So $\mathfrak{r}_{_{P}%
}^{\triangleright}\subseteq I.$

By Corollary \ref{C4.2n}, there exists an ascending transfinite $\ll_{_{P}}%
$-series of ideals from $\{0\}$ to $\mathfrak{r}_{_{P}}^{\triangleright}$.

Part (ii) follows from (i).\bigskip
\end{proof}

We show now that the C*-properties $\mathsf{G}P$ and $\mathsf{NG}P$ are also
extension stable.

\begin{proposition}
\label{P3.3P} If $P$ is upper stable then the \emph{C*}-properties
$\mathsf{G}P$ and $\mathsf{NG}P$ are extension stable \emph{(}Definition
\emph{\ref{D3.3}).}
\end{proposition}

\begin{proof}
Let $A\in\mathfrak{A}$ and let $I\in$ Id$_{A}$ and $A/I$ be \textsf{G}%
$P$-algebras. If $\mathfrak{r}_{_{P}}^{\triangleright}\neq A$ then, by Theorem
\ref{T3.1}(i), $I\subseteq\mathfrak{r}_{_{P}}^{\triangleright}$ and
$A/\mathfrak{r}_{_{P}}^{\triangleright}$ is an $\mathsf{NG}P$-algebra. On the
other hand, $A/\mathfrak{r}_{_{P}}^{\triangleright}$ is a quotient of the
\textsf{G}$P$-algebra $A/I$. Since the C*-property \textsf{G}$P$ is upper
stable, all the quotients of $A/I$ are \textsf{G}$P$-algebras. Hence
$A/\mathfrak{r}_{_{P}}^{\triangleright}$ is a \textsf{G}$P$-algebra -- a
contradiction. Thus $A=\mathfrak{r}_{_{P}}$ is a \textsf{G}$P$-algebra. So
\textsf{G}$P$ is extension stable.

Let $A\in\mathfrak{A}$ and let $I\in$ Id$_{A}$ and $A/I$ be $\mathsf{NG}%
P$-algebras. Then, by Theorem \ref{T3.1}(i), $\mathfrak{r}_{_{P}%
}^{\triangleright}\subseteq I$ and there is an ascending transfinite
$\ll_{_{P}}$-series of ideals $\left(  I_{\lambda}\right)  _{1\leq\lambda
\leq\gamma}$ from $\{0\}$ to $\mathfrak{r}_{_{P}}^{\triangleright}$. As
$I_{1}=\{0\}$ and $I_{1}\ll_{_{P}}I_{2},$ we have $I_{2}=I_{2}/I_{1}$ is a
$P$-ideal in $\mathfrak{r}_{_{P}}^{\triangleright}.$ If $\mathfrak{r}_{_{P}%
}^{\triangleright}\neq\{0\}$ then $\{0\}\neq I_{2}\subset I$ -- a
contradiction, since $I$ is an $\mathsf{NG}P$-algebra and, therefore, has no
$P$-ideals. Thus $\mathfrak{r}_{_{P}}^{\triangleright}=\{0\}.$ By Theorem
\ref{T3.1}, $A$ is an $\mathsf{NG}P$-algebra. Thus the C*-property
$\mathsf{NG}P$ is extension stable.\bigskip
\end{proof}

If $P$ is both lower and upper stable then $\mathsf{G}P$ is also lower and
upper stable$.$

\begin{proposition}
\label{P3.4}If $P$ is a lower and upper stable \emph{C*}-property\emph{ }then
$\mathsf{G}P$ is also lower and upper stable$.$
\end{proposition}

\begin{proof}
As \textsf{G}$P$ is upper stable by Lemma \ref{L3.6}, we only need to show
that it is lower stable.

Let $A$ be a \textsf{G}$P$-algebra and $J\in$ Id$_{A}.$ We have to prove that
$J\in$ \textsf{G}$P.$ By Theorem \ref{T3.1}, $A$ has an ascending transfinite
$\ll_{_{P}}$-series $\left(  I_{\lambda}\right)  _{1\leq\lambda\leq\gamma}$ of
ideals, $I_{1}=\{0\},$ $I_{\gamma}=A$, $I_{\lambda+1}/I_{\lambda}$ are
$P$-algebras for all $\lambda.$ Then $\left(  I_{\lambda}\cap J\right)
_{1\leq\lambda\leq\gamma}$ is an ascending transfinite series of ideals of $J$
and%
\[
(I_{\lambda+1}\cap J)/(I_{\lambda}\cap J)=(I_{\lambda+1}\cap J)/(I_{\lambda
}\cap(I_{\lambda+1}\cap J))\overset{(\ref{4.01})}{\approx}((I_{\lambda+1}\cap
J)+I_{\lambda})/I_{\lambda}.
\]
As $(I_{\lambda+1}\cap J)+I_{\lambda}$ is an ideal of $I_{\lambda+1},$
$((I_{\lambda+1}\cap J)+I_{\lambda})/I_{\lambda}$ is an ideal of
$I_{\lambda+1}/I_{\lambda}.$ As $I_{\lambda+1}/I_{\lambda}$ is a $P$-algebra
and since the C*-property $P$ is lower stable, $((I_{\lambda+1}\cap
J)+I_{\lambda})/I_{\lambda}$ is a $P$-algebra. Hence $(I_{\lambda+1}\cap
J)/(I_{\lambda}\cap J)$ is a $P$-algebra. So $(I_{\lambda}\cap J)\ll_{_{P}%
}(I_{\lambda+1}\cap J).$

Let $\beta$ be a limit ordinal and $x\in I_{\beta}\cap J.$ As $I_{\beta
}=\overline{\cup_{\lambda<\beta}I_{\lambda}},$ for each $\varepsilon>0,$ there
is $\lambda_{\varepsilon}<\beta$ such that $\min\{\left\Vert x-y\right\Vert $:
$y\in I_{\lambda_{\varepsilon}}\}<\varepsilon.$ As $J/(I_{\lambda
_{\varepsilon}}\cap J)\approx(J+I_{\lambda_{\varepsilon}})/I_{\lambda
_{\varepsilon}}$ by $(\ref{4.01}),$%
\[
\min\{\left\Vert x-z\right\Vert \text{: }z\in I_{\lambda_{\varepsilon}}\cap
J\}=\left\Vert x\right\Vert _{J/(I_{\lambda_{\varepsilon}}\cap J)}=\left\Vert
x\right\Vert _{(J+I_{\lambda_{\varepsilon}})/I_{\lambda_{\varepsilon}}}%
=\min\{\left\Vert x-y\right\Vert \text{: }y\in I_{\lambda_{\varepsilon}%
}\}<\varepsilon.
\]
Hence $I_{\beta}\cap J=\overline{\cup_{\lambda<\beta}(I_{\lambda}\cap J)}.$
Thus $\left(  I_{\lambda}\cap J\right)  _{1\leq\lambda\leq\gamma}$ is an
ascending transfinite $\ll_{_{P}}$-series of ideals in $J,$ $I_{0}\cap
J=\{0\},$ $I_{\gamma}\cap J=J.$ By Corollary \ref{C4.2n}(i), $J$ is a
\textsf{G}$P$-algebra.
\end{proof}

\begin{corollary}
Let $P$\textbf{ }be a \emph{C*}-property,\textbf{ }$A\in\mathfrak{A}$ be
non-unital and $\widehat{A}=A\dotplus\mathbb{C}\mathbf{1.\smallskip}$

\emph{(i) }Let $P$\textbf{ }be upper stable. Then\textbf{ }$A\in$
$\mathsf{G}P$ implies $\widehat{A}\in$ $\mathsf{G}P$ if and only if\textbf{
}$\mathbb{C}\mathbf{1}\in P$\textbf{.}$\smallskip$

\emph{(ii) }If $P$ is lower and upper stable then $\widehat{A}\in$
$\mathsf{G}P$ implies $A\in$ $\mathsf{G}P.$
\end{corollary}

\begin{proof}
(i) By Proposition \ref{P3.3P}, the C*-property \textsf{G}$P$ is extension
stable. As $A\in$ Id$_{\widehat{A}},$ $A\in$ \textsf{G}$P$ and $\widehat{A}%
/A\approx\mathbb{C}\mathbf{1,}$ we have that $\mathbb{C}\mathbf{1}\in P$
implies $\widehat{A}\in$ \textsf{G}$P.$

Conversely, if $\widehat{A}\in$ \textsf{G}$P$ then $\widehat{A}/A\approx
\mathbb{C}\mathbf{1}\in P.$

(ii) If $P$ is lower stable, \textsf{G}$P$ is lower stable by Proposition
\ref{P3.4}. Hence $\widehat{A}\in$ \textsf{G}$P$ implies that its ideal $A$
belongs to \textsf{G}$P.\bigskip$
\end{proof}

The set $\mathcal{P}_{\text{up}}$ of all upper stable C*-properties in
$\mathfrak{A}$ is a complete lattice with $P\leq P_{1}$ if $P\subseteq P_{1}.$
For a subset $P_{\Lambda}=\{P_{\lambda}\}_{\lambda\in\Lambda}$ in
$\mathcal{P}_{\text{up}},$ set
\begin{equation}
\wedge P_{\Lambda}=\cap_{_{\lambda\in\Lambda}}P_{\lambda}\text{ and }\vee
P_{\Lambda}=\cup_{_{\lambda\in\Lambda}}P_{\lambda}. \label{3.8}%
\end{equation}

\begin{proposition}
\label{L3.1}The map $P\mapsto\mathsf{G}P$ is a closure operator in
$\mathcal{P}_{\text{\emph{up}}},$ that is\emph{,}%
\begin{equation}
P\subseteq\mathsf{G}P=\mathsf{G}(\mathsf{G}P)\text{ for }P\in\mathcal{P}%
_{\text{\emph{up}}}\text{ and }\mathsf{G}P\subseteq\mathsf{G}P_{1}.
\label{3.13}%
\end{equation}
if $P\subseteq P_{1}$ in $\mathcal{P}_{\text{\emph{up}}}.$ In this case
$\mathfrak{r}_{_{P}}^{\triangleright}\subseteq\mathfrak{r}_{_{P_{1}}%
}^{\triangleright}$ in $\emph{Id}_{A}$ for each $A\in\mathfrak{A}.$

If $P\subseteq P_{1}\subseteq\mathsf{G}P,$ then $\mathsf{G}P=\mathsf{G}P_{1}$
and $\mathfrak{r}_{_{P}}^{\triangleright}=\mathfrak{r}_{_{P_{1}}%
}^{\triangleright}$ in $\emph{Id}_{A}$ for all $A\in\mathfrak{A}.$
\end{proposition}

\begin{proof}
If $A\in P$ and $P$ is upper stable, all quotients of $A$ are $P$-algebras. So
$A\in\mathsf{G}P$. Thus $P\subseteq\mathsf{G}P.$ Hence \textsf{G}%
$P\subseteq\mathsf{G}(\mathsf{G}P).$

Let $A\in\mathsf{G}(\mathsf{G}P)$ and $I\in$ Id$_{A},$ $I\neq A.$ Then $A/I$
has a \textsf{G}$P$-ideal $J\neq\{0\}.$ Hence $J$ itself has a $P$-ideal
$K\neq\{0\}.$ As $K$ is also an ideal of $A/I$, we have that $A$ is a
\textsf{G}$P$-algebra. So \textsf{G}$P=\mathsf{G}(\mathsf{G}P).$

Let $\{0\}\neq A\in\mathsf{G}P.$ Then each non-zero quotient of $A$ has a
non-zero $P$-ideal. As $P\subseteq P_{1},$ this ideal is also an $P_{1}%
$-ideal. So $A\in\mathsf{G}P_{1}.$ Thus $\mathsf{G}P\subseteq\mathsf{G}P_{1}.$
So (\ref{3.13}) holds and the map $P\mapsto\mathsf{G}P$ is a closure operator
in $\mathcal{P}_{\text{up}}$ (see (\ref{2.6})).

If $\mathfrak{r}_{_{P}}^{\triangleright}\nsubseteq\mathfrak{r}_{_{P_{1}}%
}^{\triangleright}$ in Id$_{A}$ then $\mathfrak{r}_{_{P_{1}}}^{\triangleright
}$ has a $\ll_{_{P}}$-successor by Corollary \ref{C4.2n}. Thus, as $P\subseteq
P_{1},$ $\mathfrak{r}_{_{P_{1}}}^{\triangleright}$ has a $\ll_{_{P_{1}}}%
$-successor which contradicts Corollary \ref{C4.2n}. So $\mathfrak{r}_{_{P}%
}^{\triangleright}\subseteq\mathfrak{r}_{_{P_{1}}}^{\triangleright}$.

As $P\mapsto\mathsf{G}P$ is a closure operator, $P\subseteq P_{1}%
\subseteq\mathsf{G}P$ implies $\mathsf{G}P=\mathsf{G}P_{1}$ by Lemma
\ref{L2.1}. So $\mathfrak{r}_{_{\mathsf{G}P}}=\mathfrak{r}_{_{\mathsf{G}P_{1}%
}}$ in Id$_{A}$ for $A\in\mathfrak{A}.$ By Theorem \ref{T3.0}, we have
$\mathfrak{r}_{_{\mathsf{G}P}}=\mathfrak{r}_{_{P}}^{\triangleright}.$ Thus
$\mathfrak{r}_{_{P}}^{\triangleright}=\mathfrak{r}_{_{P_{1}}}^{\triangleright
}.\bigskip$
\end{proof}

Note that if $P$ is not upper stable then $P\nsubseteq$ $\mathsf{G}P,$ as
$\mathsf{G}P$ is upper stable by Lemma \ref{L3.6}.

\begin{example}
\emph{Let }$P_{un}$\emph{ and }$P_{n\text{-}un}$\emph{ be the C*-properties of
all unital and non-unital C*-algebras, respectively}.\smallskip

$1)$ $P_{n\text{-}un}$\textbf{ }\emph{is neither lower, nor upper stable,
since} $A=C(H)\oplus\mathbb{C}\mathbf{1}$ \emph{is non-unital, while the ideal
}$\mathbb{C}\mathbf{1}$ \emph{and the quotient} $A/C(H)\approx\mathbb{C}%
\mathbf{1}$ \emph{are unital}. \emph{Moreover,} $P_{un}\nsubseteq$
$\mathsf{G}P_{n\text{-}un},$ \emph{since} $B=C(H)+\mathbb{C}\mathbf{1}_{H}\in
P_{un},$ \emph{but} $B\neq$ $\mathsf{G}P_{n\text{-}un},$ \emph{as}
$B/C(H)\approx\mathbb{C}\mathbf{1}$ \emph{does not have non-zero non-unital
ideals.\smallskip}

$2)$ \emph{The C*-property }$P_{un}$\emph{ is upper stable, so that }%
$P_{un}\subseteq$ $\mathsf{G}P_{un}.$ \emph{However, }$P_{un}\neq$
$\mathsf{G}P_{un}.$ \emph{Indeed, let }$\{A_{\lambda}\}_{\lambda\in\Lambda}%
$\emph{ be unital C*-algebras with identities }$e_{\lambda}$. \emph{The
C*-algebra }$A(\Lambda)$\emph{ of all sequences }$(a_{\lambda})_{\lambda
\in\Lambda},$\emph{ }$a_{\lambda}\in A_{\lambda}$\emph{ for all }$\lambda
\in\Lambda,$\emph{ with }$\left\Vert (a_{\lambda})\right\Vert =\sup
\{\left\Vert a_{\lambda}\right\Vert $\emph{: }$\lambda\in\Lambda\}<\infty
$\emph{ is unital.}

\emph{Let }$A_{0}$\emph{ be the set of all }$(a_{\lambda})_{\lambda\in\Lambda
}\in A(\Lambda)$\emph{ such that only finite number of them are non-zero. Then
the closure }$A_{0}(\Lambda)$\emph{ of }$A_{0}$\emph{ in }$A(\Lambda)$\emph{
is a non-unital C*-algebra. Let us show that }$A_{0}(\Lambda)\in$
$\mathsf{G}P_{un}.$

\emph{Let }$I\in$\emph{ Id}$_{A_{0}(\Lambda)}.$\emph{ For each }$\lambda
\in\Lambda,$\emph{ }$e_{\lambda}I\subset I$\emph{ and can be considered as an
ideal of }$A_{\lambda}.$\emph{ If }$e_{\lambda}I=A_{\lambda}$\emph{ for all
}$\lambda\in\Lambda,$\emph{ then }$I=A_{0}(\Lambda).$\emph{ If }$I\neq
A_{0}(\Lambda),$\emph{ }$e_{\lambda_{0}}I\neq A_{\lambda_{0}}$\emph{ for some
}$\lambda_{0}\in\Lambda.$\emph{ Then }$A_{\lambda_{0}}/(e_{\lambda_{0}}%
I)$\emph{ is unital and can be considered as an ideal of }$A_{0}(\Lambda
)/I.$\emph{ Thus }$A_{0}(\Lambda)\in$\emph{ }$\mathsf{G}P_{un}.$ \emph{So
}$P_{un}\neq$ $\mathsf{G}P_{un}.$ \ \ \ $\blacksquare$
\end{example}

By Proposition \ref{L3.1}, the C*-property \textsf{G}$P$ can me much larger
than $P,$ if $P$ is upper stable. We consider now the conditions under which
they coincide, so that $\mathfrak{r}_{_{P}}^{\triangleright}(A)\in P$ for all
$A\in\mathfrak{A}$ (by Theorem \ref{T3.1} (i)).

In fact, any closure operator $f$ is completely characterized by stable
elements: $f(P)$ is the smallest of all $f$-stable C*-properties containing
$P$. So the following result gives important information about the map
$P\rightarrow\mathsf{G}P$.

\begin{corollary}
\label{C3.2}Let $P$ be an upper stable \emph{C*}-property. The following
conditions are equivalent.\smallskip

\emph{(i) \ \ }$\mathsf{G}P=P;$\smallskip

\emph{(ii) \ }$\ll_{_{P}}=$ $\ll_{_{P}}^{\triangleright}$ in \emph{Id}$_{A}$
for each $A\in\mathfrak{A};$\smallskip

\emph{(iii) }$\ll_{_{P}}$ is an $\mathbf{R}$-order in \emph{Id}$_{A}$ for each
$A\in\mathfrak{A};$\smallskip

\emph{(iv) }$P$ is extension stable and the closure of the union of any
ascending transfinite series of $P$-ideals is a $P$-ideal in each
$A\in\mathfrak{A.}$

\emph{(v) }$P$ is extension stable and the closure of the sum of any family of
$P$-ideals in each $A\in\mathfrak{A}$ is a $P$-ideal.
\end{corollary}

\begin{proof}
(i) $\Rightarrow$ (ii) follows from Theorem \ref{T3.0}(i).\smallskip

(ii) $\Leftrightarrow$ (iii) follows from Theorem \ref{inf}(i).\smallskip

(iii) $\Rightarrow$ (iv). Let $\ll_{_{P}}$ be an $\mathbf{R}$-order in
Id$_{A}.$ If $I,$ $A/I\in P$ then $\{0\}\ll_{_{P}}I$ and $I\ll_{_{P}}A.$ As
$\ll_{_{P}}$ is transitive, $\{0\}\ll_{_{P}}A,$ so that $A\in P.$ Thus (see
(\ref{3.7})) $P$ is extension stable.

Let $\left(  I_{\lambda}\right)  _{1\leq\lambda<\beta}$ be an ascending
transfinite series of $P$-ideals of $A$ for some limit ordinal $\beta.$ In
this case (see Section 11.4 \cite{KR}) $I_{\beta}=\overline{\cup
_{\lambda<\beta}I_{\lambda}}$ is its inductive limit . We have $\{0\}\ll
_{_{P}}I_{\lambda}$ for all $\lambda<\beta,$ so that $I_{\lambda}\in
\lbrack\{0\},\ll_{_{P}}]$ (see (\ref{1})). As $\ll_{_{P}}$ is an $\mathbf{R}%
$-order in Id$_{A}$, the set $[\{0\},\ll_{_{P}}]$ is $\vee$-complete (see
(\ref{1.5})). Hence $\vee\left(  I_{\lambda}\right)  _{1\leq\lambda<\beta
}=\overline{\cup_{\lambda<\beta}I_{\lambda}}=I_{\beta}\in\lbrack
\{0\},\ll_{_{P}}].$ So $\{0\}\ll_{_{P}}I_{\beta},$ i.e., $I_{\beta}\in
P.\smallskip$

(iv) $\Rightarrow$ (i). By Proposition \ref{L3.1}, $P\subseteq\mathsf{G}P.$
Let $A\in\mathsf{\ }P.$ By Theorem \ref{T3.1}, there is an ascending
transfinite $\ll_{_{P}}$-series of ideals $\left(  I_{\lambda}\right)
_{1\leq\lambda\leq\gamma}$ such that $I_{1}=\{0\},$ $I_{\gamma}=A$ and
$I_{\lambda+1}/I_{\lambda}\in P$ for all $\lambda.$ Suppose that $I_{\lambda
}\in P$ for some $\lambda$. Then $I_{\lambda+1}\in P,$ as $P$ is extension
stable. Let $\beta$ be a limit ordinal and all $I_{\lambda}\in P,$
$\lambda<\beta$. By (\ref{a}), $I_{\beta}=\overline{\cup_{\lambda<\beta
}I_{\lambda}}$, i.e.,\ $I_{\beta}$ is an inductive limit of $P$-algebras.
Hence $I_{\beta}\in P$. Thus, by transfinite induction, $I_{\gamma}=A\in P$.
So \textsf{G}$P=P.$

(i) $\Rightarrow$ (v). It suffices to show that if $\{I_{\lambda}%
\}_{\lambda\in\Lambda}$ are $P$-ideals then $I=\overline{\sum_{\lambda
\in\Lambda}I_{\lambda}}\in P.$ Indeed, let $J$ be a proper ideal of $I$ and
$q$: $I\rightarrow I/J$ be the standard epimorphism. If $q(I_{\lambda})=0$ for
all ${\lambda}$, then $q(I)=0$. So $J=I$, a contradiction. Hence
$q(I_{\lambda})\neq0$ for some ${\lambda}$. So $q(I_{\lambda})\approx
(I_{\lambda}+J)/J\overset{(\ref{4.01})}{\approx}I_{\lambda}/(I_{\lambda}\cap
J)$ is an ideal of $I/J.$ As $\mathsf{G}P=P,$ $I_{\lambda}$ is a $\mathsf{G}%
P$-ideal. Thus $q(I_{\lambda})$ has a non-zero $P$-ideal which is a $P$-ideal
of $I/J.$ So $I$ is a $\mathsf{G}P$-algebra. As $\mathsf{G}P=P,$ $I$ is a
$P$-algebra.\bigskip
\end{proof}

We consider now an analogue of Corollary \ref{C3.2} when $P\neq\mathsf{G}P,$
but $P=\mathsf{G}P\cap\mathfrak{M}$ for some lower and upper stable
C*-property $\mathfrak{M}.$

\begin{corollary}
\label{C3.3}Let $P$ be an upper stable C*-property contained in a lower and
upper stable C*-property $\mathfrak{M}.$ The following conditions are
equivalent.\smallskip

\emph{(i) \ \ }$\mathsf{G}P\cap\mathfrak{M}=P;$\smallskip

\emph{(ii) \ }$\ll_{_{\mathsf{G}P}}=$ $\ll_{_{P}}=$ $\ll_{_{P}}%
^{\triangleright}$ in \emph{Id}$_{A}$ for each $A\in\mathfrak{M};$\smallskip

\emph{(iii) }$\ll_{_{P}}$ is an $\mathbf{R}$-order in \emph{Id}$_{A}$ for each
$A\in\mathfrak{M}$.
\end{corollary}

\begin{proof}
(i) $\Rightarrow$ (ii). Let $A\in\mathfrak{M.}$ First, let us show that
$\ll_{_{\mathsf{G}P}}=$ $\ll_{_{P}}$ in Id$_{A}.$ If $I\ll_{_{P}}J$ in
Id$_{A},$ then $J/I\in P.$ As $P\subseteq$ \textsf{G}$P$, $J/I\in$
\textsf{G}$P.$ So $I\ll_{_{\mathsf{G}P}}J.$

Conversely, let $I\ll_{_{\mathsf{G}P}}J.$ Then $J/I\in\mathsf{G}P.$ As
$\mathfrak{M}$ is lower stable, $J\in\mathfrak{M.}$ As $\mathfrak{M}$ is upper
stable, $J/I\in\mathfrak{M.}$ Hence $J/I\in\mathsf{G}P\cap\mathfrak{M}=P.$ So
$I\ll_{_{P}}J.$ Thus $\ll_{_{\mathsf{G}P}}=$ $\ll_{_{P}}$ in Id$_{A}$ for each
$A\in\mathfrak{M.}$

By Theorem \ref{T3.0}(i), \emph{\ }$\ll_{_{\mathsf{G}P}}=$ $\ll_{_{P}%
}^{\triangleright}$ in Id$_{A}$ for each $A\in\mathfrak{A}$ which completes
the proof of (ii).\smallskip

(ii) $\Leftrightarrow$ (iii) follows from Theorem \ref{inf}(i).

(ii) $\Rightarrow$ (i). We have $P\subseteq\mathsf{G}P\cap\mathfrak{M.}$ Let
$A\in\mathsf{G}P\cap\mathfrak{M.}$ Then $\{0\}\ll_{_{\mathsf{G}P}}A.$ As
$A\in\mathfrak{M,}$ it follows from (ii) that $\{0\}\ll_{_{P}}A.$ So $A\in P.$
Thus \textsf{G}$P\cap\mathfrak{M}=P.$
\end{proof}

\subsection{The closure operator $P\mapsto\mathsf{dG}P$ on lower stable
properties}

Following the duality principle, we introduce now two operations dual to the
operations $P\rightarrow$ \textsf{G}$P$ and $P\rightarrow\mathsf{NG}P$. We
showed in Theorem \ref{T3.1} that the $\ll_{_{P}}^{\triangleright}$-radical
$\mathfrak{r}_{_{P}}^{\triangleright}$ is the "dividing line" between
\textsf{G}$P$ and $\mathsf{NG}P$ parts in each $A\in\mathfrak{A}$. In this
section we will show that the dual $\ll_{_{P}}^{\triangleleft}$-radical
$\mathfrak{p}_{_{P}}^{\triangleleft}$ is the "dividing line" between
$\mathsf{dG}P$ and \textsf{d}$\mathsf{NG}P$ parts of $A$.

\begin{definition}
\label{D3.2}Let $P$ be a \emph{C*}-property in $\mathfrak{A.}$ We call a
$C^{\ast}$-algebra $A\smallskip$

\emph{(i) }a $\mathsf{dG}P$\textbf{-algebra} \emph{(}dual generalized
$P$-algebra\emph{) }if either $A=\{0\},$ or each $\{0\}\neq J\in$
\emph{Id}$_{A}$ has a non-zero quotient which is a $P$-algebra\emph{: }there
is $I\in\lbrack\{0\},J)$ such that $J/I\in P,$\emph{ }i.e.\emph{,} $I\ll
_{_{P}}J.\smallskip$

\emph{(ii) }a $\mathsf{dNG}P$\textbf{-algebra} if each non-zero quotient of
$A$ is not a $P$-algebra.
\end{definition}

Denote by $\mathsf{dG}P$ and \textsf{d}$\mathsf{NG}P$ the classes of all
$\mathsf{dG}P$- and all \textsf{d}$\mathsf{NG}P$-algebras in $\mathfrak{A.}$
They are C*-properties. As in (\ref{4.1}), we consider the corresponding
relations $\ll_{_{\text{\textsf{d}}\mathsf{G}P}}$ and $\ll_{_{\text{\textsf{d}%
}\mathsf{NG}P}}$ in Id$_{A}.$

\begin{lemma}
\label{L3.7}For each \emph{C*}-property $P,$ the \emph{C*}-property
$\mathsf{dG}P$ is lower stable and $\mathsf{dNG}P$ is upper stable.
\end{lemma}

\begin{proof}
Let $A\in$ $\mathsf{dG}P$ and $\{0\}\neq K\in$ Id$_{A}.$ Each $\{0\}\neq J\in$
Id$_{K}$ also belongs to Id$_{A}.$ As $A\in$ $\mathsf{dG}P,$ there is $I\in$
Id$_{A}$ such that $\{0\}\neq J/I\in P.$ As $I\in$ Id$_{K},$ we conclude that
$K\in$ $\mathsf{dG}P.$ Thus $\mathsf{dG}P$ is a lower stable C*-property.

Let $A\in$ \textsf{d}$\mathsf{NG}P$ and $J\in$ Id$_{A}.$ As each non-zero
quotient of $A/J$ is isomorphic to a non-zero quotient of $A,$ it is not a
$P$-algebra. Thus $A/J$ is a \textsf{d}$\mathsf{NG}P$-algebra. So
\textsf{d}$\mathsf{NG}P$ is upper stable.
\end{proof}

\begin{example}
\label{exDG} \emph{1) }If $P=Comm$ then $\mathsf{dG}P=Comm$\emph{ (see Example
\ref{exG})}.

\emph{Indeed, as $P$ is lower stable, $P\subseteq\mathsf{dG}P$. Let
$A\in\mathsf{dG}P,$ }$\Pi_{1}(A)$ \emph{be the set of all one-dimensional
representations of $A$ and }$J=\cap\{\ker\pi$\emph{: }$\pi\in\Pi_{1}%
(A)\}.$\emph{ If $J\neq\{0\}$ then, by Definition \ref{D3.2}, }$\{0\}\neq
B:=J/I\in P$\emph{ for some }$I\in$ \emph{Id}$_{J}.$ \emph{Let} \emph{$f$:
$J\rightarrow B$ be the standard epimorphism, }$0\neq\pi\in\Pi_{1}(B)$\emph{
and }$\tau:=\pi$\emph{$\circ f.$ Then }$\tau\in\Pi_{1}($\emph{$J)$. It extends
to a representation }$\tau^{\prime}$\emph{ of $A$ on the same one-dimensional
space. As }$\tau^{\prime}(J)=\tau(G)\neq\{0\},$\emph{ we have a contradiction
with definition of }$J.$\emph{ So $J=\{0\}$. }

\emph{To prove that }$A$\emph{ is commutative, assume that $xy\neq yx$ for
some $x,y\in A$. Then there is $\pi\in\Pi_{1}(A)$ with $\pi(xy-yx)\neq0$.
However, $\pi(xy-yx)=\pi$}$(x)\pi(y)-\pi(y)\pi(x)=0,$ \emph{a contradiction.
Thus }$A\in Comm.$\emph{ So }$\mathsf{dG}P=Comm.\smallskip$

\emph{2) }If\emph{ $P=\mathbb{C}\mathbf{1}$ }is the class of
one-dimensional\emph{ C*-}algebras\emph{, }then\emph{ $\mathsf{dG}P=Comm$}.
\emph{As $P\subset Comm$, $\mathsf{dG}P\subset\mathsf{dG}Comm=Comm$. On the
other hand, if $A\in Comm$, then any }$I\in$ \emph{Id}$_{A}$\emph{ is
commutative and, therefore, has a one-dimensional representation. So
$A\in\mathsf{dG}P$. Thus $\mathsf{dG}P=Comm$.\smallskip}

\emph{3) }Let\emph{ }$\Pi_{f}(A)=\{\pi\in\Pi(A)$\emph{: }$\dim\pi<\infty\}$
and%
\begin{equation}
P=RFD=\{A\in\mathfrak{A}\text{\emph{:}}\emph{\ }\cap\{\ker\pi\text{\emph{:}%
}\emph{\ }\pi\in\Pi_{f}(A)\}=\{0\}\} \label{3.30}%
\end{equation}
\emph{ }be the set of all C*-algebras $A$ with a separating set $\Pi_{f}(A).$
Then\emph{ $\mathsf{dG}P=P$.}

\emph{ Note that $P$ is lower stable. Indeed, if $I\in$ Id}$_{A}$\emph{ and
}$A\in\emph{P}$\emph{ then, by definition, for each $a\in I$, there is $\pi
\in\Pi$}$_{f}(A)$\emph{ with $\pi(a)\neq0$. So the restrictions to }$I$\emph{
of all $\pi\in\Pi$}$_{f}(A)$\emph{ separate }$I$\emph{. So $I\in P$. Hence
$P\subset\mathsf{dG}P$. Let $A\in\mathsf{dG}P$ and }$J=\cap\{\ker\pi$\emph{:
}$\pi\in\Pi_{f}(A)\}.$ \emph{Repeating the proof of 1), we get that
}$J=\{0\}.$ \emph{So }$A\in\emph{P}$\emph{. Thus $\mathsf{dG}P=P$.}
\end{example}

If $P$ is lower stable, $\ll_{_{P}}$ is a dual $\mathbf{H}$-relation and
$\ll_{_{P}}^{\triangleleft}$ is a dual $\mathbf{R}$-order in Id$_{A}$ by
Theorem \ref{T4.1}.

We have that $I\ll_{_{\text{\textsf{d}}\mathsf{G}P}}J$ in Id$_{A}$ if and only
if $J/I\in$ $\mathsf{dG}P.$ So, by Definition \ref{D3.2},%
\begin{equation}
I\ll_{_{\text{\textsf{d}}\mathsf{G}P}}J\Longleftrightarrow\text{ for each
}K\in(I,J],\text{ there is }L\in\lbrack I,K)\text{ such that }L\ll_{_{P}}K.
\label{3.22}%
\end{equation}

\begin{theorem}
\label{T3.7}Let $P$ be a lower stable \emph{C*}-property in $\mathfrak{A}$ and
$A\in\mathfrak{A}.$ Then$\smallskip$

\emph{(i) \ \ }$\ll_{_{\mathsf{dG}P}}=$ $\ll_{_{P}}^{\triangleleft}$ in
\emph{Id}$_{A},$ so that $\ll_{_{\mathsf{dG}P}}$ is a dual $\mathbf{R}$-order
in \emph{Id}$_{A}\emph{;}$\smallskip

\emph{(ii) \ }$\ll_{_{\mathsf{dNG}P}}=$ $\overrightarrow{\ll_{_{P}}}$ in
\emph{Id}$_{A},$ so that $\ll_{_{\mathsf{dNG}P}}$ is an $\mathbf{R}$-order in
\emph{Id}$_{A}$\emph{;\smallskip}

\emph{(iii)} $\mathfrak{p}_{_{P}}^{\triangleleft}=\mathfrak{p}_{_{\mathsf{dG}%
P}}=\mathfrak{r}_{_{\mathsf{dNG}P}}.\smallskip$

\emph{(iv) }If $\mathfrak{p}_{_{P}}^{\triangleleft}\cap K=\{0\}$ for $K\in$
\emph{Id}$_{A},$ then $K$ is a $\mathsf{dG}P$-algebra.
\end{theorem}

\begin{proof}
As $P$ is lower stable, $\ll_{_{P}}$ is a dual $\mathbf{H}$-relation in
Id$_{A}$ by Theorem \ref{T4.1}.

(i) Let $I\ll_{_{\text{\textsf{d}}\mathsf{G}P}}J,$ $I\neq J,$ in Id$_{A}.$ By
(\ref{2.4}) and (\ref{3.22}), $I\ll_{_{\text{\textsf{d}}\mathsf{G}P}}J$ if and
only if each $K\in(I,J]$ has a $\ll_{_{P}}$-predecessor$,$ i.e., $[I,J]$ is an
lower $\ll_{_{P}}$-set. Thus, by (\ref{2.1}), $I\ll_{_{\text{\textsf{d}%
}\mathsf{G}P}}J\Longleftrightarrow I\ll_{_{P}}^{\text{lo}}J$. So
$\ll_{_{\text{\textsf{d}}\mathsf{G}P}}=$ $\ll_{_{P}}^{\text{lo}}.$ As
$\ll_{_{P}}$ is a dual $\mathbf{H}$-relation, $\ll_{_{P}}^{\text{lo}}=$
$\ll_{_{P}}^{\triangleleft}$ by Theorem \ref{inf}. Thus $\ll
_{_{\text{\textsf{d}}\mathsf{G}P}}=$ $\ll_{_{P}}^{\triangleleft}$ is a dual
$\mathbf{R}$-order in Id$_{A}$ by Theorem \ref{T4.1}.\smallskip

(ii) Let $I\ll_{_{\text{\textsf{d}}\mathsf{NG}P}}J$ and $I\neq J.$ Then $J/I$
is a \textsf{d}$\mathsf{NG}P$-algebra, so that each its non-zero quotient is
not a $P$-algebra. Hence, by (\ref{3.22}), for each $K\in\lbrack I,J),$ $J/K$
is not a $P$-algebra, i.e., $K\not \ll _{_{P}}J.$ Thus $I\ll
_{_{\text{\textsf{d}}\mathsf{NG}P}}J$ if and only if $[\ll_{_{P}}%
,J]\cap\lbrack I,J]=J.$ So, by (\ref{A1}), $I\ll_{_{\text{\textsf{d}%
}\mathsf{NG}P}}J\Longleftrightarrow I$ $\overrightarrow{\ll_{_{P}}}$ $J.$
Hence $\ll_{_{\text{\textsf{d}}\mathsf{NG}P}}$ $=$ $\overrightarrow{\ll_{_{P}%
}}$ is an $\mathbf{R}$-order in Id$_{A}$ by Theorem \ref{inf}(ii).\smallskip

(iii) By (i), $\mathfrak{p}_{_{P}}^{\triangleleft}=\mathfrak{p}%
_{_{\text{\textsf{d}}\mathsf{G}P}}.$ As $\ll_{_{\text{\textsf{d}}\mathsf{NG}%
P}}$ $=$ $\overrightarrow{\ll_{_{P}}}$ is an $\mathbf{R}$-order by (ii),
$\mathfrak{r}_{_{\text{\textsf{d}}\mathsf{NG}P}}$ coincides with the
$\overrightarrow{\ll_{_{P}}}$-radical which, in turn, coincides with
$\mathfrak{p}_{_{P}}^{\triangleleft}$ by Theorem \ref{inf}(ii).$\smallskip$

(iv) If $\{0\}\neq I\in$ Id$_{K}$ then $I\in$ Id$_{A}$ and $\mathfrak{p}%
_{_{P}}\cap I=\{0\}.$ By Corollary \ref{C4.2n}(ii), there is $J\in$ Id$_{A},$
$J\subsetneqq I,$ such that $I/J\neq\{0\}$ is a $P$-algebra. Thus $K$ is a
$\mathsf{dG}P$-algebra.$\bigskip$
\end{proof}

The set $\mathcal{P}_{\text{lo}}$ of all lower stable C*-properties in
$\mathfrak{A}$ is a complete lattice with $P\leq P_{1}$ if $P\subseteq P_{1}$
and $\wedge$ and $\vee$ defined in (\ref{3.8}).

\begin{proposition}
\label{P3.7}\emph{(i)}\textbf{ }The map $P\mapsto\mathsf{dG}P$ is a closure
operator in $\mathcal{P}_{\text{\emph{lo}}},$ that is\emph{,}%
\[
P\subseteq\mathsf{dG}P=\mathsf{dG}(\mathsf{dG}P)\text{ for }P\in
\mathcal{P}_{\text{\emph{lo}}},\text{ and }\mathsf{dG}P\subseteq
\mathsf{dG}P_{1},
\]
if $P\subseteq P_{1}$ in $\mathcal{P}_{\text{\emph{lo}}}.$ In this case
$\mathfrak{p}_{_{P_{1}}}^{\triangleleft}\subseteq\mathfrak{p}_{_{P}%
}^{\triangleleft}$ in $\emph{Id}_{A}$ for each $A\in\mathfrak{A}.\smallskip$

\emph{(ii) }If $P,P_{1}\in\mathcal{P}_{\text{\emph{lo}}}$ and $P\subseteq
\mathsf{dG}P_{1},$ then $\mathsf{dG}P\subseteq\mathsf{dG}P_{1}.\smallskip$

\emph{(iii) }If $P\subseteq P_{1}\subseteq\mathsf{dG}P$ in $\mathcal{P}%
_{\text{\emph{lo}}},$ then $\mathsf{dG}P=\mathsf{dG}P_{1}$ and $\mathfrak{p}%
_{_{P}}^{\triangleleft}=\mathfrak{p}_{_{P_{1}}}^{\triangleleft}$ in
$\emph{Id}_{A}$ for $A\in\mathfrak{A}.$
\end{proposition}

\begin{proof}
(i) The inclusions $P\subseteq$ $\mathsf{dG}P$ and $\mathsf{dG}P\subseteq
\mathsf{dG}P_{1}$ for $P\subseteq P_{1}$ in $\mathcal{P}_{\text{lo}}$ are obvious.

As $P\subseteq$ $\mathsf{dG}P,$ we have $\mathsf{dG}P\subseteq$ $\mathsf{dG}%
(\mathsf{dG}P).$ Let $A\in$ $\mathsf{dG}(\mathsf{dG}P)$ and $J\in$ Id$_{A}.$
Then $J/I$ is a $\mathsf{dG}P$-algebra for some $I\subsetneqq J.$ Hence there
is an ideal $\widehat{K}\subsetneqq J/I$ such that $(J/I)/\widehat{K}$ is a
$P$-algebra$.$ So there is an ideal $K\in\lbrack I,J),$ $\widehat{K}=K/J,$
such that $J/K\approx(J/I)/\widehat{K}$ is a $P$-algebra. Thus $A\in$
$\mathsf{dG}P.$ So $\mathsf{dG}(\mathsf{dG}P)=\mathsf{dG}P.$ So the map
$P\mapsto\mathsf{dG}P$ is a closure operator in $\mathcal{P}_{\text{lo}}$ (see
(\ref{2.6})).

If $\mathfrak{p}_{_{P_{1}}}^{\triangleleft}\nsubseteq\mathfrak{p}_{_{P}%
}^{\triangleleft}$ in Id$_{A}$ then, by Corollary \ref{C4.2n}, $\mathfrak{p}%
_{_{P_{1}}}^{\triangleleft}$ has a $\ll_{_{P}}$-predecessor $I$ in Id$_{A}$.
As $P\subseteq P_{1},$ $I$ is also a $\ll_{_{P_{1}}}$-predecessor of
$\mathfrak{p}_{_{P_{1}}}^{\triangleleft}$ which contradicts Corollary
\ref{C4.2n}. So $\mathfrak{p}_{_{P_{1}}}^{\triangleleft}\subseteq
\mathfrak{p}_{_{P}}^{\triangleleft}.$

(ii) Let $P,P_{1}\in\mathcal{P}_{\text{lo}}.$ By Lemma \ref{L3.7},
$\mathsf{dG}P_{1}$ is lower stable. If $P\subseteq\mathsf{dG}P_{1}$ then, by
above, $\mathsf{dG}P\subseteq\mathsf{dG(dG}P_{1})=\mathsf{dG}P_{1}.$

(iii) As $P\mapsto\mathsf{dG}P$ is a closure operator, $P\subseteq
P_{1}\subseteq\mathsf{dG}P$ implies $\mathsf{dG}P=\mathsf{dG}P_{1}$ by Lemma
\ref{L2.1}. So $\mathfrak{p}_{_{\text{\textsf{d}}\mathsf{G}P}}=\mathfrak{p}%
_{_{\text{\textsf{d}}\mathsf{G}P_{1}}}$ in Id$_{A}$ for $A\in\mathfrak{A}.$ By
Theorem \ref{T3.7}(iii), $\mathfrak{p}_{_{P}}^{\triangleleft}=\mathfrak{p}%
_{_{\text{\textsf{d}}\mathsf{G}P}}$ in Id$_{A}.$ Thus $\mathfrak{p}_{_{P}%
}^{\triangleleft}=\mathfrak{p}_{_{P_{1}}}^{\triangleleft}$.\bigskip
\end{proof}

The following theorem is an analogue of Theorem \ref{T3.1} for lower stable C*-properties.

\begin{theorem}
\label{T3.2}Let $P$ be a lower stable \emph{C*}-property in $\mathfrak{A}$ and
$A\in\mathfrak{A}$. Then\smallskip

\emph{(i) }\ The radical $\mathfrak{p}_{_{P}}^{\triangleleft}$ is the largest
$\mathsf{dNG}P$-ideal of $A$ and the smallest ideal with $\mathsf{dG}%
P$-quotient$.$ There is a descending transfinite $\ll_{_{P}}$-series of ideals
from $A$ to $\mathfrak{p}_{_{P}}^{\triangleleft}$.\smallskip

\emph{(ii) }$A$ is a $\mathsf{dG}P$-algebra if and only if $\mathfrak{p}%
_{_{P}}^{\triangleleft}=\{0\};$ it is a $\mathsf{dNG}P$-algebra if and only if
$\mathfrak{p}_{_{P}}^{\triangleleft}=A.$
\end{theorem}

\begin{proof}
(i) By (\ref{3.4}), $\{0\}$ $\overrightarrow{\ll_{_{P}}^{\triangleleft}}$
$\mathfrak{p}_{_{P}}^{\triangleleft}.$ As $\ll_{_{\text{\textsf{d}}%
\mathsf{NG}P}}$ $=$ $\overrightarrow{\ll_{_{P}}^{\triangleleft}}$ by Theorem
\ref{T3.7}, we have $\{0\}$ $\ll_{_{\text{\textsf{d}}\mathsf{NG}P}}$
$\mathfrak{p}_{_{P}}^{\triangleleft}.$ So $\mathfrak{p}_{_{P}}^{\triangleleft
}$ is a \textsf{d}$\mathsf{NG}P$-ideal in $A$. If $I$ is a \textsf{d}%
$\mathsf{NG}P$-ideal$,$ all quotients of $I$ are not $P$-algebras. Thus $I$
has no $\ll_{_{P}}$-predecessor. By Corollary \ref{C4.2n}, $I\subseteq
\mathfrak{p}_{_{P}}^{\triangleleft}.$ So $\mathfrak{p}_{_{P}}^{\triangleleft}$
contains all \textsf{d}$\mathsf{NG}P$-ideals of $A.$ By (\ref{3.4}),
$\mathfrak{p}_{_{P}}^{\triangleleft}$ $\ll_{_{P}}^{\triangleleft}$ $A.$ As
\emph{\ }$\ll_{_{\text{\textsf{d}}\mathsf{G}P}}=$ $\ll_{_{P}}^{\triangleleft}$
by Theorem \ref{T3.7}, we have $\mathfrak{p}_{_{P}}^{\triangleleft}$
$\ll_{_{\text{\textsf{d}}\mathsf{G}P}}$ $A.$ So $A/\mathfrak{p}_{_{P}%
}^{\triangleleft}$ is a $\mathsf{dG}P$-algebra.

Let $A/I$ be a $\mathsf{dG}P$-algebra. Then $I$ $\ll_{_{\text{\textsf{d}%
}\mathsf{G}P}}$ $A.$ As $\ll_{_{\text{\textsf{d}}\mathsf{G}P}}=$ $\ll_{_{P}%
}^{\triangleleft},$ we have $I$ $\ll_{_{P}}^{\triangleleft}$ $A.$ Then it
follows from Corollary \ref{C4.2n}(ii) that $\mathfrak{p}_{_{P}}%
^{\triangleleft}\subseteq I.$

The existence of a descending transfinite $\ll_{_{P}}$-series of ideals from
$A$ to $\mathfrak{p}_{_{P}}^{\triangleleft}$ also follows from Corollary
\ref{C4.2n}. Part (ii) follows from (i).\bigskip
\end{proof}

The following result is an analogue of Proposition \ref{P3.3P} for
$\mathsf{dG}P$ and \textsf{d}$\mathsf{NG}P$ C*-properties.

\begin{proposition}
\label{P3.5}If $P$ is lower stable$,$\emph{\ }the \emph{C*}-properties
$\mathsf{dG}P$ and $\mathsf{dNG}P$ are extension stable.
\end{proposition}

\begin{proof}
For $A\in\mathfrak{A},$ let $I\in$ Id$_{A}$ and the quotient $A/I$ be
$\mathsf{dG}P$-algebras. It follows from Theorem \ref{T3.2}(i) that
$\mathfrak{p}_{_{P}}^{\triangleleft}\subseteq I.$ Let $\mathfrak{p}_{_{P}%
}^{\triangleleft}\neq\{0\}.$ As $I\in$ $\mathsf{dG}P$ and $\mathfrak{p}_{_{P}%
}^{\triangleleft}$ is an ideal of $I,$ it follows that there is an ideal
$J\subsetneqq\mathfrak{p}_{_{P}}^{\triangleleft}$ such that $\mathfrak{p}%
_{_{P}}^{\triangleleft}/J\in P,$ i.e., $J$ is a $\ll_{_{P}}$-predecessor of
$\mathfrak{p}_{_{P}}^{\triangleleft}$ which contradicts Corollary \ref{C4.2n}.
Thus $\mathfrak{p}_{_{P}}^{\triangleleft}=\{0\}$ and $A\in$ $\mathsf{dG}P$ by
Theorem \ref{T3.2}. So the C*-property $\mathsf{dG}P$ is extension
stable.\smallskip

Let $A\in\mathfrak{A}$ and let its ideal $J$ and the quotient $A/J$ be
\textsf{d}$\mathsf{NG}P$-algebras. Suppose that $\mathfrak{p}_{_{P}%
}^{\triangleleft}\neq A.$ As $J$ is \textsf{d}$\mathsf{NG}P$-algebra,
$J\subseteq\mathfrak{p}_{_{P}}^{\triangleleft}$ and $A/\mathfrak{p}_{_{P}%
}^{\triangleleft}$ is a $\mathsf{dG}P$-algebra by Theorem \ref{T3.2}. Hence
$A/\mathfrak{p}_{_{P}}^{\triangleleft}$ has a quotient $(A/\mathfrak{p}_{_{P}%
}^{\triangleleft})/\widehat{I}\neq\{0\}$ which is a $P$-algebra. Therefore
there is $I\in$ Id$_{A}$ such that $\mathfrak{p}_{_{P}}^{\triangleleft
}\subsetneqq I\subsetneqq A,$ $\widehat{I}=I/\mathfrak{p}_{_{P}}$ and
$A/I\approx(A/\mathfrak{p}_{_{P}}^{\triangleleft})/\widehat{I}$ is a $P$-algebra.

Hence the quotient $(A/J)/(I/J)\approx A/I$ of $A/J$ is a $P$-algebra which
contradicts the assumption that $A/J$ is a \textsf{d}$\mathsf{NG}P$-algebra.
Thus $\mathfrak{p}_{_{P}}^{\triangleleft}=A.$ So $A$ is a \textsf{d}%
$\mathsf{NG}P$-algebra by Theorem \ref{T3.2}. Thus the C*-property
\textsf{d}$\mathsf{NG}P$ is extension stable.\bigskip
\end{proof}

As Example \ref{exDG} shows, the C*-property $\mathsf{dG}P$ can be much larger
than $P.$ Let us consider now the case when they coincide.\textbf{ }

\begin{theorem}
\label{dstable}A \emph{C*}-property $P$ coincides with $\mathsf{dG}P$ if and
only if $P$ is lower stable\emph{,} extension stable and satisfies the
following condition\emph{: }for each $A\in\mathfrak{A}$ and%
\begin{equation}
\text{for any family of ideals }\{I_{t}\}_{t\in T}\text{ with }A/I_{t}\in
P,\text{ the algebra }A/(\cap_{t\in T}I_{t})\in P. \label{coind}%
\end{equation}

\end{theorem}

\begin{proof}
Let $P=\mathsf{dG}P.$ Then $P$ is lower and extension stable by Lemma L3.7 and
Proposition \ref{P3.5}. For $A\in\mathfrak{A},$ let $A/I_{t}\in P$, $t\in T$.
Set $I=\cap_{t}I_{t}$. Let us prove that $A/I\in\mathsf{dG}P$. By Definition
\ref{D3.2}, we have to show that, for each $J\in$ Id$_{A}$, $I\subsetneqq J,$
there is $K\in$ Id$_{A}$ such that $I\subseteq K\subsetneqq J$ and $J/K\in P.$
As $I\subsetneqq J,$ there is $t$ with $K:=J\cap I_{t}\neq J$. Then $\{0\}\neq
J/K\overset{(\ref{4.01})}{\approx}(J+I_{t})/I_{t}$. Since $(J+I_{t})/I_{t}$ is
an ideal of the $P$-algebra $A/I_{t}$, it belongs to $P$, as $P$ is lower
stable. So $J/K\in P$ and $A/I$ is a $\mathsf{dG}P$-algebra. Thus it is a $P$-algebra.

Conversely, assume that $P$ is lower and extension stable and (\ref{coind})
holds. For $A\in\mathsf{dG}P$, let $W$ be the set of all ideals $I$ with
$A/I\in P$. By (\ref{coind}), the ideal $I_{0}:=\cap_{I\in W}I\in W$ and is
the smallest element of $W$. Suppose that $I_{0}\neq\{0\}$. As $A\in
\mathsf{dG}P$, there is $J\subsetneqq I_{0}$ with $I_{0}/J\in P$. Since
$I_{0}/J$ is an ideal of $A/J$ and $(A/J)/(I_{0}/J)\overset{(\ref{4.01}%
)}{\approx}A/I_{0}\in P$, it follows from the extension stability of $P$ that
$A/J\in P$. But this contradicts the minimality of $I_{0}$. So $I_{0}=\{0\}$
and $A\in P$. Thus $\mathsf{dG}P\subseteq G$; the converse inclusion follows
from Proposition \ref{P3.7}.
\end{proof}

\begin{remark}
\emph{The proof shows that in (\ref{coind}) one could consider only linearly
ordered families of ideals. In this case Zorn's Lemma gives that there is a
minimal ideal }$I_{0}\in W$\emph{ which is sufficient for further arguments.}
\end{remark}

\begin{corollary}
\label{C3.4}Let $P$ be a lower stable \emph{C*}-property. The following
conditions are equivalent.\smallskip

\emph{(i) \ }$\mathsf{dG}P=P;\smallskip$

\emph{(ii) \ \ }$\ll_{_{\mathsf{dG}P}}=$ $\ll_{_{P}}=$ $\ll_{_{P}%
}^{\triangleleft}$ in \emph{Id}$_{A}$ for each $A\in\mathfrak{A};\smallskip$

\emph{(iii) }$\ll_{_{P}}$ is a dual $\mathbf{R}$-order in \emph{Id}$_{A}$ for
each $A\in\mathfrak{A};$
\end{corollary}

\begin{proof}
(i) $\Rightarrow$ (ii). By Theorem \ref{T3.7}(i), $\ll_{_{\text{\textsf{d}%
}\mathsf{G}P}}=$ $\ll_{_{P}}^{\triangleleft}$ for each $A\in\mathfrak{A.}$ If
$\mathsf{dG}P=P$ then $\ll_{_{\text{\textsf{d}}\mathsf{G}P}}=$ $\ll_{_{P}},$
so that $\ll_{_{\text{\textsf{d}}\mathsf{G}P}}=$ $\ll_{_{P}}=$ $\ll_{_{P}%
}^{\triangleleft}.\smallskip$

(ii) $\Leftrightarrow$ (iii) follows from Theorem \ref{inf}(ii).\smallskip

(ii) $\Rightarrow$ (i). If $A\in$ $\mathsf{dG}P$ then $\{0\}\ll
_{_{\text{\textsf{d}}\mathsf{G}P}}A.$ By (ii), $\{0\}\ll_{_{P}}A.$ So $A\in
P.$ Thus $\mathsf{dG}P=P.$
\end{proof}

\begin{remark}
\label{R3.1}\emph{Let }$P$\emph{ be a lower and upper stable C*-property.
Then}

\emph{(i) For each} $A\in\mathfrak{A,}$ $A/\mathfrak{r}_{_{P}}^{\triangleright
}$ \emph{is a }$\mathsf{NG}P$\emph{-algebra}, $A/\mathfrak{p}_{_{P}%
}^{\triangleleft}$ \emph{is a }$\mathsf{dG}P$\emph{-algebra},%
\begin{equation}
\mathfrak{r}_{_{P}}^{\triangleright}=\mathfrak{r}_{_{\mathsf{G}P}%
}=\mathfrak{p}_{_{\mathsf{NG}P}}\text{ \emph{is a }}\mathsf{G}%
P\text{\emph{-ideal and} }\mathfrak{p}_{_{P}}^{\triangleleft}=\mathfrak{p}%
_{_{\mathsf{dG}P}}=\mathfrak{r}_{_{\mathsf{dNG}P}}\text{ \emph{is a }%
}\mathsf{dNG}P\text{-\emph{ideal}.} \label{3.5}%
\end{equation}

\emph{(ii) By Proposition \ref{P3.4}, }$\mathsf{G}P$ \emph{is a }lower and
upper stable \emph{C*}-property\emph{. By Lemmas \ref{L3.6} and \ref{L3.7},
}$\mathsf{NG}P$ \emph{and }$\mathsf{dG}P$ \emph{are }lower stable\emph{
C*-properties, and }$\mathsf{dNG}P$ \emph{is} an upper stable
\emph{C*-property}.

\emph{However, the C*-properties }$\mathsf{NG}P$ \emph{and }$\mathsf{dG}P$
\emph{are not always} upper stable\emph{ and }$\mathsf{dNG}P$ \emph{need not
be} lower stable. \emph{Indeed, let }$P=CCR$\emph{ be the C*-property that
consists of all }$CCR$\emph{ algebras. Then\smallskip\ }

$1)$ \emph{the C*-property }$\mathsf{NG}P$ \emph{consists of} $NGCR$\emph{
calgebras. It is not upper stable, since (see 4.7.4 c) \cite{D}) some
quotients of }$NGCR$\emph{-algebras may not be }$NGCR$%
\emph{-algebras;\smallskip}

$2)$ \emph{the C*-property }$\mathsf{dG}P$\emph{ consists of }$\mathsf{d}%
GCR$\emph{ algebras. It is not upper stable: in Remark \ref{R1} 1) below we
consider an example of a }$\mathsf{d}$\emph{GCR algebra that has quotients
which are not }$\mathsf{d}$\emph{GCR algebras;\smallskip}

$3)$ \emph{the C*-property }$\mathsf{dNG}P$\emph{ consists of }$\mathsf{d}%
NGCR$\emph{ algebras. It is not lower stable: in Remark \ref{R1} 2) below we
consider an example of a }$\mathsf{d}NGCR$\emph{ algebra with an ideal which
is not a }$\mathsf{d}NGCR$\emph{ algebra.\smallskip}

\emph{(iii) The algebras in }$\mathsf{dNG}P$\emph{,} \emph{for the C*-property
}$P=Comm$ \emph{of all commutative algebras, can be described as algebras
without multiplicative functionals. Indeed, if }$f$\emph{ is a multiplicative
functional, then }$A/\ker(f)$\emph{ is commutative, so }$A$\emph{ is not}
$\mathsf{dNG}P$. \emph{Conversely, if} $A\notin$ $\mathsf{dNG}P$ \emph{then
there is }$J\in\text{ }$\emph{Id}$_{A}$\emph{ with commutative }$A/J$\emph{.
Then any multiplicative functional on }$A/J$\emph{ defines a multiplicative
functional on} $A$.
\end{remark}

\subsection{The closure operators $P\mapsto\mathsf{R}(P)$ and $P\mapsto
\mathsf{G}_{_{\Pi}}(P)$ on upper stable properties}

Let $\Pi(A)$ be the set of all equivalence classes of non-zero irreducible
*-representations of $A\in\mathfrak{A}.$ For $I\in$ Id$_{A},$ let $p$:
$A\rightarrow A/I$ be the quotient map. Then (see Corollary 1.8.3, Proposition
2.10.4, Lemma 2.10.3 \cite{D})%
\begin{align}
\text{each }\pi &  \in\Pi(I)\text{ on }\mathcal{H}\text{ uniquely extends to
}\pi^{\prime}\in\Pi(A)\text{ on }\mathcal{H}\text{ and }\pi=\pi^{\prime}%
|_{I};\label{4.13}\\
\text{for }\pi &  \in\Pi(A\mathcal{)},\text{ either }I\subseteq\ker\pi\text{
or }\pi|_{I}\in\Pi(I),\text{ and }\pi(A)\approx A/\ker\pi;\label{4.13e}\\
\text{if }\pi &  \in\Pi(A/I)\text{ then }\pi^{\prime}=\pi\circ p\in\Pi(A).
\label{4.7}%
\end{align}

\begin{definition}
\label{D4.1}Let $P$ be a \emph{C*}-property in $\mathfrak{A}.$\emph{\ }We call
$A\in\mathfrak{A}$ an \textbf{$\mathsf{R}$}$(P)$\textbf{-algebra} if either
$A=\{0\},$ or $\pi(A)\in P$ for each $\pi\in\Pi(A).$
\end{definition}

Clearly, the class of all $\mathsf{R}(P)$-algebras is a C*-property.

Recall that $\mathcal{P}_{\text{up}}$ is the complete lattice of all upper
stable C*-properties in $\mathfrak{A}$. For $P\in\mathcal{P}_{\text{up}},$
consider the C*-property $P_{\mathsf{G}P}=\{A\in\mathfrak{A}$: there is
$\pi\in\Pi(A)$ with $\ker\pi\in\mathsf{G}P\}.$

\begin{theorem}
\label{T3.6}\emph{(i) }For any \emph{C*}-property $P$ in $\mathfrak{A},$ the
\emph{C*}-property $\mathsf{R}(P)$ is upper stable and $\mathsf{R}%
(P)\subseteq\mathsf{G}(\mathsf{R}(P)).\smallskip$

\emph{(ii) }The map $P\rightarrow\mathsf{R}(P)$ is a closure operator
\emph{(}see \emph{(\ref{2.6}))} on $\mathcal{P}_{\text{\emph{up}}%
}(\mathfrak{A})$\emph{, }that is\emph{,}%
\begin{equation}
P\subseteq\mathsf{R}(P)=\mathsf{R}(\mathsf{R}(P))\text{ and }\mathsf{R}%
(P)\subseteq\mathsf{R}(P_{1})\text{ for }P\subseteq P_{1}\text{ in
}\mathcal{P}_{\text{\emph{up}}}. \label{3.31}%
\end{equation}

If $P\subseteq P_{1}\subseteq\mathsf{R}(P)$ then $\mathsf{R}(P)=\mathsf{R}%
(P_{1})$. Moreover\emph{, }for $P\in\mathcal{P}_{\text{\emph{up}}},$
\begin{align*}
\mathsf{R}(\mathsf{G}P)\cap P_{\mathsf{G}P}  &  =\mathsf{G}P\cap
P_{\mathsf{G}P}\subseteq\mathsf{G}(\mathsf{R}(P))\cap P_{\mathsf{G}P}\text{
and}\\
\mathsf{G}P\cup\mathsf{R}(P)  &  \subseteq\mathsf{G}(\mathsf{R}(P))\cap
\mathsf{R}(\mathsf{G}P).
\end{align*}

\emph{(iii) }Let $P$ be lower stable\emph{. }Then\emph{ }$\mathsf{R}(P)$ is
lower stable\emph{,} $\mathsf{G}(\mathsf{R}(P))$ is lower and upper
stable\emph{, }and $\mathsf{R}(P)\subseteq\mathsf{dG}(\mathsf{R}%
(P))\subseteq\mathsf{dG}P.$
\end{theorem}

\begin{proof}
(i) Let $A\in\mathsf{R}(P)$ and $I\in$ Id$_{A}.$ Let $p$: $A\rightarrow A/I$
be the standard epimorphism. By (\ref{4.7}), $\pi^{\prime}=\pi\circ p\in
\Pi(A)$ for each $\pi\in\Pi(A/I)$. So $\pi(A/I)=\pi^{\prime}(A)\in P$. Thus
$A/I$ is a $\mathsf{R}(P)$-algebra. Hence the C*-property $\mathsf{R}(P)$ is
always upper stable. So, by Proposition \ref{L3.1},\emph{ }$\mathsf{R}%
(P)\subseteq\mathsf{G}(\mathsf{R}(P)).$

(ii) Let $P$ be upper stable and $A\in P.$ By (\ref{4.13e}), $\pi(A)\approx
A/\ker\pi$ for each $\pi\in\Pi(A).$ As $P$ is upper stable, $A/\ker\pi\in P.$
Hence $\pi(A)\in P.$ So $A\in\mathsf{R}(P).$ Thus $P\subseteq\mathsf{R}(P).$

As $\mathsf{R}(P)$ is upper stable, $\mathsf{R}(P)\subseteq\mathsf{R}%
(\mathsf{R}(P))$ by the above. Let $A\in\mathsf{R}(\mathsf{R}(P)).$ For each
$\pi\in\Pi(A)$, $\pi(A)\in\mathsf{R}(P).$ Let $id$ be the identity map of
$\pi(A)$ on itself. Then $id\in\Pi(\pi(A)).$ As $\pi(A)$ is a $\mathsf{R}%
(P)$-algebra, $\pi(A)=id(\pi(A))\in P.$ Thus $A\in\mathsf{R}(P).$ So
$\mathsf{R}(\mathsf{R}(P))\subseteq\mathsf{R}(P).$ Hence $\mathsf{R}%
(\mathsf{R}(P))=\mathsf{R}(P).\smallskip$

Let $P\subseteq P_{1}$ and $A\in\mathsf{R}(P).$ For each $\pi\in\Pi(A),$
$\pi(A)\in P\subseteq P_{1}.$ Hence $A\in\mathsf{R}(P_{1}).$ So $\mathsf{R}%
(P)\subseteq\mathsf{R}(P_{1}).$ Thus the map $P\mapsto\mathsf{R}(P)$ is a
closure operator in $\mathcal{P}_{\text{up}}$ (see (\ref{2.6})). By Lemma
\ref{L2.1}, $P\subseteq P_{1}\subseteq\mathsf{R}(P)$ implies $\mathsf{R}%
(P)=\mathsf{R}(P_{1})$.

Let $P\in\mathcal{P}_{\text{up}}$. For $A\in\mathsf{R}(\mathsf{G}P)\cap
P_{\mathsf{G}P},$ let $\pi\in\Pi(A)$ be such that $\ker\pi\in$ \textsf{G}$P.$
By Definition \ref{D4.1}, $A/\ker\pi\in$ \textsf{G}$P.$ By Proposition
\ref{P3.3P}, the C*-property \textsf{G}$P$ is extension stable (Definition
\ref{D3.3}). Thus $A\in$ \textsf{G}$P.$ So $\mathsf{R}(\mathsf{G}P)\cap
P_{\mathsf{G}P}\subseteq$ \textsf{G}$P\cap P_{\mathsf{G}P}.$ As $\mathsf{G}%
P\subseteq\mathsf{R}(\mathsf{G}P)$ and $P\subseteq\mathsf{R(}P)$ by
(\ref{3.31}), we have $\mathsf{R}(\mathsf{G}P)\cap P_{\mathsf{G}P}%
=\mathsf{G}P\cap P_{\mathsf{G}P}\subseteq\mathsf{G}(\mathsf{R(}P))\cap
P_{\mathsf{G}P}.$

By Proposition \ref{L3.1}, $P\subseteq\mathsf{G}P.$ So $\mathsf{R}%
(P)\subseteq\mathsf{R}(\mathsf{G}P)$ by (\ref{3.31}). On the other hand,
$\mathsf{R}(P)\subseteq\mathsf{G}(\mathsf{R}(P))$ by Proposition \ref{L3.1}.
So $\mathsf{R}(P)\subseteq\mathsf{G}(\mathsf{R}(P))\cap\mathsf{R}%
(\mathsf{G}P).$ Similarly, we get $\mathsf{G}P\subseteq\mathsf{G}%
(\mathsf{R}(P))\cap\mathsf{R}(\mathsf{G}P)$. So $\mathsf{G}P\cup
\mathsf{R}(P)\subseteq\mathsf{G}(\mathsf{R}(P))\cap\mathsf{R}(\mathsf{G}P).$

(iii) Let $P$ be lower stable, $A$ be a $\mathsf{R}(P)$-algebra and $I\in$
Id$_{A}.$ By (\ref{4.13}), each $\pi\in\Pi(I)$ extends to some $\pi^{\prime
}\in\Pi(A).$ So $\pi^{\prime}(A)\in P$. As $\pi(I)$ is an ideal of
$\pi^{\prime}(A)$ and the C*-property $P$ is lower stable, $\pi(I)\in P$. Thus
$I$ is a $\mathsf{R}(P)$-algebra. So the C*-property $\mathsf{R}(P)$ is lower stable.

Let $A\in\mathsf{R}(P).$ As $\mathsf{R}(P)$ is lower stable, each $\{0\}\neq
I\in$ Id$_{A}$ belongs to $\mathsf{R}(P).$ Hence, for every $\pi\in\Pi(I),$
$I/\ker\pi\approx\pi(I)\in P$ by (\ref{4.13e}). Thus $A\in$ $\mathsf{dG}P.$ So
$\mathsf{R}(P)\subseteq$ $\mathsf{dG}P.$ Hence, as $\mathsf{R}(P)$ is lower
stable$,$ we have from Proposition \ref{P3.7} that $\mathsf{R}(P)\subseteq
\mathsf{dG(R}(P)\mathsf{)}\subseteq\mathsf{dG}P.$

If $P$ is lower stable, $\mathsf{R}(P)$ is a lower and upper stable
C*-property (see above). Hence, by Proposition \ref{P3.4}, \textsf{G}%
$(\mathsf{R}(P))$ is lower and upper stable.
\end{proof}

\begin{corollary}
\label{C3.12}Let $P$ be a lower and an upper stable \emph{C*}-property.
Then\smallskip

\emph{(i) \ \ }$\mathsf{dG}(\mathsf{R}(P))=\mathsf{dG}P.\smallskip$

\emph{(ii)} $\ A\in\mathsf{dG}P$\emph{ }if and only if\emph{, }for each
$\{0\}\neq I\in$ \emph{Id}$_{A},$ there is $\pi^{\prime}\in\Pi(I)$ with
$\pi^{\prime}(I)\in P.\smallskip$

\emph{(iii) }The relation $\ll_{_{P}}^{\triangleleft}=$ $\ll_{_{\mathsf{dG}P}%
}=$ $\ll_{_{\mathsf{dG}(\mathsf{R}(P))}}=$ $\ll_{_{\mathsf{R}(P)}%
}^{\triangleleft}$ is a dual $\mathbf{R}$-order in each \emph{Id}$_{A}%
$\emph{,} so that
\[
\mathfrak{p}_{_{P}}^{\triangleleft}=\mathfrak{p}_{_{\mathsf{dG}P}%
}=\mathfrak{p}_{_{\mathsf{dG}(\mathsf{R}(P))}}=\mathfrak{p}_{_{\mathsf{R}(P)}%
}^{\triangleleft}.
\]

\end{corollary}

\begin{proof}
(i) It follows from Theorem \ref{T3.6} that $\mathsf{dG}(\mathsf{R}%
(P))\subseteq\mathsf{dG}P,$ $P\subseteq\mathsf{R}(P)$ and $\mathsf{R}(P)$ is
lower stable. Thus $\mathsf{dG}P\subseteq\mathsf{dG}(\mathsf{R}(P))$ by
Proposition \ref{P3.7}. So $\mathsf{dG}(\mathsf{R}(P))=\mathsf{dG}P.$

(ii) Let $A\in$ $\mathsf{dG}P$ and $\{0\}\neq I\in$ Id$_{A}.$ By (i), $A\in$
$\mathsf{dG}(\mathsf{R}(P)).$ Then $I/J\in\mathsf{R}(P)$ for some
$J\subsetneqq I.$ Hence $\pi(I/J)\in P$ for all $\pi\in\Pi(I/J).$ By
(\ref{4.7}), there is $\pi^{\prime}\in\Pi(I)$ such that $\pi^{\prime}%
(I)=\pi(I/J)\in P.$

Conversely, if $\pi^{\prime}(I)\in P$ for each $I\in$ Id$_{A}$ and some
$\pi^{\prime}\in\Pi(I),$ then $I/\ker\pi^{\prime}\approx\pi^{\prime}(I)\in P.$
So $A\in$ $\mathsf{dG}P$.

Part (iii) follows from (i) and Theorem \ref{T3.7}.\bigskip
\end{proof}

We will show now that all closure operators \emph{$\mathsf{R}$}, $\mathsf{G}$,
$\mathsf{G\circ}$\emph{$\mathsf{R}$}, \emph{$\mathsf{R\circ}$}$\mathsf{G,}$
$\mathsf{dG}$ are different by applying them to the C*-property $C$ consisting
of C*-algebras isomorphic to $C(H)$ for arbitrary $H.$

\begin{example}
\label{3.24}\emph{(i)} \emph{$\mathsf{R\neq}$ }$\mathsf{G\circ}$%
\emph{$\mathsf{R\neq}$ }$\mathsf{G.}$\emph{ Indeed, by definition,
$\mathsf{R}(C)=CCR,$ $\mathsf{G}(C)=Sc$ -- the class of scattered algebras
(Example \ref{exG}) and $\mathsf{G}(\mathsf{R}(C))=\mathsf{G}(CCR)=GCR$. As
}$C(0,1)\in CCR$ \emph{and} $C(0,1)\notin Sc,$ \emph{we have}
\emph{$\mathsf{R}(C)$ $\mathsf{\neq G}(C)\mathsf{.}$ So $\mathsf{R}$
$\mathsf{\neq G.}$ Clearly, $GCR\neq CCR,$ so that $\mathsf{R\neq}$
}$\mathsf{G\circ}$\emph{$\mathsf{R.}$ As }$B(H)\in GCR$ \emph{and }$B(H)\in
Sc,$ $\mathsf{G\circ}$\emph{$\mathsf{R\neq}$ }$\mathsf{G.}$

\emph{(ii) $\mathsf{R\neq}$ $\mathsf{R}$}$\circ\mathsf{G}$ \emph{$\mathsf{\neq
}$ }$\mathsf{G.}$ \emph{Indeed,} $\mathsf{G}C\cup\mathsf{R}(C)\subseteq
\mathsf{R}(\mathsf{G}C)$ \emph{by Corollary \ref{C3.12}, and} $\mathsf{R}%
(C)\neq\mathsf{G}C$\emph{ by (i).}

\emph{(iii) $\mathsf{R}$}$\circ\mathsf{G}$ $\mathsf{\neq}$ $\mathsf{G\circ}%
$\emph{$\mathsf{R.}$ By (i), $\mathsf{R}(\mathsf{G}(C))=\mathsf{R}(Sc)$ and
}$\mathsf{G(}$\emph{$\mathsf{R}(C))=$ }$\mathsf{G(}$$CCR)=GCR.$ \emph{Thus it
suffices to find a $GCR$-algebra that does not belong to $\mathsf{R}(Sc)$. The
Toeplitz algebra }$T$ \emph{(the C*-algebra generated by the unilateral shift
$U$) is what we need. It is a GCR-algebra. As it is irreducible, the identity
representation }$id\in\Pi(T)\emph{,}$ \emph{but }$id(T)=T\notin$ \emph{$Sc$,
since the spectrum $\sigma(U)$ is uncountable}. \emph{So} $T\notin$
\emph{$\mathsf{R}(Sc).$}

\emph{(iv) If }$P$ \emph{is lower and upper stable, then }$\mathsf{R}%
(P)\subseteq\mathsf{dG}(\mathsf{R}(P))=\mathsf{dG}P$ \emph{by Theorem
\ref{T3.6} and Corollary \ref{C3.12}}. \emph{Let us show that,} in general,
$\mathsf{R}(P)\neq\mathsf{dG}P.$ \emph{Let }$P=C.$ \emph{Then }$\mathsf{R}%
(C)\neq\mathsf{dG}C.$ \emph{Indeed, }$A=C(H)+\mathbb{C}\mathbf{1}_{H}$\emph{
is not a }$\mathsf{R}(C)$\emph{-algebra (}$CCR$-\emph{algebra), as its
identity representation }$id$\emph{ is irreducible and }$id(A)=A\notin
C.$\emph{ On the other hand, }$A$\emph{ is a }$\mathsf{dG}C$\emph{-algebra, as
each its non-zero ideal }$(C(H)$\emph{ and }$A)$\emph{ has a quotient from
that belongs to }$C.$
\end{example}

If in Definition \ref{D3.1} of the C*-property $\mathsf{G}P$ one only
considers primitive ideals instead of all ideals, then one obtains a new
important closure operator $\mathsf{G}_{_{\Pi}}$ on $P_{\text{up}}(A).$

\begin{definition}
\label{D3.4}Let $P$ be a \emph{C*}-property in $\mathfrak{A.}$ We call
$A\in\mathfrak{A}$\emph{\ }a $\mathsf{G}_{_{\Pi}}P$-algebra if either
$A=\{0\},$ or each quotient $A/\ker\pi,$ $\pi\in\Pi(A),$ has a non-zero
$P$-ideal\emph{,} that is\emph{,} $\pi(A)$ has a non-zero $P$-ideal.
\end{definition}

We will show now that $\mathsf{G}_{_{\Pi}}$ is a closure operator on the class
$\mathcal{P}_{\text{up}}$ of upper stable C*-properties.

\begin{theorem}
\label{T3.4}\emph{(i) }For any \emph{C*}-property $P,$ the \emph{C*}-property
$\mathsf{G}_{_{\Pi}}P$ is upper stable and%
\[
\mathsf{G}(\mathsf{R}(P))\cup\mathsf{R}(\mathsf{G}P)\subseteq\mathsf{G}%
_{_{\Pi}}P.
\]

\emph{(ii) \ }The map $P\rightarrow\mathsf{G}_{_{\Pi}}P$ is a closure operator
\emph{(}see \emph{(\ref{2.6}))} on $\mathcal{P}_{\text{\emph{up}}}$\emph{,
}i.e.\emph{,}%
\[
P\subseteq\mathsf{G}_{_{\Pi}}P=\mathsf{G}_{_{\Pi}}(\mathsf{G}_{_{\Pi}}P)\text{
and }\mathsf{G}_{_{\Pi}}P\subseteq\mathsf{G}_{_{\Pi}}P_{1}\text{ for
}P\subseteq P_{1}\text{ in }\mathcal{P}_{\text{\emph{up}}}.
\]

\emph{(iii)\ }The \emph{C*}-property $\mathsf{G}_{_{\Pi}}P$ is extension
stable \emph{(}Definition \emph{\ref{D3.3}).\smallskip}
\end{theorem}

\begin{proof}
(i) Let $A\in$ \textsf{G}$_{_{\Pi}}P,$ $I\in$ Id$_{A}$ and $p$: $A\rightarrow
A/I.$ For $\pi\in\Pi(A/I),$ $\pi^{\prime}=\pi\circ p\in\Pi(A)$ and
$\pi^{\prime}(A)=\pi(A/I)$ by (\ref{4.7}). As $A\in$ \textsf{G}$_{_{\Pi}}P,$
$\pi^{\prime}(A)$ has a non-zero $P$-ideal. Hence $\pi(A/I)$ has a non-zero
$P$-ideal. So $A/I\in$ \textsf{G}$_{_{\Pi}}P.$ Thus the C*-property
\textsf{G}$_{_{\Pi}}P$ is upper stable.

Let $A$ be a \textsf{G}$(\mathsf{R}(P))$-algebra and $\pi\in\Pi(A).$ As
$\pi(A)\approx A/\ker\pi,$ it follows from Definition \ref{D3.1} that the
operator algebra $\pi(A)$ on a Hilbert space $\mathcal{H}_{\pi}$ has an
$\mathsf{R}(P)$-ideal $I_{\pi}\neq\{0\}$. By (\ref{4.13e}), the identity
representation $id$ of $I_{\pi}$ on $\mathcal{H}_{\pi}$ is irreducible. So
$id\in\Pi(I_{\pi}).$ Hence, by Definition \ref{D4.1}, $I_{\pi}=id(I_{\pi})\in
P.$ Thus $A\in$ \textsf{G}$_{_{\Pi}}P$ by Definition \ref{D3.4}. So
\textsf{G}$(\mathsf{R}(P))\subseteq\mathsf{G}_{_{\Pi}}P$.

Let $A\in\mathsf{R}(\mathsf{G}(P))$ and $\pi\in\Pi(A).$ By Definition
\ref{D4.1}, $\pi(A)\in\mathsf{G}P.$ By Definition \ref{D3.1}, each non-zero
quotient of $\pi(A)$ has a non-zero $P$-ideal. In particular, $\pi(A)$ has a
non-zero $P$-ideal. So, by Definition \ref{D3.4}, $A\in\mathsf{G}_{_{\Pi}}P.$
Thus $\mathsf{R}(\mathsf{G}(P))\subseteq\mathsf{G}_{_{\Pi}}P$ which completes
the proof of (i).\smallskip

(ii) By (i) and Theorem \ref{T3.6}, $P\subseteq\mathsf{R}(P)\subseteq$
\textsf{G}$(\mathsf{R}(P))\subseteq$ \textsf{G}$_{_{\Pi}}P.$ So \textsf{G}%
$_{_{\Pi}}P\subseteq$ \textsf{G}$_{_{\Pi}}($\textsf{G}$_{_{\Pi}}P).$

Let $A\in$ \textsf{G}$_{_{\Pi}}($\textsf{G}$_{_{\Pi}}P)$ and $\pi\in\Pi(A).$
By Definition \ref{D3.4}, $\pi(A)$ has a non-zero \textsf{G}$_{_{\Pi}}P$-ideal
$I.$ By (\ref{4.13e}), the identity representation $id$ of $I$ on
$\mathcal{H}_{\pi}$ belongs to $\Pi(I).$ As $I\in$ \textsf{G}$_{_{\Pi}}P,$ it
follows from Definition \ref{D3.4} that $I=id(I)$ has a non-zero $P$-ideal
which is also an ideal of $\pi(A).$ Hence $A\in$ \textsf{G}$_{_{\Pi}}P.$ Thus
\textsf{G}$_{_{\Pi}}($\textsf{G}$_{_{\Pi}}P)\subseteq$ \textsf{G}$_{_{\Pi}}P.$
So \textsf{G}$_{_{\Pi}}P=$ \textsf{G}$_{_{\Pi}}($\textsf{G}$_{_{\Pi}}P).$

Let $P\subseteq P_{1}$ and $A\in$ \textsf{G}$_{_{\Pi}}P.$ For each $\pi\in
\Pi(A),$ $\pi(A)$ has a $P$-ideal $I\neq\{0\}.$ Then $I$ is a $P_{1}$-ideal of
$\pi(A).$ So $A\in$ \textsf{G}$_{_{\Pi}}P_{1}.$ Thus \textsf{G}$_{_{\Pi}%
}P\subseteq$ \textsf{G}$_{_{\Pi}}P_{1}.$ Hence $P\rightarrow$ \textsf{G}%
$_{_{\Pi}}P$ is a closure operator on $\mathcal{P}_{\text{up}}(\mathfrak{A}).$

(iii) Let $A\in\mathfrak{A}$ and $I\in$ Id$_{A}.$ Let $I,$ $A/I\in$
\textsf{G}$_{_{\Pi}}P,$ let $\pi\in\Pi(A)$ and $p$: $A\rightarrow A/I.$ If
$\pi|_{I}\neq\{0\}$ then $\pi|_{I}\in\Pi(I)$ by (\ref{4.13e}). As
$I\in\mathsf{G}_{_{\Pi}}P,$ $\pi(I)$ contains a $P$-ideal which is a $P$-ideal
of $\pi(A).$

Let $\pi|_{I}=\{0\}.$ Then the map $\rho$: $x\in A/I\rightarrow\pi(p^{-1}(x))$
belongs to $\Pi(A/I)$ and $\rho(A/I)=\pi(A).$ As $A/I\in$ \textsf{G}$_{_{\Pi}%
}P,$ $\rho(I)$ contains a $P$-ideal which is also a $P$-ideal of $\pi(A).$
Thus $A\in$ \textsf{G}$_{_{\Pi}}P.$ So \textsf{G}$_{_{\Pi}}P$ is extension
stable.$\bigskip$
\end{proof}

Glimm's \cite{Gl} famous result states that GCR-algebras can be equivalently
defined as C*-algebras whose images in all irreducible representations contain
compact operators. Thus the equality
\begin{equation}
\mathsf{G}_{_{\Pi}}P=\mathsf{G}(\mathsf{R}(P))\text{ holds for }P=C\text{:
}\mathsf{G}_{_{\Pi}}C=\mathsf{G}(\mathsf{R}(C))=GCR. \label{quest}%
\end{equation}
It would be interesting to find conditions on a C*-property $P$ for which
equality (\ref{quest}) is valid:

\begin{problem}
\label{P2}\emph{When }$\mathsf{G}(\mathsf{R}(P))=\mathsf{R}(\mathsf{G}%
P))$\emph{ and when one of them,\ or both coincide with }$\mathsf{G}_{_{\Pi}%
}P?$
\end{problem}

We consider below some C*-properties for which (\ref{quest}) holds.

Let $P=Comm$ and $A\in\mathsf{G}_{_{\Pi}}P.$ For $\pi\in\Pi(A),$ $\pi(A)$ has
a commutative ideal. Hence $\dim\pi=1$. Repeating the argument of Example
\ref{exG} 3), we get $A\in P.$ So $\mathsf{G}_{_{\Pi}}P=\mathsf{G}%
(\mathsf{R}(P))=\mathsf{R}(\mathsf{G}(P))=P$.

Below we generalize this result. Let $A\in\mathfrak{A}$ and $I\in$ Id$_{A}.$
It follows from (\ref{4.13}) that there is a unique map $\theta_{I,A}$:
$\Pi(I)\mapsto\Pi(A)$ such that
\begin{equation}
\mathcal{H}_{\pi^{\prime}}=\mathcal{H}_{\pi}\text{ and }\pi^{\prime}|_{I}%
=\pi,\text{ where }\pi^{\prime}=\theta_{I,A}(\pi)\text{ for }\pi\in\Pi(I).
\label{3.23}%
\end{equation}

\begin{definition}
\label{D3.5}A map $F$ that maps each $A\in\mathfrak{A}$ into a subset $F(A)$
of $\Pi(A)$ is \textbf{compatable}\emph{ }if\smallskip

$1)$ $\theta_{I,A}(\pi)\in F(A)$ for all $A\in\mathfrak{A,}$ $I\in
\emph{Id}_{A}$ and $\pi\in F(I);\smallskip$

$2)$ $A\approx B$ implies that\emph{, }for some isomorphism $\varphi$\emph{:
}$B\rightarrow A,$ $\pi\circ\varphi\in F(B)$ for all $\pi\in F(A).$%
\textbf{\smallskip}

A compatable map $F$ is \textbf{an ideal map} if $\pi\circ\rho\in F(A)$ for
all $\rho\in\Pi(A)$ and $\pi\in F(\rho(A)).$
\end{definition}

The maps in $1)-4)$ below are ideal maps. For each $A\in\mathfrak{A,}$
let\smallskip

$1)$ $F_{n}(A)=\{\pi\in\Pi(A)$: $\dim\pi\leq n\}$;\smallskip

$2)$ $F_{\infty}(A)=\{\pi\in\Pi(A)$: $\dim\pi<\infty\}=\Pi_{f}(A)$ (see
\ref{3.30}))$;$\smallskip

$3)$ $F_{_{GCR}}(A)=\{\pi\in\Pi(A)$: $C(\mathcal{H}_{\pi})\subseteq
\operatorname{Im}\pi\};$\smallskip

$4)$ $F_{_{NGR}}(A)=\{\pi\in\Pi(A)$: $\operatorname{Im}\pi\cap C(\mathcal{H}%
_{\pi})=\{0\}\}$. \smallskip

$5)$ $F_{_{CCR}}(A)=\{\pi\in\Pi(A)$: $\operatorname{Im}\pi=C(\mathcal{H}_{\pi
})\}$

The map $F_{_{CCR}}$ is not compatable (take $A=B(\mathcal{H})$ and
$I=C(\mathcal{H}))$.

For a compatable map $F,$ let $r_{_{F}}$ be the class of all $A\in
\mathfrak{A}$ for which $F(A)$ is a separting set:%
\begin{equation}
r_{_{F}}=\{A\in\mathfrak{A}\text{:\emph{ }}\cap\{\ker\pi\text{: }\pi\in
F(A)\}=\{0\}\}. \label{3.19}%
\end{equation}
It is a C*-property: if $A\in r_{_{F}}$ and $A\approx B,$ then $B\in r_{_{F}%
},$ since for some isomorphism $\varphi$\emph{: }$B\rightarrow A,$%
\begin{align*}
\cap\{\ker\rho &  \text{:}\rho\in F(B)\}\subseteq\cap\{\ker(\pi\circ
\varphi\text{: }\pi\in F(A)\}\\
&  =\{b\in B\text{: }\varphi(b)\in\cap\{\ker(\pi\text{: }\pi\in
F(A)\}\}=\{0\}.
\end{align*}
Clearly, $r_{_{F_{n}}}\subset r_{_{F_{\infty}}}=RFD\subset r(F_{_{GCR}})$ and
$GCR\subset r(F_{_{GCC}}).$

\begin{lemma}
\label{L3.12}Let $F$ be a compatable map on $\mathfrak{A}.$ If an irreducible
\emph{C*}-algebra $B$ of operators has a non-zero $r_{_{F}}$-ideal $I$\emph{,}
then $B\in r_{_{F}}.$
\end{lemma}

\begin{proof}
By Definition \ref{D3.5}, for each $\pi\in F(I),$ $\pi^{\prime}=\theta
_{I,A}(\pi)\in F(B).$ Set $J=\cap\{\ker\pi^{\prime}$: $\pi\in F(I)\}.$ Then
$J$ is an ideal of $B.$ For all $a\in I,$ $b\in J$ and $\pi\in F(I),$ we have
$ab\in I\cap J$ and $\pi(ab)=\pi^{\prime}(ab)=\pi^{\prime}(a)\pi^{\prime
}(b)=0,$ as $\pi^{\prime}(b)=0.$ As $\cap\{\ker\pi$: $\pi\in F(I)\}=\{0\}$ by
(\ref{3.19}), we have $ab=0.$ So $I\cdot J=\{0\}.$

As $B$ is irreducible, $\{0\}$ is a primitive ideal of $B.$ As $I\neq\{0\},$
it follows from Lemma 2.11.4 \cite{D} that $J=\{0\}.$ Hence%
\[
\cap\{\ker\rho\text{: }\rho\in F(B)\}\subseteq\cap\{\ker\pi^{\prime}\text{:
}\pi\in F(I)\}=J=\{0\}.
\]
So $B\in r_{_{F}}.$
\end{proof}

\begin{theorem}
\label{T3.10}Let $P=r_{_{F}}$ be the \emph{C*}-property for a compatable map
$F$ on $\mathfrak{A}.$ Then\smallskip

\emph{(i)\ \ }$\mathsf{G}_{_{\Pi}}r_{_{F}}=\mathsf{G}(\mathsf{R}(r_{_{F}%
}))=\mathsf{R}(r_{_{F}}).$\smallskip\ 

\emph{(ii) }If $F$ is an ideal map then $\mathsf{G}_{_{\Pi}}r_{_{F}%
}=\mathsf{G}(\mathsf{R}(r_{_{F}}))=\mathsf{R}(r_{_{F}})\subseteq r_{_{F}}.$

Moreover\emph{, }$\mathsf{G}_{_{\Pi}}r_{_{F}}=\mathsf{G}(\mathsf{R}(r_{_{F}%
}))=\mathsf{R}(r_{_{F}})=r_{_{F}}$ if and only if $r_{_{F}}$ is upper stable.
\end{theorem}

\begin{proof}
(i) For any C*-property $P,$ $\mathsf{G}(\mathsf{R}(P))\cup\mathsf{R}%
(\mathsf{G}P)\subseteq\mathsf{G}_{_{\Pi}}P$ by Theorem \ref{T3.4}, and
$\mathsf{R}(P)\subseteq\mathsf{G}(\mathsf{R}(P))$ by Theorem \ref{T3.6}. So
\begin{equation}
\mathsf{R}(P)\subseteq\mathsf{G}(\mathsf{R}(P))\subseteq\mathsf{G}_{_{\Pi}}P.
\label{3.33}%
\end{equation}

Let $P=r_{_{F}},$ $A\in\mathsf{G}_{_{\Pi}}P$ and $\pi\in\Pi(A).$ Then $\pi(A)$
is an irreducible operator algebra on $\mathcal{H}_{\pi}$ and, by Definition
\ref{D3.4}, it has a non-zero $P$-ideal. By Lemma \ref{L3.12}, $\pi(A)$ is a
$P$-algebra. So, by Definition \ref{D4.1}, $A\in\mathsf{R}(P).$ Thus
$\mathsf{G}_{_{\Pi}}P=\mathsf{G}(\mathsf{R}(P))=\mathsf{R}(P).$

(ii) Let $A\in\mathsf{G}_{_{\Pi}}P$ and $a\in A.$ Choose $\pi\in\Pi(A)$ with
$\pi(a)\neq0.$ Then $\pi(A)$ has an $P$-ideal $I\neq\{0\}$ by Definition
\ref{D3.4}. As the operator C*-algebra $\pi(A)$ is irreducible on
$\mathcal{H}_{\pi},$ we have from Lemma \ref{L3.12} that $\pi(A)\in P.$ By
(\ref{3.19}), there is $\tau\in F(\pi(A))$ such that $\tau(\pi(a))\neq0.$ As
$F$ is an ideal map, $\tau\circ\pi\in F(A)$ and $(\tau\circ\pi)(a)=\tau
(\pi(a))\neq0.$ So $A\in P$ by (\ref{3.19}). Thus $\mathsf{G}_{_{\Pi}}P\subset
P.$

If $P$ is upper stable then $P\subseteq\mathsf{G}_{_{\Pi}}P$ by Proposition
\ref{L3.1} and Theorem \ref{T3.4}. Thus $\mathsf{G}_{_{\Pi}}P=\mathsf{G}%
(\mathsf{R}(P))=\mathsf{R}(P)=P.$ Conversely, if this equality holds then $P$
is upper stable, as $\mathsf{G}_{_{\Pi}}P$ is upper stable by Theorem
\ref{T3.4}.\bigskip
\end{proof}

By Theorem \ref{T3.10}, for all maps $F$ in $1)-4)$ above, $\mathsf{G}_{_{\Pi
}}r_{_{F}}=\mathsf{G}(\mathsf{R}(r_{_{F}}))=\mathsf{R}(r_{_{F}})\subseteq
r_{_{F}}.$

Consider now the following C*-properties:%
\begin{align}
P_{n}  &  =\{A\in\mathfrak{A}\text{: }\dim\pi\leq n\text{ for }\pi\in
\Pi(A)\}\text{ for }1\leq n<\infty,\text{ and}\nonumber\\
P_{\infty}  &  =\{A\in\mathfrak{A}\text{: }\dim\pi<\infty\text{ for }\pi\in
\Pi(A)\}. \label{3.20}%
\end{align}

\begin{proposition}
\label{P4.6}\emph{(i) C*}-properties $P_{n},$ $1\leq n\leq\infty,$ are
lower\emph{, }upper and extension stable and%
\begin{equation}
\mathsf{G}_{_{\Pi}}P_{n}=\mathsf{G}(\mathsf{R}(P_{n}))=\mathsf{R}(P_{n})=P_{n}
\label{3.25}%
\end{equation}

\emph{(ii) }$Comm=P_{1}=r_{_{F_{1}}},$ so that $\mathsf{G}_{_{\Pi}}r_{_{F_{1}%
}}=\mathsf{G}(\mathsf{R}(r_{_{F_{1}}}))=\mathsf{R}(r_{_{F_{1}}})=r_{_{F_{1}}%
}.$

\emph{(iii) }The \emph{C*}-property $r_{_{F_{\infty}}}=RFD$ is lower and
extension stable\emph{, }but not upper stable\emph{, }so that%
\[
P_{\infty}\subsetneqq r_{_{F_{\infty}}}\text{ and }\mathsf{G}_{_{\Pi}%
}r_{_{F_{\infty}}}=\mathsf{G}(\mathsf{R}(r_{_{F_{\infty}}}))=\mathsf{R}%
(r_{_{F_{\infty}}})\subsetneqq r_{_{F_{\infty}}}.
\]

\end{proposition}

\begin{proof}
(i) By (\ref{4.13})-(\ref{4.7}),\ all $P_{n}$ are lower\emph{, }upper and
extension stable C*-properties.

By (\ref{3.33}) and Theorem \ref{T3.6}, $P\subseteq\mathsf{R}(P)\subseteq
\mathsf{G}(\mathsf{R}(P))\subseteq\mathsf{G}_{_{\Pi}}P$ for all upper stable C*-properties.

Let $A\in\mathsf{G}_{_{\Pi}}P_{n}$ and $\pi\in\Pi(A).$ Then $\pi(A)$ is an
irreducible operator algebra on $\mathcal{H}_{\pi}.$ By Definition \ref{D3.4},
$\pi(A)$ has a non-zero $P_{n}$-ideal $I$. As the identity representation of
$I$ on $\mathcal{H}_{\pi}$ is irreducible, $\dim\mathcal{H}_{\pi}\leq n.$ Thus
$A\in P_{n}.$ So $\mathsf{G}_{_{\Pi}}P_{n}\subseteq P_{n}$ and (\ref{3.25}) holds.

(ii) Clearly, $Comm\subseteq P_{1}\subseteq r_{_{F_{1}}}=\{A\in\mathfrak{A}$:
$\cap\{\ker\pi$: $\dim\pi=1,\pi\in\Pi(A)\}=\{0\}\}.$

Conversely, let $A\in r_{_{F_{1}}}$ and $xy\neq yx$ for some $x,y\in A$. Then
there is $\pi\in\Pi(A)$ with $\dim\pi=1$ and $0\neq\pi(xy-yx).$ On the other
hand, $\pi(x)$ and $\pi(y)$ commute, as $\dim\pi=1.$ So $\pi(xy-yx)=\pi
(x)\pi(y)-\pi(y)\pi(x)=0$. This contradiction shows that $A\in Comm$. So
$Comm=P_{1}=r_{_{F_{1}}}.$

(iii) Clearly, $P_{\infty}\subseteq r_{_{F_{\infty}}}.$ We get from
(\ref{4.13})-(\ref{4.7})\ that $r_{_{F_{\infty}}}=RFD$ is a lower\emph{
}stable and extension stable C*-property. It follows from Remark \ref{R1} that
$r_{_{F_{\infty}}}$ is not upper stable. So $P_{\infty}\neq r_{_{F_{\infty}}}$
and, by Theorem \ref{T3.10}, $\mathsf{G}_{_{\Pi}}r_{_{F_{\infty}}}\neq
r_{_{F_{\infty}}}.\bigskip$
\end{proof}

We consider now the case $1<n<\infty.$

\begin{proposition}
\label{P4.7}Let $A\in r_{_{F_{n}}},$ $1<n<\infty.$ If $\beta\in\Pi(A)$ is such
that $C(\mathcal{H}_{\beta})\subseteq\beta(A),$ then $\dim\beta\leq n,$
i.e.\emph{, }$\beta\in F_{n}(A).$
\end{proposition}

\begin{proof}
For each $\pi\in F_{n}(A),$ $\dim\pi\leq n.$ So $\pi(A)$ is isomorphic to the
algebra $M_{k}(\mathbb{C)}$ of all complex $k\times k$ matrices for some
$k\leq n.$ We consider $M_{k}(\mathbb{C})$ as a subalgebra of $M_{n}%
(\mathbb{C}).$

It follows from Corollary 6.2.1 \cite{H} and the comments after it that the
polinomial
\[
p(x_{1},...,x_{2n})=\sum_{\sigma\in\Sigma_{2n}}(-1)^{\sigma}x_{\sigma(1)}%
\cdot\cdot\cdot x_{\sigma(2n)}\text{ satisfies }p(X_{1},...,X_{2n})=0
\]
for all $X_{i}\in M_{n}(\mathbb{C}),$ where $\Sigma_{2n}$ is the symmetric
group of degree $2n.$ For each $\pi\in F_{n}(A),$ we have $\pi(a_{1}%
),...,\pi(a_{2n})\in M_{n}(\mathbb{C}),$ so that%
\[
0=p(\pi(a_{1}),...,\pi(a_{2n}))=\pi(p(a_{1},...,a_{2n}))\text{ for all }%
a_{i}\in A.
\]
Hence, by (\ref{3.19}), $p(a_{1},...,a_{2n})=0$ for all $a_{i}\in A,$ and
deg$(p)=2n.$

Let $\beta\in\Pi(A).$ Then%
\[
0=\beta(p(a_{1},...,a_{2n}))=p(\beta(a_{1}),...,\beta(a_{2n}))\text{ for all
}a_{i}\in A.
\]
Hence $p(X_{1},...,X_{2n})=0$ for all $X_{i}\in\beta(A)\subseteq
B(\mathcal{H}_{\beta}).$

Set $m=\dim\beta.$ If $n+1\leq m<\infty$ then $\beta(A)=M_{m}(\mathbb{C}),$ as
$\beta$ is irreducible. So $p(X_{1},...,X_{2n})=0$ for all $X_{i}\in
M_{m}(\mathbb{C}).$ However, it follows from Lemma 6.3.1 \cite{H} that
$M_{m}(\mathbb{C})$ does not satisfy a polynomial identity of degree less than
$2m.$ Since deg$(p)=2n<2m,$ this contradiction shows that $\dim\beta\leq n.$
So $\beta\in F_{n}(A).$

Let $m=\infty$ and $C(\mathcal{H}_{\beta})\subseteq\beta(A).$ Then there is a
subalgebra $B_{n+1}\subset C(\mathcal{H}_{\beta})\subseteq\beta(A)$ isomorphic
to $M_{n+1}(\mathbb{C}).$ Hence $p(X_{1},...,X_{2n})=0$ for all $X_{i}\in
B_{n+1}=M_{n+1}(\mathbb{C}).$ As above, we get a contradiction. Thus
$\dim\beta\leq n$ and $\beta\in F_{n}(A).$
\end{proof}

\begin{corollary}
\label{C4.0}If a $GCR$-algebra $A\in r_{_{F_{n}}},$ $n<\infty,$ then $A\in
P_{n}.$ So $P_{n}=r_{_{F_{n}}}\cap GCR.$
\end{corollary}

\begin{problem}
\label{P4.2}\emph{By Proposition \ref{P4.6},} $P_{1}=r_{_{F_{1}}}$ \emph{and}
$P_{\infty}\subsetneqq r_{_{F_{\infty}}}.$ \emph{By Corollary \ref{C4.0},}
$P_{n}=r_{_{F_{n}}}\cap GCR.$ \emph{Does }$P_{n}=r_{_{F_{n}}}$\emph{ for}
$1<n<\infty?$
\end{problem}

\section{Relations generated by special C*-properties}

In this section we illustrate the results of the previous sections by
considering various (mostly well-known) classes of C*-algebras. Some of them
are wider than but similar to the classes of $CCR$- and $GCR$-algebras and the
other are generated by the real rank zero, AF, nuclear and exact C*-algebras.
Many results of this section are well known. However, our main purpose is to
show that they arise naturally from the theory of relations in lattices
developed in the previous sections.

\subsection{C*-properties that consist of simple algebras}

In what follows by $S$ we denote any C*-property that consists of
\textbf{simple} C*-algebras$.$ As concrete examples one can have in mind the
class of all simple algebras, or the class $C$ of all algebras isomorphic to
the algebras $C(H)$ of all compact operators on separable Hilbert spaces $H$
with $\dim H\leq\infty$, or the class of Cuntz algebras $O_{n}$, or uniformly
hyperfinite algebras, etc.

Using $S$ and $\mathsf{R}(S)$ as "basic" C*-properties, we construct wider
C*-properties using Definitions \ref{D3.1}\emph{ }and \ref{D3.2}. It follows
from Definitions \ref{D3.1} and \ref{D3.2} that a C*-algebra $A$ is\medskip

\qquad$1)$\emph{\ }a \textsf{G}$(\mathsf{R}(S))$\textit{-algebra}
(\textsf{G}$S$-algebra) if each its non-zero quotient has a non-zero
$\mathsf{R}(S)$-ideal ($S$-ideal)$;\smallskip$

\qquad$2)$ an $\mathsf{NG}(\mathsf{R}(S))$\textit{-algebra} ($\mathsf{NG}%
S$\textit{-}algebra)\textit{ }if it has no $\mathsf{R}(S)$-ideals
($S$-ideals);\emph{\smallskip}

\qquad$3)$ a $\mathsf{dG}(\mathsf{R}(S)$-\textit{algebra} ($\mathsf{dG}%
S$-algebra) if each ideal of $A$ has a $\mathsf{R}(S)$-quotient ($S$%
-quotient);\emph{\smallskip}

\qquad$4)$ a \textsf{d}$\mathsf{NG}(\mathsf{R}(S)$-\textit{algebra}
(\textsf{d}$\mathsf{NG}S$-algebra) if $A$ has no $\mathsf{R}(S)$-quotients
($S$-quotients)\emph{.\medskip}

As in (\ref{4.1}), for each C*-property and $A\in\mathfrak{A}$, we consider
the corresponding relations. Since each C*-property $S$ is lower and upper
stable, Theorems \ref{T4.1} and \ref{T3.6} and Corollary \ref{C3.12} yield

\begin{corollary}
\label{C3.9}\emph{(i) }The C*-property $\mathsf{R}(S)$ is lower and upper
stable$,$ and%
\begin{equation}
S\subseteq\mathsf{R}(S)=\mathsf{R}(\mathsf{R}(S)\subseteq\mathsf{dG}%
S=\mathsf{dG}\mathsf{R}(S)\text{ and }\mathsf{G}S\cup\mathsf{R}(S)\subseteq
\mathsf{G}(\mathsf{R}(S))\cap\mathsf{R}(\mathsf{G}S). \label{4.5}%
\end{equation}

\emph{(ii) }The relations $\ll_{_{S}}$ and $\ll_{_{\mathsf{R}(S)}}$ are
$\mathbf{H}$- and dual $\mathbf{H}$-relations in each \emph{Id}$_{A}%
.\smallskip$

\emph{(iii) }$\ll_{_{S}}^{\triangleleft}=$ $\ll_{_{\mathsf{dG}S}}=$
$\ll_{_{\mathsf{dG}\mathsf{R}(S)}}=$ $\ll_{_{\mathsf{R}(S)}}^{\triangleleft}$
is a dual $\mathbf{R}$-order\emph{,} so that\emph{ }$\mathfrak{p}_{_{S}%
}^{\triangleleft}=\mathfrak{p}_{_{\mathsf{dG}S}}=\mathfrak{p}_{_{\mathsf{dG}%
\mathsf{R}(S)}}=\mathfrak{p}_{_{\mathsf{R}(S)}}^{\triangleleft}.\smallskip$

\emph{(iv) }A C$^{\ast}$-algebra $A$ is a $\mathsf{dG}S$-algebra\emph{ }if and
only if\emph{, }for each $\{0\}\neq I\in$ \emph{Id}$_{A},$ there is
$\pi^{\prime}\in\Pi(I)$ such that $\pi^{\prime}(I)\in S.$
\end{corollary}

We write \textsf{G}$\mathsf{R}(S)$ for \textsf{G}$(\mathsf{R}(S)),$
$\mathsf{NG}\mathsf{R}(S)$ for $\mathsf{NG}(\mathsf{R}(S)),$ $\mathsf{dG}%
\mathsf{R}(S)$ for $\mathsf{dG}(\mathsf{R}(S))$ and \textsf{d}$\mathsf{NG}%
\mathsf{R}(S)$ for \textsf{d}$\mathsf{NG}(\mathsf{R}(S))$. For $S=C,$ the
C*-property $\mathsf{R}(C)$ is usually denoted by $CCR,$
\begin{equation}
\mathsf{GR}(C)\text{ by }GCR,\text{ }\mathsf{NG}\mathsf{R}(C)\text{ by
}NGCR,\text{ }\mathsf{dG}\mathsf{R}(C)\text{ by }\mathsf{d}GCR,\text{
}\mathsf{dNG}\mathsf{R}(C)\text{ by }\mathsf{d}NGCR. \label{6}%
\end{equation}

\paragraph{The C*-properties of GSR\textbf{-} and NGSR-algebras.}

Combining some previous results yields

\begin{theorem}
\label{T4.2}\emph{(i) }$\mathsf{G}\mathsf{R}(S)$ is a lower and upper stable
C*-property$;$ $\mathsf{NG}\mathsf{R}(S)$ is lower stable\emph{.}$\smallskip$

\emph{(ii)} $\ll_{_{\mathsf{G}\mathsf{R}(S)}}=$ $\ll_{_{\mathsf{R}(S)}%
}^{\triangleright}$ is an $\mathbf{R}$-order\emph{,} $\ll_{_{\mathsf{NG}%
\mathsf{R}(S)}}=$ $\overleftarrow{\ll_{_{P}}}$ is a dual $\mathbf{R}$-order
and $\mathfrak{r}_{_{\mathsf{G}\mathsf{R}(S)}}=\mathfrak{r}_{_{\mathsf{R}(S)}%
}^{\triangleright}=\mathfrak{p}_{_{\mathsf{NG}\mathsf{R}(S)}}$ in
\emph{Id}$_{A}$ for each $A\in\mathfrak{A}.\medskip$

\emph{(iii)}\ If $\mathfrak{r}_{_{\mathsf{R}(S)}}^{\triangleright}\nsubseteq
I\neq A$ then there is $J\in$ \emph{Id}$_{A}$ such that $J/I$ is a
$\mathsf{R}(S)$-algebra.\medskip

\emph{(iv)\ }If $\mathfrak{r}_{_{\mathsf{R}(S)}}^{\triangleright}+I=A$ then
$A/I$ is a $\mathsf{G}\mathsf{R}(S)$-algebra.
\end{theorem}

\begin{proof}
(i) follows from Lemma \ref{L3.6} and Theorem \ref{T3.6}. Part (ii) follows
from Theorem \ref{T3.0}. Part (iii) follows from Corollary \ref{C4.2n} and
(iv) from Proposition \ref{P3.1}(ii).\bigskip
\end{proof}

The reflexive relation $\ll_{_{\mathsf{R}(S)}},$ generally, is not transitive.
For example, if $S=C$ then%
\[
\{0\}\ll_{_{CCR}}C(\mathcal{H})\ll_{_{CCR}}C(\mathcal{H})+\mathbb{C}%
\mathbf{1,}\text{ while }\{0\}\not \ll _{_{CCR}}C(\mathcal{H})+\mathbb{C}%
\mathbf{1,}%
\]
if $\dim\mathcal{H}=\infty.$ Theorem \ref{T3.1} and Proposition \ref{P3.3P}
give the following extension of well-known results for $GCR$- and
$NGCR$-algebras (\cite{D}).

\begin{theorem}
\label{C4.6}\emph{(i)} The radical $\mathfrak{r}_{_{\mathsf{R}(S)}%
}^{\triangleright}$ is the largest $\mathsf{G}\mathsf{R}(S)$-ideal of $A$ and
the smallest ideal with $N\mathsf{G}\mathsf{R}(S)$-quotient$.$ There is an
ascending transfinite $\ll_{_{\mathsf{R}(S)}}$-series of ideals of $A$ from
$\{0\}$ to $\mathfrak{r}_{_{\mathsf{R}(S)}}^{\triangleright}$.\smallskip

\emph{(ii) \ }$A$ is a $\mathsf{G}\mathsf{R}(S)$-algebra if and only if
$A=\mathfrak{r}_{_{\mathsf{R}(S)}}^{\triangleright};$ it is a $\mathsf{NG}%
\mathsf{R}(S)$-algebra if and only if $\mathfrak{r}_{_{\mathsf{R}(S)}%
}^{\triangleright}=\{0\}.\smallskip$

\emph{(iii) }If\emph{ }some $I\in$ \emph{Id}$_{A}$ and $A/I$ are
$\mathsf{G}\mathsf{R}(S)$-algebras then $A$ is a $\mathsf{G}\mathsf{R}%
(S)$-algebra.\smallskip

\emph{(iv)\ }If\emph{ }some $I\in$ \emph{Id}$_{A}$ and $A/I$ are
$\mathsf{NG}\mathsf{R}(S)$-algebras then $A$ is a $\mathsf{NG}\mathsf{R}(S)$-algebra.
\end{theorem}

\paragraph{The C*-properties of \textsf{d}GSR\textbf{-} and \textsf{d}%
NGSR-algebras.}

Clearly, $C(\mathcal{H})+\mathbb{C}\mathbf{1}$ is a $\mathsf{d}GCR$-algebra,
and $B(\mathcal{H})$ is a \textsf{d}$NGCR$-algebra, if $\dim\mathcal{H}%
=\infty.$

\begin{lemma}
\label{L4.5}\emph{(i) }$A\in$ $\mathsf{dG}\mathsf{R}(S)$ if and only if$,$ for
each $I\in$ \emph{Id}$_{A},$ there is $\pi^{\prime}\in\Pi(I)$ such that
$\pi^{\prime}(I)\in S.\smallskip$

\emph{(ii) } If there are representations $\{\pi_{\lambda}\}_{\lambda
\in\Lambda}$ in $\Pi(A)$ such that $\cap_{\lambda\in\Lambda}\ker\pi_{\lambda
}=\{0\}$ and $\pi_{\lambda}(A\mathcal{)}\in S$ for each $\lambda\in\Lambda,$
then $A$ is a $\mathsf{dG}\mathsf{R}(S)$-algebra.\smallskip

\emph{(iii) }$A$ is a $\mathsf{dNG}\mathsf{R}(S)$-algebra if and only if
$\pi(A)\notin S$ for all $\pi\in\Pi(A).$
\end{lemma}

\begin{proof}
(i) If $A\in\mathsf{dG}\mathsf{R}(S)$ and $\{0\}\neq I\in$ Id$_{A},$ then
$I/J$ is a $\mathsf{R}(S)$-algebra for some $J\subsetneqq I.$ Hence
$\pi(I/J)\in S$ for all $\pi\in\Pi(I/J).$ By (\ref{4.13}), each $\pi$ extends
to $\pi^{\prime}\in\Pi(I)$ with $\pi^{\prime}(I)=\pi(I/J)\in S.$

Conversely, if $\pi^{\prime}(I)\in S$ for each $I\in$ Id$_{A}$ and some
$\pi^{\prime}\in\Pi(I),$ then $I/\ker\pi^{\prime}\in S\subseteq\mathsf{R}(S)$.
So $A\in\mathsf{dG}\mathsf{R}(S)$.

(ii) Let $\{0\}\neq I\in$ Id$_{A}.$ As $S$ consists of simple algebras, either
$\pi_{\lambda}(I)=\pi_{\lambda}(A)\in S,$ or $I\subseteq\ker\pi_{\lambda}$ for
each $\lambda\in\Lambda.$ As $\cap_{\lambda\in\Lambda}\ker\pi_{\lambda
}=\{0\},$ $\pi_{\lambda}(I)\in S$ for some $\lambda\in\Lambda.$ So $A$ is a
$\mathsf{dG}\mathsf{R}(S)$-algebra by (i).

(iii) Let $A$ be a $\mathsf{dNG}\mathsf{R}(S)$-algebra. If $\pi(A)\in
S\subset\mathsf{R}(S)$ for some $\pi\in\Pi(A),$ then $A/\ker\pi\approx\pi(A)$
is a $\mathsf{R}(S)$-algebra, a contradiction (Definition \ref{D3.2}).

Conversely, let $\pi(A)\notin S$ for all $\pi\in\Pi(A).$ If $A$ is not
$\mathsf{dNG}\mathsf{R}(S)$-algebra then $A/I$ is a $\mathsf{R}(S)$-algebra
for some $I\in$ Id$_{A}.$ Then $\pi(A/I)\in S$ for all $\pi\in\Pi(A/I).$ By
(\ref{4.13e}), each $\pi$ extends to $\pi^{\prime}\in\Pi(A)$ such that
$\pi^{\prime}(A)=\pi(A/I)\in S,$ a contradiction.\bigskip
\end{proof}

If all $\dim\pi_{\lambda}<\infty$ in Lemma \ref{L4.5}(ii) then $A$ is called
\textit{residually finite-dimensional} (RFD) (see Example \ref{exDG}). Hence,
by Lemma \ref{L4.5}, RFD-algebras are $\mathsf{d}GCR$-algebras.

\begin{example}
\label{E5}$1)$ \emph{The group C*-algebra C}$^{\ast}$\emph{(}$\mathbb{F}_{2}%
)$\emph{ of the free group }$\mathbb{F}_{2}$\emph{ on }$2$\emph{ generators
has finite dimensional representations }$\{\pi_{k}\}_{k\geq1}$\emph{ with
}$\cap_{k\geq1}\ker\pi_{k}=\{0\}$ \emph{(Proposition VII.6.1 \cite{Da}). Thus
C}$^{\ast}$\emph{(}$\mathbb{F}_{2})$\emph{ is an RFD-algebra. So it is a
}$\mathsf{d}GCR$\emph{-algebra.\smallskip}

$2)$ \emph{The C*-algebra }$M$ \emph{of all bounded sequences }$(a_{1}%
,...,a_{n},...),$\emph{ }$a_{n}\in M_{k}(\mathbb{C}),$\emph{ is an
RFD-algebra.Thus it is a }$\mathsf{d}GCR$\emph{-algebra.\smallskip}

$3)$ \emph{The GCR-algebra $A$}$_{\infty}$ \emph{in Proposition \ref{P5.5} is
a }$\mathsf{d}NGCR$\emph{-algebra.}
\end{example}

Let $A\in\mathfrak{A.}$ As $\mathsf{R}(S)$ is a lower and upper stable,\emph{
}$\ll_{_{\mathsf{R}(S)}}$ is an $\mathbf{H}$- and a dual $\mathbf{H}$-relation
in Id$_{A}.$ Let $\mathfrak{p}_{_{\mathsf{R}(S)}}^{\triangleleft}$ be the dual
$\ll_{_{\mathsf{R}(S)}}^{\triangleleft}$-radical in Id$_{A}.$ Corollary
\ref{C4.2n}, Lemma \ref{L3.7} and Theorem \ref{T3.7} yield

\begin{theorem}
\label{T4.4}\emph{(i) }$\mathsf{dG}\mathsf{R}(S)$ is a lower and
$\mathsf{dNG}\mathsf{R}(S)$ is an upper stable \emph{C*}%
-properties.$\smallskip$

\emph{(ii)\ }$\ll_{_{\mathsf{dG}\mathsf{R}(S)}}=$ $\ll_{_{\mathsf{R}(S)}%
}^{\triangleleft}$ is a dual $\mathbf{R}$-order\emph{, }$\ll_{_{\mathsf{dNG}%
\mathsf{R}(S)}}=$ $\overleftarrow{\ll_{_{\mathsf{R}(S)}}}$ is an $\mathbf{R}%
$-order\emph{,} $\mathfrak{p}_{_{\mathsf{dN}\mathsf{R}(S)}}=\mathfrak{p}%
_{_{\mathsf{R}(S)}}^{\triangleleft}=\mathfrak{r}_{_{\mathsf{dNG}\mathsf{R}%
(S)}}^{\triangleright}$ in \emph{Id}$_{A}$ for each $A\in\mathfrak{A.}%
\smallskip$

\emph{(iii)}\ \ If $\{0\}\neq I\nsubseteq\mathfrak{p}_{_{\mathsf{R}(S)}%
}^{\triangleleft}$ then there is $J\in$ \emph{Id}$_{A}$ such that $I/J$ is a
$\mathsf{R}(S)$-algebra.\smallskip

\emph{(iv)\ }If $\mathfrak{p}_{_{\mathsf{R}(S)}}^{\triangleleft}\cap I=\{0\}$
for $I\in$ \emph{Id}$_{A},$ then $I$ is a $\mathsf{dG}\mathsf{R}(S)$-algebra.
\end{theorem}

\begin{remark}
\label{R1}$1)$ The \emph{C*}-properties $\mathsf{d}GCR$ and \emph{RFD}
\emph{(}of all \emph{RFD}-algebras\emph{)} are not upper stable.

\emph{The group C*-algebra C}$^{\ast}$\emph{(}$F_{2})$\emph{ of the free group
}$F_{2}$\emph{ is a RFD-algebra and, therefore, a }$\mathsf{d}GCR$%
-\emph{algebra (Example \ref{E5}). Let C}$_{r}^{\ast}(F_{2})=\overline
{\pi(l^{1}(F_{2}))}$\emph{ be the reduced group C*-algebra of }$F_{2}$\emph{.
It is the norm closure of the image of the representation }$\pi$\emph{ of the
algebra }$l^{1}(F_{2})$\emph{ on the Hilbert space }$l^{2}(F_{2}).$\emph{ Then
C}$_{r}^{\ast}$\emph{(}$F_{2})\approx$\emph{ C}$^{\ast}$\emph{(}$F_{2})/I$
\emph{for some ideal }$I$ \emph{of C}$^{\ast}$\emph{(}$F_{2})$\emph{. It is a
simple algebra (Corollary VII.7.5 and Theorem VII.8.6 \cite{Da}) which is not
isomorphic to }$C(\mathcal{H})$ \emph{for} $\dim\mathcal{H}\leq\infty.$
\emph{Thus C}$_{r}^{\ast}$\emph{(}$F_{2})$\emph{ is neither an RFD- nor a
}$\mathsf{d}GCR$-\emph{algebra. So the quotient C}$^{\ast}$\emph{(}$F_{2}%
)/I$\emph{ is neither an RFD- nor a }$\mathsf{d}GCR$-\emph{algebra}.\emph{
Thus the C*-properties RFD and }$\mathsf{d}GCR$ are not upper
stable.\emph{\smallskip}

$2)$\emph{ }The \emph{C*}-property $\mathsf{dNG}\mathsf{R}(S)$ is not lower
stable\emph{. For example, not each ideal of a }$\mathsf{d}NGCR$\emph{-algebra
is a }$\mathsf{d}NGCR$\emph{-algebra. Indeed,} $B(\mathcal{H})$ \emph{is a
}$\mathsf{d}NGCR$\emph{-algebra, but its ideal }$C(\mathcal{H})$\emph{ is
not}. \emph{\ }$\blacksquare$
\end{remark}

Using Theorem \ref{T3.2} and Proposition \ref{P3.5}, we have

\begin{theorem}
\label{T5.1}\emph{(i) }The dual radical $\mathfrak{p}_{_{\mathsf{R}(S)}%
}^{\triangleleft}$ is the largest $\mathsf{dNG}\mathsf{R}(S)$-ideal and the
smallest ideal with $\mathsf{dG}\mathsf{R}(S)$-quotient$.$ There is a
descending transfinite $\ll_{_{\mathsf{R}(S)}}$-series of ideals from $A$ to
$\mathfrak{p}_{_{\mathsf{R}(S)}}^{\triangleleft}$.\smallskip

\emph{(ii) }$A$ is a $\mathsf{dG}\mathsf{R}(S)$-algebra if and only if
$\mathfrak{p}_{_{\mathsf{R}(S)}}^{\triangleleft}=\{0\};$\emph{\ }it is a
$\mathsf{dNG}\mathsf{R}(S)$-algebra iff $\mathfrak{p}_{_{\mathsf{R}(S)}%
}^{\triangleleft}=A.\smallskip$

\emph{(iii)\ }If\emph{ }some $I\in$ \emph{Id}$_{A}$ and $A/I$ are
$\mathsf{dG}\mathsf{R}(S)$-algebras then $A$ is a $\mathsf{dG}\mathsf{R}%
(S)$-algebra.\smallskip

\emph{(iv)\ }If\emph{ }some $I\in$ \emph{Id}$_{A}$ and $A/I$ are
$\mathsf{dNG}\mathsf{R}(S)$-algebras then $A$ is a $\mathsf{dNG}\mathsf{R}(S)$-algebra.
\end{theorem}

Summarizing the results as in (\ref{3.5}), we get that, for each
$A\in\mathfrak{A,}$%
\begin{align}
\mathfrak{r}_{_{\mathsf{G}\mathsf{R}(S)}}  &  =\mathfrak{r}_{_{\mathsf{R}(S)}%
}^{\triangleright}=\mathfrak{p}_{_{\mathsf{NG}\mathsf{R}(S)}}\in
\mathsf{G}\mathsf{R}(S)\text{ and}\nonumber\\
\mathfrak{p}_{_{\mathsf{G}\mathsf{R}(S)}}  &  \subseteq\mathfrak{p}_{_{S}%
}^{\triangleleft}=\mathfrak{p}_{_{\mathsf{dG}S}}=\mathfrak{p}_{_{\mathsf{dG}%
\mathsf{R}(S)}}=\mathfrak{p}_{_{\mathsf{R}(S)}}^{\triangleleft}=\mathfrak{r}%
_{_{\mathsf{dNG}\mathsf{R}(S)}}\in\mathsf{dNG}\mathsf{R}(S). \label{5.1}%
\end{align}

\paragraph{The classes of $\mathsf{G}S$-$,$ $\mathsf{NG}S$-$,$ $\mathsf{dG}S$-
and \textsf{d}$\mathsf{NG}S$-algebras.}

In the previous subsections we took the class of $\mathsf{R}(S)$-algebras as
the "basic" class and "constructed" wider classes of \textsf{G}$\mathsf{R}%
(S),$ $\mathsf{NG}\mathsf{R}(S),$ $\mathsf{dG}\mathsf{R}(S)$ and
\textsf{d}$\mathsf{NG}\mathsf{R}(S)$-algebras.

Take now $S$ as the "basic" class and consider the C*-properties of
\textsf{G}$S,$ $\mathsf{NG}S,$ $\mathsf{dG}S$ and \textsf{d}$\mathsf{NG}%
S$-algebras. They define the corresponding relations in Id$_{A}$ for all
$A\in\mathfrak{A.}$ Since $S$ is a lower and upper stable C*-property, all the
results of Theorems \ref{T4.2}, \ref{C4.6}, \ref{T4.4}, \ref{T5.1} hold with
$\mathsf{R}(S)$ replaced by $S.$

Unlike the dual $\ll_{_{S}}^{\triangleleft}$- and the dual $\ll_{_{\mathsf{R}%
(S)}}^{\triangleleft}$-radicals $\mathfrak{p}_{_{S}}^{\triangleleft}$ and
$\mathfrak{p}_{_{\mathsf{R}(S)}}^{\triangleleft}$ which are always equal$,$
the $\ll_{_{S}}^{\triangleright}$- and $\ll_{_{\mathsf{R}(S)}}^{\triangleright
}$-radicals $\mathfrak{r}_{_{S}}^{\triangleright}$ and $\mathfrak{r}%
_{_{\mathsf{R}(S)}}^{\triangleright}$ may differ for some algebras $A.$

\begin{corollary}
\label{C4.1}\emph{(i)\ }$\mathfrak{r}_{_{S}}^{\triangleright}\subseteq
\mathfrak{r}_{_{\mathsf{R}(S)}}^{\triangleright}$ in each $A\in\mathfrak{A}%
.\smallskip$

\emph{(ii)} If $A$ is a $\mathsf{dG}S$-\emph{,} or an $\mathsf{R}%
(S)$-algebra\emph{, }then $\mathfrak{p}_{_{S}}^{\triangleleft}=\{0\}$ and
there is a descending transfinite $\ll_{_{S}}$-series\textit{ }of ideals from
$A$ to $\{0\}.$
\end{corollary}

\begin{proof}
(i) As $S\subseteq\mathsf{R}(S)$, we have $\mathfrak{r}_{_{S}}^{\triangleright
}\subseteq\mathfrak{r}_{_{\mathsf{R}(S)}}^{\triangleright}$ by (\ref{3.13})
for all C*-algebras $A.$

(ii) If $A\in$ $\mathsf{dG}S$ then $\mathfrak{p}_{_{\text{\textsf{d}%
}\mathsf{G}S}}=\{0\}.$ If $A\in\mathsf{R}(S)$ then, as $\mathsf{R}%
(S)\subseteq$ $\mathsf{dG}S$ by (\ref{4.5}), $\mathfrak{p}_{_{\text{\textsf{d}%
}\mathsf{G}S}}=\{0\}.$ By (\ref{5.1}), $\mathfrak{p}_{_{S}}^{\triangleleft
}=\mathfrak{p}_{_{\text{\textsf{d}}\mathsf{G}S}}.$ So $\mathfrak{p}_{_{S}%
}^{\triangleleft}=\{0\}$ and the transfinite series\textit{ }exists by Theorem
\ref{T3.2}.$\bigskip$
\end{proof}

For $S=C,$ Corollary \ref{C4.1}(i) gives a well-known result that each
$CCR$-algebra $A$ has a descending transfinite $\ll_{_{C}}$-series\textit{ }of
ideals from $A$ to $\{0\}$.

For many C*-algebras $A,$ the radicals $\mathfrak{r}_{_{S}}^{\triangleright}$
and $\mathfrak{r}_{_{\mathsf{R}(S)}}^{\triangleright}$ differ in Id$_{A}.$ For
example, let $S=C.$ Clearly $A=C(0,1)$ is a $CCR$-algebra, so that
$\mathfrak{r}_{_{CCR}}^{\triangleright}=A,$ while $\mathfrak{r}_{_{C}%
}^{\triangleright}=\{0\},$ as $A$ has no ideals isomorphic to $C(\mathcal{H}%
).$ On the other hand, for C*-algebras $A$ with separable conjugate space,
$\mathfrak{r}_{_{C}}^{\triangleright}=\mathfrak{r}_{_{CCR}}^{\triangleright}$
in Id$_{A}$ (see Corollary \ref{Tom}).

Denote by $\mathfrak{A}_{\text{sep}}$ the C*-property that consists of all
separable C*-algebras, and by $\mathfrak{A}_{\text{sep}}^{\ast}$ the
C*-property that consists of all C*-algebras with separable conjugate space.

\begin{lemma}
\label{L4.1}The \emph{C*}-properties $\mathfrak{A}_{\text{\emph{sep}}}$ and
$\mathfrak{A}_{\text{\emph{sep}}}^{\ast}$ are lower and upper stable.
\end{lemma}

\begin{proof}
Let $A\in\mathfrak{A}$, $I\in$ Id$_{A}$ and $p$: $A\rightarrow A/I$ be the
standard epimorphism. If $A\in\mathfrak{A}_{\text{sep}}$ then, clearly, $I$
and $A/I$ are seprable. So $\mathfrak{A}_{\text{sep}}$ is a lower and upper
stable C*-property.

Now let $A\in\mathfrak{A}_{\text{sep}}^{\ast}.$ Set $B=A/I$. Each functional
$g\in B^{\ast}$ extends to a functional $g^{\ast}\in A^{\ast}$ by $g^{\ast
}(x)=g(p(x))$ for $x\in A.$ The map $g\rightarrow g^{\ast}$ is an isometric
isomorphism from $B^{\ast}$ onto the closed subspace $I^{\bot}=\{f\in A^{\ast
}$: $f|_{I}=0\}$ of $A^{\ast}$ (\cite{DS} II.4.18(b)). As $A^{\ast}$ is
separable, $I^{\bot}$ is also separable. So $B^{\ast}$ is separable. Thus
$A/I\in\mathfrak{A}_{\text{sep}}^{\ast}.$

We also have that $I^{\ast}\approx A^{\ast}/I^{\bot}.$ As $A^{\ast}$ is
separable, $I^{\ast}$ is separable. Thus $I\in\mathfrak{A}_{\text{sep}}^{\ast
}.$ Hence $\mathfrak{A}_{\text{sep}}^{\ast}$ is a lower and upper stable
C*-property.\bigskip
\end{proof}

It follows from Theorem \ref{T3.4}, from (\ref{4.5}) and Example \ref{3.24}
that
\begin{equation}
\mathsf{G}C\subsetneqq GCR\subseteq\mathsf{G}_{_{\Pi}}C. \label{4.6}%
\end{equation}
However, in $\mathfrak{A}_{\text{sep}}$ the C*-properties $GCR$ and
\textsf{G}$_{_{\Pi}}C$ coincide, and in $\mathfrak{A}_{\text{sep}}^{\ast}$ the
C*-properties \textsf{G}$C,$ $GCR$ and \textsf{G}$_{_{\Pi}}C$ coincide.

\begin{proposition}
\label{P4.1}\emph{(i) }$GCR\cap\mathfrak{A}_{\text{\emph{sep}}}=\mathsf{G}%
_{_{\Pi}}C\cap\mathfrak{A}_{\text{\emph{sep}}}$.\smallskip

\emph{(ii) }$\mathfrak{A}_{\text{\emph{sep}}}^{\ast}\subseteq\mathsf{G}C,$ so
that $\mathsf{G}C\cap\mathfrak{A}_{\text{\emph{sep}}}^{\ast}=GCR\cap
\mathfrak{A}_{\text{\emph{sep}}}^{\ast}=\mathsf{G}_{_{\Pi}}C\cap
\mathfrak{A}_{\text{\emph{sep}}}^{\ast}=\mathfrak{A}_{\text{\emph{sep}}}%
^{\ast}.$
\end{proposition}

\begin{proof}
(i) Let $A\in\mathsf{G}_{_{\Pi}}C\cap\mathfrak{A}_{\text{sep}}$. By Definition
\ref{D3.4}, $\pi(A)\supseteq C(\mathcal{H}_{\pi})$ for each $\pi\in\Pi(A).$ As
$A$ is separable, it follows from Theorem 9.1 \cite{D} that $A$ is a
$GCR$-algebra. Hence, by (\ref{4.6}), $GCR\cap\mathfrak{A}_{\text{sep}%
}\subseteq\mathsf{G}_{_{\Pi}}C\cap\mathfrak{A}_{\text{sep}}\subseteq
GCR\cap\mathfrak{A}_{\text{sep}}$. So $GCR\cap\mathfrak{A}_{\text{sep}%
}=\mathsf{G}_{_{\Pi}}C\cap\mathfrak{A}_{\text{sep}}.$

(ii) Let $A\in\mathfrak{A}_{\text{sep}}^{\ast}$. By Lemma \ref{L4.1}, for each
$I\in$ Id$_{A},$ the quotient $A/I\in\mathfrak{A}_{\text{sep}}^{\ast}.$
Tomiyama \cite{To} proved that each algebra in $\mathfrak{A}_{\text{sep}%
}^{\ast}$ has an ideal isomorphic to $C(\mathcal{H})$ for a separable
$\mathcal{H}.$ Hence $A/I$ contains an ideal isomorphic to $C(\mathcal{H)}.$
Therefore $A$ is a $\mathsf{G}C$-algebra. Thus $\mathfrak{A}_{\text{sep}%
}^{\ast}\subseteq\mathsf{G}C.$ From this and from (\ref{4.6}) follows the rest
of (ii).
\end{proof}

\begin{corollary}
\label{Tom}Each $A\in\mathfrak{A}_{\text{\emph{sep}}}^{\ast}$ is a
$\mathsf{G}C$-algebra and a $GCR$-algebra\emph{,}
\[
\mathfrak{r}_{_{C}}^{\triangleright}=\mathfrak{r}_{_{CCR}}^{\triangleright
}=A\text{ and }\ll_{_{\mathsf{G}C}}\text{ }=\text{ }\ll_{_{GCR}}\text{
}=\text{ }\subseteq\text{ in \emph{Id}}_{A}.
\]
The algebra $A$ has a\emph{ }countable\emph{ }ascending transfinite $\ll
_{_{C}}$-series $\left(  I_{\lambda}\right)  _{1\leq\lambda\leq\gamma}$ of
ideals such that $I_{1}=\{0\},$ $I_{\gamma}=A$ and $I_{\lambda+1}/I_{\lambda
}\approx C(\mathcal{H}_{\lambda})$ for some separable $\mathcal{H}_{\lambda}.$
\end{corollary}

\begin{proof}
By Proposition \ref{P4.1}, if $A\in\mathfrak{A}_{\text{sep}}^{\ast},$ we have
$A$ is a \textsf{G}$C$-algebra and a $GCR$-algebra. It follows from Theorems
\ref{T3.1}(ii) and \ref{C4.6}(ii) that $\mathfrak{r}_{_{C}}^{\triangleright
}=\mathfrak{r}_{_{CCR}}^{\triangleright}=A.$

As \textsf{G}$C\subseteq GCR,$ we have $I\ll_{_{\mathsf{G}C}}J$ $\Rightarrow$
$I\ll_{_{GCR}}J$ $\Rightarrow$ $I\subseteq J$ in Id$_{A}.$ If $I\subseteq J$
then $J/I\in\mathfrak{A}_{\text{sep}}^{\ast}$ by Lemma \ref{L4.1}. Hence, by
Proposition \ref{P4.1}, $J/I\in\mathsf{G}C.$ So $I\ll_{_{\mathsf{G}C}}J.$ Thus
$\ll_{_{\mathsf{G}C}}=$ $\ll_{_{GCR}}=$ $\subseteq$ in Id$_{A}.$

As $\mathfrak{r}_{_{C}}^{\triangleright}=A,$ we have from Theorem
\ref{T3.1}(i) that there is an ascending transfinite $\ll_{_{C}}$-series of
ideals from $\{0\}$ to $A$. As $A^{\ast}$ is separable, $A$ is separable. So,
by 4.3.8 \cite{D}, the series is countable.\bigskip
\end{proof}

We will now construct a $GCR$-algebra which is a \textsf{d}$NGCR$-algebra.

Let $R=C(K\otimes H),$ where $\dim H=\dim K=\infty,$ and $A$ be a C$^{\ast}%
$-algebra in $B(K).$ Then $R\cap(A\otimes\mathbf{1}_{H})=\{0\},$
$B=R+A\otimes\mathbf{1}_{H}$ is a C*-algebra in $B(K\otimes H)$ by Corollary
1.8.4 \cite{D}, and%
\begin{equation}
\text{Id}_{B}=\{\{0\},R+I\otimes\mathbf{1}_{H}\text{: }I\in\text{ Id}_{A}\}.
\label{5.13}%
\end{equation}
Indeed, if $\{0\}\neq J\in$ Id$_{B},$ $J\cap R$ is either $\{0\}$ or $R.$ If
$J\cap R=\{0\}$ then $JR=\{0\},$ so $J=\{0\}.$ Thus $R\subseteq J$ and
$J=R+I\otimes\mathbf{1}_{H}$ where $I\in$ Id$_{A}.$ Conversely, $R+I\otimes
\mathbf{1}_{H}\in$ Id$_{B}$ for $I\in$ Id$_{A}.$

Denote by $H_{n}$ the tensor product of $n$ copies of $H.$ Let $A_{0}=\{0\},$
$A_{1}=C(H)$ and%
\begin{equation}
A_{n}=C(H_{n})+A_{n-1}\otimes\mathbf{1}_{H}=C(H_{n})+\sum_{i=1}^{n-1}%
C(H_{n-i})\otimes\mathbf{1}_{H_{i}}\mathbf{,} \label{5.12}%
\end{equation}
where $\mathbf{1}_{H_{i}}$ is the tensor product of $i$ copies of
$\mathbf{1}_{H}.$ As above, $A_{n}$ are C*-algebras and, by (\ref{5.13}),
Id$_{A_{n}}=\{I_{k}\}_{k=0}^{n}=\{\{0\}=I_{0}\subset I_{1}\subset...\subset
I_{n-1}\subset I_{n}=A_{n}\}$ where%
\begin{equation}
I_{k}=C(H_{n})+\sum_{i=1}^{k-1}C(H_{n-i})\otimes\mathbf{1}_{H_{i}}\subset
I_{k+1}=I_{k}+C(H_{n-k})\otimes\mathbf{1}_{H_{k}}\text{ for }1\leq k\leq n-1.
\label{5.b}%
\end{equation}

For an isometry $U$ from $H$ onto $H\otimes H$, set $U_{n}=U\otimes
\mathbf{1}_{H_{n-1}},$ $n\geq2.$ Then $U_{n}$ is an isometry from $H_{n}$ onto
$H_{n+1}.$ Thus $\theta_{n}$: $B\rightarrow U_{n}BU_{n}^{-1}$ is an
isomorphism of $B(H_{n})$ onto $B(H_{n+1})$ and $\theta_{n}(C(H_{n}%
))=C(H_{n+1}).$ As $U_{n}=U_{n-i}\otimes\mathbf{1}_{H_{i}}$ for all $\ i<n,$
we have $\theta_{n}=\theta_{n-i}\otimes\mathbf{1}_{H_{i}},$ so that
\begin{equation}
\theta_{n}(C(H_{n-i})\otimes\mathbf{1}_{H_{i}})=\theta_{n-i}(C(H_{n-i}%
))\otimes\mathbf{1}_{H_{i}}=C(H_{n+1-i})\otimes\mathbf{1}_{H_{i}}. \label{5.c}%
\end{equation}
Then Id$_{A_{n+1}}=\{J_{k}\}_{k=0}^{n+1},$ where $J_{0}=\{0\},$ $J_{n+1}%
=A_{n+1}$ and%
\[
J_{k}\overset{(\ref{5.b})}{=}C(H_{n+1})+\sum_{i=1}^{k-1}C(H_{n+1-i}%
)\otimes\mathbf{1}_{H_{i}}\overset{(\ref{5.c})}{=}\theta_{n}(I_{k})\text{ for
}0\leq k\leq n.
\]
Thus $\varphi_{kn}=\theta_{n}...\theta_{k}$ is an isomorphism of $A_{k}$ onto
the ideal $J_{k}.$ Identifying $A_{k}$ and $J_{k},$ we have%
\begin{equation}
\text{Id}_{A_{n+1}}=\{J_{k}\}_{k=0}^{n+1}=\{\{0\},A_{1},A_{2},...,A_{n}%
,A_{n+1}\}. \label{5.16}%
\end{equation}

The union $\cup_{n}A_{n}$ with relation: $a\in A_{k}\sim b\in A_{n+1},$ if
$\varphi_{kn}(a)=b,$ is a *-algebra with the C*-norm $\left\Vert a\right\Vert
=\left\Vert a\right\Vert _{A_{k}},$ if $a\in A_{k}.$ Its completion
$A_{\infty}\mathcal{=}$ $\overline{\cup_{n}A_{n}}$ -- the inductive limit of
$\{A_{n}\}_{n=0}^{\infty}$ -- is a C*-algebra and $A_{n}$ can be considered as
ideals of $A_{\infty}$: $\{0\}=A_{0}\subset A_{1}\subset A_{2}\subset
...\subset A_{\infty}.$

\begin{proposition}
\label{P5.5}Let $A_{\infty}=\overline{\cup_{n}A_{n}}.$ Then \emph{Id}%
$_{A_{\infty}}=\{A_{n}\}_{n=0}^{\infty}$ and $\{A_{n}\}_{n=0}^{\infty}$ is a
countable ascending $\ll_{_{C}}$-series of ideals from $\{0\}$ to $A_{\infty}%
$\emph{: } $A_{n+1}/A_{n}\approx C(H).$\smallskip

The algebra $A_{\infty}$ is a $\mathsf{G}C$-\emph{,} $GCR$- and $\mathsf{d}%
NGCR$-algebra and
\[
\mathfrak{r}_{_{C}}^{\triangleright}=\mathfrak{r}_{_{CCR}}^{\triangleright
}=\mathfrak{p}_{_{CCR}}^{\triangleleft}=\mathfrak{p}_{_{C}}^{\triangleleft
}=A_{\infty}..
\]

\end{proposition}

\begin{proof}
Let $J\in$ Id$_{A_{\infty}}.$ By Lemma III.4.1 \cite{Da}, $J=\overline
{\cup_{n}(J\cap A_{n})}.$ Each $J\cap A_{n}$ is an ideal in $A_{n}.$ By
(\ref{5.16}), $J\cap A_{n}$ coincides with some $A_{k_{n}}$ in Id$_{A_{n}},$
$k_{n}\leq n.$ Thus $J=\overline{\cup_{n}A_{k_{n}}}.$ If the sequence
$\{k_{n}\}$ is unbounded then $J=A_{\infty}.$ If it is bounded and
$k=\max\{k_{n}\}$ then $J=A_{k}.$ Thus Id$_{A_{\infty}}=\{A_{n}\}_{n=0}%
^{\infty}.$ By (\ref{5.12}), $A_{k+1}=\theta_{k}(A_{k})+C(H)\otimes
\mathbf{1}_{H_{k}},$so that $A_{k+1}/A_{k}\approx A_{k+1}/\theta_{k}%
(A_{k})\approx C(H).$ Hence $\{A_{n}\}_{n=0}^{\infty}$ is a countable
ascending $\ll_{_{C}}$-series of ideals from $\{0\}$ to $A_{\infty}.$

As the ascending $\ll_{_{C}}$-series of ideals from $\{0\}$ to $A$ exists,
$\mathfrak{r}_{_{C}}^{\triangleright}=A_{\infty}.$ As $\mathfrak{r}_{_{C}%
}^{\triangleright}\subseteq\mathfrak{r}_{_{CCR}}^{\triangleright}$ by
Corollary \ref{C4.1}, $\mathfrak{r}_{_{C}}^{\triangleright}=\mathfrak{r}%
_{_{CCR}}^{\triangleright}=A_{\infty}.$ Hence, by Theorems \ref{T3.1} and
\ref{C4.6}, $A_{\infty}$ is a \textsf{G}$C$-\emph{,} $GCR$-algebra.

By (\ref{5.1}), $\mathfrak{p}_{_{CCR}}^{\triangleleft}=\mathfrak{p}_{_{C}%
}^{\triangleleft}.$ As $A_{\infty}$ has no ideal $I$ such that $A_{\infty
}/I\approx C(H)$, we have $\mathfrak{p}_{_{C}}^{\triangleleft}=A_{\infty}$ by
Corollary \ref{C4.2n}(ii). Hence $\mathfrak{p}_{_{CCR}}^{\triangleleft
}=\mathfrak{p}_{_{C}}^{\triangleleft}=A_{\infty}$. Thus it is a \textsf{d}%
$NGCR$-algebra by Theorem \ref{T5.1}.
\end{proof}

\subsection{C*-algebras with continuous trace and dual C$^{\ast}$-algebras}

\paragraph{Continuous trace C*-property.}

Let $A$ be a $GCR$-algebra. Set
\begin{align}
B_{A}  &  =\{a\in A^{+}\text{: the function }T_{a}(\pi)=\text{ Tr }%
\pi(a)\text{ is finite and continuous on }\Pi(A)\},\nonumber\\
N_{A}  &  =\{a\in A\text{: }aa^{\ast}\in B_{A}\}\text{ and }M_{A}=N_{A}^{2}.
\label{5.0}%
\end{align}
By Lemma 4.5.1 \cite{D}, the subspaces $N_{A},$ $M_{A}$ are (non-closed)
*-ideals of $A$,%
\begin{equation}
M_{A}=\text{ span}(B_{A}),\text{ }B_{A}=M_{A}\cap A^{+}\text{ and }%
\overline{M_{A}}=\overline{N_{A}}. \label{5.a}%
\end{equation}

A C*-algebra $A$ is called a \textit{continuous trace }(c.t.)\textit{
}algebra, if $\overline{M_{A}}=A.$ Each c.t. algebra is a $CCR$-algebra
(Proposition 4.5.3 \cite{D}) and $C(H)$ is a c.t. algebra. Denote by
$P_{c.t.}$ the class of all \thinspace c.t. algebras in $\mathfrak{A}$. It is
a C*-property. Let $\ll_{_{P_{c.t.}}}$ be the corresponding relation in
Id$_{A}$ for $A\in\mathfrak{A.}$

\begin{proposition}
\label{L5.5}$P_{c.t.}$ is a lower and upper stable \emph{C*}-property\emph{;}
$\ll_{_{P_{c.t.}}}$ is an $\mathbf{H}$- and a dual $\mathbf{H}$-relation.
\end{proposition}

\begin{proof}
Let $A$ be a c. t. algebra, $I\in$ Id$_{A}$ and $p$: $A\rightarrow A/I$. Then,
by (\ref{5.a}),%
\begin{equation}
\overline{p(M_{A})}=p(\overline{M_{A}})=p(A)=A/I. \label{5.14}%
\end{equation}
By Proposition 3.2.1 \cite{D}, $\theta$: $\pi\rightarrow\pi\circ p,$ $\pi
\in\Pi(A/I),$ is a homeomorphism from $\Pi(A/I)$ onto the closed subset
$\Pi(A)_{I}=\{\sigma\in\Pi(A)$: $\sigma(I)=0\}$ of $\Pi(A).$ Let $a\in B_{A}.$
As $T_{a}(\sigma)=$ Tr $\sigma(a),$ $\sigma\in\Pi(A),$ is a finite and
continuous function on $\Pi(A),$ it is finite and continuous on $\Pi(A)_{I}$.
As $\theta$ is a homeomorphism from $\Pi(A/I)$ onto $\Pi(A)_{I},$ the function
$T_{p(a)}$ on $\Pi(A/I)$ defined by%
\[
T_{p(a)}(\pi)=\text{ Tr(}\pi(p(a))=\text{ Tr((}\pi\circ p)(a))=T_{a}(\pi\circ
p)=T_{a}(\theta(\pi))\text{ for }\pi\in\Pi(A/I),
\]
is finite and continuous on $\Pi(A/I).$ Thus $p(a)\in B_{A/I},$ so that
$p(B_{A})\subseteq B_{A/I}.$ Hence
\[
p(M_{A})\overset{(\ref{5.a})}{=}p(\text{span}(B_{A}))=\text{ span}%
(p(B_{A}))\subseteq\text{ span}(B_{A/I})=M_{A/I}.
\]
So, by (\ref{5.14}), $A/I=\overline{p(M_{A})}\subseteq\overline{M_{A/I}%
}\subseteq A/I.$ Thus $\overline{M_{A/I}}=A/I,$ so that $A/I$ is a c. t.
algebra. Hence $P_{c.t.}$ is upper stable.

As $\overline{M_{A}}=\overline{N_{A}}=A$ by (\ref{5.a}), the ideal $N_{A}$ has
a bounded a. i. $\{a_{\lambda}\}_{\lambda\in\Lambda}$ of $A.$ Thus $N_{A}\cap
I\ni a_{\lambda}x\rightarrow x$ for each $x\in I.$ Hence $N_{A}\cap I$ is a
dense *-ideal of $I.$ Set $b_{\lambda,x}=a_{\lambda}x(a_{\lambda}x)^{\ast}.$
Then $b_{\lambda,x}\in B_{A}\cap I$ by (\ref{5.0}), and $b_{\lambda
,x}\rightarrow xx^{\ast}.$ So $B_{\Lambda,I}=\{b_{\lambda,x}$: $\lambda
\in\Lambda,$ $x\in I\}\subseteq B_{A}\cap I$ is dense in $I^{+}.$

Each $\pi\in\Pi(I)$ extends to $\pi^{\prime}\in\Pi(A)$ and $\varphi$:
$\pi\rightarrow\pi^{\prime}$ is a homeomorphism from $\Pi(I)$ to the open set
$\Pi(A\mathcal{)}^{I}=\{\sigma\in\Pi(A)$: $\sigma(I)\neq0\}$ in $\Pi(A)$
(Proposition 3.2.1 \cite{D}). As $b_{\lambda,x}\in B_{A},$ the function
$T_{b_{\lambda,x}}$: $\sigma\rightarrow$ Tr $\sigma(b_{\lambda,x}),$
$\sigma\in\Pi(A),$ is finite and continuous on $\Pi(A).$ So it is continuous
on $\Pi(A\mathcal{)}^{I}$ and the function $\widehat{T}_{b_{\lambda,x}}$ on
$\Pi(I)$ defined by%
\[
\widehat{T}_{b_{\lambda,x}}(\pi)=\text{ Tr(}\pi(b_{\lambda,x}))=\text{ Tr(}%
\pi^{\prime}(b_{\lambda,x}))=T_{b_{\lambda,x}}(\varphi(\pi)),\text{ }\pi\in
\Pi(I),
\]
is finite and continuous on $\Pi(I).$ Thus $b_{\lambda,x}\in B_{I},$ so that
$B_{\Lambda,I}\subseteq B_{I}.$ As $B_{\Lambda,I}$ is dense in $I^{+},$
$B_{I}$ is dense in $I^{+}.$ So $M_{I}=$ span$(B_{I})$ is dense in $I.$ Thus
$I$ is a c. t. algebra. Hence $P_{c.t.}$ is lower stable.

It follows from Theorem \ref{T4.1} that $\ll_{_{P_{c.t.}}}$ is an $\mathbf{H}%
$- and a dual $\mathbf{H}$-relation.\bigskip
\end{proof}

Let $\mathfrak{r}_{_{P_{c.t.}}}^{\triangleright}$ be the $\ll_{_{P_{c.t.}}%
}^{\triangleright}$-radical and $\mathfrak{r}_{_{CCR}}^{\triangleright}$ be
the $\ll_{_{CCR}}^{\triangleright}$-radical in Id$_{A}.$ We generalize now
Theorems 4.5.5 and 4.7.12 a), b) \cite{D}.

\begin{proposition}
\label{T5.5}\emph{(i)} $GCR=\mathsf{G}P_{c.t.}$\emph{,} and the equality
$\mathfrak{r}_{_{CCR}}^{\triangleright}=\mathfrak{r}_{_{P_{c.t.}}%
}^{\triangleright}$ holds in \emph{Id}$_{A}$ for all $A\in\mathfrak{A.}%
\smallskip$

\emph{(ii) \ }$A$ is a $GCR$-algebra if and only if $A=\mathfrak{r}%
_{_{P_{c.t.}}}^{\triangleright},$ i.e.$,$ $A$ has an ascending transfinite
$\ll_{_{P_{c.t.}}}$-series $\left(  I_{\lambda}\right)  _{1\leq\lambda
\leq\gamma}$ of ideals\emph{,} $I_{1}=\{0\},$ $I_{\gamma}=A$ and all
$I_{\lambda+1}/I_{\lambda}$ are continuous trace algebras.\medskip

\emph{(iii) }If $\mathfrak{r}_{_{P_{c.t.}}}^{\triangleright}\nsubseteqq I\neq
A$ then $J/I$ is a continuous trace algebra for some $J\in$ \emph{Id}$_{A},$
$I\subsetneqq J$.
\end{proposition}

\begin{proof}
(i) By Proposition 4.5.3 \cite{D}, $P_{c.t.}\subseteq CCR.$ Let $B\in CCR$ and
$I\in$ Id$_{B}.$ Then $B/I\in CCR.$ By Lemma 4.4.4 \cite{D}, each
$GCR$-algebra has a c.t.-ideal. Hence $B/I$ has such an ideal$.$ Then, by
Definition \ref{D3.1}, $B\in$\textsf{G}$P_{c.t.}.$ So $P_{c.t.}\subseteq
CCR\subseteq\mathsf{G}P_{c.t.}.$ The rest follows from Proposition \ref{L3.1}.

(ii) By (i) and Theorem \ref{C4.6}, $A$ is a $GCR$-algebra if and only if
$A=\mathfrak{r}_{_{P_{c.t.}}}^{\triangleright}.$ The rest of (ii) and (iii)
follow from Corollary \ref{C4.2n}(i).\bigskip
\end{proof}

Let $\mathfrak{p}_{_{P_{c.t.}}}^{\triangleleft}$ be the dual $\ll_{_{P_{c.t.}%
}}^{\triangleleft}$-radical and $\mathfrak{p}_{_{CCR}}^{\triangleleft}$ be the
dual $\ll_{_{CCR}}^{\triangleleft}$-radical of $A.$ Recall (Definition
\ref{D3.2}) that $A$ is a $\mathsf{d}GCR$-algebra, if each $I\in$ Id$_{A}$ has
a $CCR$ quotient. For dual radicals we obtain results similar to the results
of Proposition \ref{T5.5}.

\begin{proposition}
\label{P5.6}\emph{(i)} $\mathsf{d}GCR=$ $\mathsf{dG}P_{c.t.}$\emph{,} and the
equality $\mathfrak{p}_{_{CCR}}^{\triangleleft}=\mathfrak{p}_{_{P_{c.t.}}%
}^{\triangleleft}$ holds in \emph{Id}$_{A}$ for all $A\in\mathfrak{A.}%
\smallskip$

\emph{(ii) }$A$ is a $\mathsf{d}GCR$-algebra if and only if $\mathfrak{p}%
_{_{P_{c.t.}}}^{\triangleleft}=\{0\}$\emph{, }i.e\emph{.,} $A$ has a
descending transfinite $\ll_{_{P_{c.t.}}}$-series $\left(  I_{\lambda}\right)
_{1\leq\lambda\leq\gamma}$ of ideals\emph{,} $I_{1}=A,$ $I_{\gamma}=\{0\}$ and
all $I_{\lambda}/I_{\lambda+1}$ are continuous trace algebras.\medskip

\emph{(iii) }If $\{0\}\neq I\nsubseteqq\mathfrak{p}_{_{P_{c.t.}}%
}^{\triangleleft}$ then $I/J$ is a continuous trace algebra for some $J\in$
\emph{Id}$_{A},$ $J\subsetneqq I.$
\end{proposition}

\begin{proof}
(i) By Proposition 4.5.3 \cite{D}, $P_{c.t.}\subseteq CCR.$ Let $B\in CCR$ and
$I\in$ Id$_{B}.$ Then $I\in CCR.$ For $\pi\in\Pi(I),$ $\pi(I)\in C$ by
Definition \ref{D4.1}. Hence $I/\ker\pi\approx\pi(I)\in C\subseteq P_{c.t.}.$
Thus, by Definition \ref{D3.2}, $B\in$ $\mathsf{dG}P_{c.t.}.$ So
$P_{c.t.}\subseteq CCR\subseteq$ $\mathsf{dG}P_{c.t.}.$ The rest follows from
Proposition \ref{P3.7}.

(ii) By (i) and Theorem \ref{T3.2}, $A$ is a $GCR$-algebra if and only if
$\mathfrak{p}_{_{P_{c.t.}}}^{\triangleleft}=\{0\}.$ The rest of (ii) and (iii)
follow from Corollary \ref{C4.2n}(ii).
\end{proof}

\paragraph{The C*-property of dual algebras.}

The left and right annihilators of a subset $E$ of an algebra $A$ are defined
by lan$(E)=\{a\in A$: $aE=\{0\}\}$ and ran$(E)=\{a\in A$: $Ea=\{0\}\}.$

A C*-algebra $A$ is \textit{dual }(see 4.7.20 \cite{D}) if, for each closed
left ideal $L$ and right ideal $R,$
\[
\text{lan(ran}(L))=L\text{ and ran(lan}(R))=R.
\]
If $I\in$ Id$_{A}$ then lan$(I)=$ ran$(I)$ are closed ideals (Lemma 32.4
\cite{BD}). Set an$(I):=$ lan$(I)=$ ran$(I).$ Kaplansky \cite{Ka} proved that
$A\in\mathfrak{A}$ is dual if and only if it is a C*($\infty)$-sum of the
algebras $C(\mathcal{H}_{\lambda})$ on some Hilbert spaces (restricted product
of C*-algebras \cite[1.9.14]{D}), i.e.,
\begin{equation}
A=\{a=(a_{\lambda})_{\lambda\in\Lambda}\text{: }a_{\lambda}\in C(\mathcal{H}%
_{\lambda})\text{; for any }\varepsilon>0,\text{ the set }\{\lambda\in
\Lambda\text{: }\left\Vert a_{\lambda}\right\Vert \geq\varepsilon\}\text{ is
finite}\}. \label{7.1}%
\end{equation}

We denote by $D$ the C*-property that consists of all dual C*-algebras. It can
be also characterized (\cite{AW}) as the class of all compact C*-algebras (a
Banach algebra $A$ is compact if, for each $a\in A,$ the map $x\mapsto axa$ is
a compact operator on $A$).

\begin{proposition}
\label{P7.1}The \emph{C*}-property $D$ is lower and upper stable.
\end{proposition}

\begin{proof}
If $A\in D$ and $I\in$ Id$_{A}.$ Then, for each $a\in I,$ the operator
$x\mapsto axa$ is compact on $I$. So $I$ is a compact algebra whence $I\in D$.
Thus $D$ is a lower stable C*-property.

By Lemma 32.4 \cite{BD}, $I\oplus$an$(I)=A.$ Hence $A/I\approx$ an$(I).$ As
an$(I)\in D$ by the above, $A/I\ $is dual. Thus $D$ is upper stable.\bigskip
\end{proof}

For $A\in\mathfrak{A,}$ define the relation $\ll_{_{D}}$ on Id$_{A}$ as in
(\ref{4.1}). It is not transitive, as $\{0\}\ll_{_{D}}C(\mathcal{H})\ll_{_{D}%
}C(\mathcal{H})+\mathbb{C}\mathbf{1}$, if $\dim\mathcal{H}=\infty$, while
$C(\mathcal{H})+\mathbb{C}\mathbf{1}$ is not dual. By Theorem \ref{T4.1} and
Proposition \ref{P7.1}, $\ll_{_{D}}$ is an $\mathbf{H}$- and a dual
$\mathbf{H}$-relation. Let $\mathfrak{r}_{_{D}}^{\triangleright}$ be the
$\ll_{_{D}}^{\triangleright}$-radical and $\mathfrak{p}_{_{D}}^{\triangleleft
}$ the dual $\ll_{_{D}}^{\triangleleft}$-radical in $A.$

Recall that each $A\in\mathfrak{A}$ has the $\ll_{_{C}}^{\triangleright}%
$-radical $\mathfrak{r}_{_{C}}^{\triangleright}$ and the dual $\ll_{_{C}%
}^{\triangleleft}$-radical $\mathfrak{p}_{_{C}}^{\triangleleft}$.

\begin{proposition}
\label{P7.2}\emph{(i) }Each $A\in D$ is a $\mathsf{G}C$- and a $\mathsf{dG}%
C$-algebra\emph{, }so that\emph{ }$D\subseteq\mathsf{G}C$ $\cap$
$\mathsf{dG}C.\smallskip$

\emph{(ii) }For $A\in\mathfrak{A},$ $\mathfrak{r}_{_{D}}^{\triangleright
}=\mathfrak{r}_{_{C}}^{\triangleright}$ and $\mathfrak{p}_{_{D}}%
^{\triangleleft}=\mathfrak{p}_{_{C}}^{\triangleleft}.$ If $A$ is dual then
$\mathfrak{r}_{_{C}}^{\triangleright}=A$ and $\mathfrak{p}_{_{C}%
}^{\triangleleft}=\{0\}.$
\end{proposition}

\begin{proof}
(i) Let $A\in D$ and $I\in$ Id$_{A}.$ By Proposition \ref{P7.1}, $I$ and $A/I$
are dual. By (\ref{7.1}), dual C*-algebras have ideals isomorphic to
$C(\mathcal{H)}.$ So $A/I$ has such an ideal. Thus $A$ is a \textsf{G}$C$-algebra.

As $I$ is dual, it has a quotient isomorphic to $C(\mathcal{H)}$ by
(\ref{7.1}). So $A$ is a $\mathsf{dG}C$-algebra.

(ii) If $B\in C$ then $B\approx C(\mathcal{H}).$ So $B$ is dual$.$ Thus
$C\subseteq D.$ By (i), $C\subseteq D\subseteq\mathsf{G}C$ and $C\subseteq
D\subseteq$ $\mathsf{dG}C.$ Hence, by Propositions \ref{L3.1} and \ref{P3.7},
$\mathfrak{r}_{_{D}}^{\triangleright}=\mathfrak{r}_{_{C}}^{\triangleright}$
and $\mathfrak{p}_{_{D}}^{\triangleleft}=\mathfrak{p}_{_{C}}^{\triangleleft}$
for all $A\in\mathfrak{A.}$\bigskip
\end{proof}

Combining the above results with (\ref{5.1}), we have for each $A\in
\mathfrak{A,}$%
\begin{align*}
\mathfrak{p}_{_{\mathsf{NG}C}}^{\triangleleft}  &  =\mathfrak{r}%
_{_{\mathsf{G}C}}^{\triangleright}=\mathfrak{r}_{_{C}}^{\triangleright
}=\mathfrak{r}_{_{D}}^{\triangleright}\subseteq\mathfrak{r}_{_{CCR}%
}^{\triangleright}=\mathfrak{r}_{_{P_{c.t.}}}^{\triangleright}=\mathfrak{r}%
_{_{GCR}}^{\triangleright}=\mathfrak{p}_{_{NGCR}}^{\triangleleft},\\
\mathfrak{p}_{_{GCR}}^{\triangleleft}  &  \subseteq\mathfrak{p}_{_{CCR}%
}^{\triangleleft}=\mathfrak{p}_{_{C}}^{\triangleleft}=\mathfrak{p}%
_{_{P_{c.t.}}}^{\triangleleft}=\mathfrak{p}_{_{D}}^{\triangleleft
}=\mathfrak{p}_{_{\mathsf{d}GCR}}^{\triangleleft}=\mathfrak{r}_{_{\mathsf{d}%
\ }}^{\triangleright}.
\end{align*}

\begin{corollary}
\label{vstav}\emph{(i)} $\mathsf{G}D=Sc$ - the \emph{C*}-property of scattered
algebras \emph{(}see Example \emph{\ref{exG})}.\smallskip

\emph{(ii)} $\ \mathsf{dG}D=\mathsf{dG}C\supset RFD\cup CCR.\smallskip$

\emph{(iii)} $\mathsf{R}(D)=CCR$.
\end{corollary}

\begin{proof}
(i) Since $C\subset D\subset GC$ by Proposition \ref{P7.2}, $\mathsf{G}%
D=\mathsf{G}C$ by Proposition \ref{L3.1}. So $\mathsf{G}D=Sc,$ since
$\mathsf{G}C=Sc$ by Example \ref{exG}.

(ii) Similarly, $\mathsf{dG}D=\mathsf{dG}C$. The inclusion $\mathsf{dG}%
C\supset CCR$ follows from Theorem \ref{T3.6} for $P=C.$ If $A\in RFD$ then
$\cap\{\ker\pi\emph{:\ }\pi\in\Pi(A),\dim\pi<\infty\}=\{0\}$ by (\ref{3.30}).
Hence, for each $J\in$ Id$_{A}$, there is $\pi\in\Pi(A),$ $\dim\pi<\infty,$
with $\pi(J)\neq0$. So the algebra $J/(J\cap\ker(\pi))\in C$. Thus
$A\in\mathsf{dG}C.$

(iii) Let $A\in\mathsf{R}(D).$ Then, for each $\pi\in\Pi(A),$ the algebra
$\pi(A)\in D$. As $\pi(A)$ is irreducible, its center is trivial. So it
follows from (\ref{7.1}) that $\pi(A)\approx C(\mathcal{H}).$ Thus $\pi(A)\in
C$. Hence $\mathsf{R}(D)\subset\mathsf{R}(C)=CCR$. The converse inclusion is evident.
\end{proof}

\subsection{Real rank zero, AF, nuclear and exact C*-algebras}

A unital C*-algebra $A$ has \textbf{real rank zero }(\cite{BP}, \cite{Da}%
)\textbf{ }if its invertible selfadjoint elements are dense in the set of all
selfadjoint elements of $A.$ A non-unital algebra is real rank zero if its
unitization is real rank zero. Denote by $R_{0}$ the class of all real rank
zero algebras. Clearly, $R_{0}$ is a C*-property.

Let $I\in$ Id$_{A}.$ By Theorem 3.14 \cite{BP}, $A$ is a $R_{0}$-algebra if
and only if $I$ and $A/I$ are $R_{0}$-algebras and all projections in $A/I$
lift to projections in $A$. In Remark 3.17 \cite{BP} Brown and Pedersen write
that "split extensions of real rank zero algebras by other real rank zero
algebras produce algebras of real rank zero. For general extensions this is no
longer true: Every Bunce-Deddens algebra $A_{_{BD}}$ has a one-dimensional
extension $\widehat{A}_{_{BD}}$ (determined by a nonliftable projection in the
corona), with real rank one." Combining this with Theorem 3.14 \cite{BP} yields

\begin{corollary}
\label{C8.1}The \emph{C*}-property $R_{0}$ is lower and upper stable\emph{,
}but not extension stable.
\end{corollary}

By Theorem \ref{T4.1}, Proposition \ref{P3.3} and Corollary \ref{C8.1}, the
relation $\ll_{_{R_{0}}}$ in Id$_{A},$ $A\in\mathfrak{A},$ (see (\ref{4.1})$)$
is an $\mathbf{H}$- and a dual $\mathbf{H}$-relation, but not necessarily
transitive. Let $\mathfrak{r}_{_{R_{0}}}^{\triangleright}$ be the
$\ll_{_{R_{0}}}^{\triangleright}$-radical and $\mathfrak{p}_{_{R_{0}}%
}^{\triangleleft}$ be the dual $\ll_{_{R_{0}}}^{\triangleleft}$-radical in
Id$_{A}.$ The ideal $\mathfrak{r}_{_{R_{0}}}^{\triangleright}$ is not
necessarily an $R_{0}$-algebra. Indeed, we have $0\ll_{_{R_{0}}}A_{_{BD}}%
\ll_{_{R_{0}}}\widehat{A}_{_{BD}},$ so that $\mathfrak{r}_{_{R_{0}}%
}^{\triangleright}=\widehat{A}_{_{BD}}$ which is not an $R_{0}$-algebra.

\begin{theorem}
\label{T8.1}\emph{(i) }For each $A\in\mathfrak{A},$ $\mathfrak{r}_{_{C}%
}^{\triangleright}\subseteq\mathfrak{r}_{_{R_{0}}}^{\triangleright}$ and
$\mathfrak{p}_{_{R_{0}}}^{\triangleleft}\subseteq\mathfrak{p}_{_{C}%
}^{\triangleleft}$ in \emph{Id}$_{A}.\smallskip$

\emph{(ii) } $\mathfrak{r}_{_{R_{0}}}^{\triangleright}$ is the largest
$\mathsf{G}R_{0}$-ideal of $A$\emph{ }and contains all $R_{0}$-ideals$.$ It is
also the smallest ideal with $\mathsf{NG}R_{0}$-quotient$.\smallskip$

\emph{(iii) }$A$ is a $\mathsf{G}R_{0}$-algebra if and only if $A=\mathfrak{r}%
_{_{R_{0}}}^{\triangleright};$ it is a $\mathsf{NG}R_{0}$-algebra if and only
if $\mathfrak{r}_{_{R_{0}}}^{\triangleright}=\{0\}.\smallskip$

\emph{(iv) }$\mathfrak{p}_{_{R_{0}}}^{\triangleleft}$ is the largest
$\mathsf{dNG}R_{0}$-ideal of $A.$ It is the smallest ideal with $\mathsf{dG}%
R_{0}$-quotient.\smallskip

\emph{(v) }$A$ is a $\mathsf{dG}R_{0}$-algebra if and only if $\mathfrak{p}%
_{_{R_{0}}}^{\triangleleft}=\{0\}\emph{;}$ it is a $\mathsf{dNG}R_{0}$-algebra
if and only if $\mathfrak{p}_{_{R_{0}}}^{\triangleleft}=A.$
\end{theorem}

\begin{proof}
(i) If $A\in C,$ $A\approx C(\mathcal{H})$ for some separable $\mathcal{H}.$
So $A\in R_{0}.$ Thus $C\subset R_{0}$. So (i) follows from Propositions
\ref{L3.1} and \ref{P3.7}. Parts (ii) -- (v) follow from Theorems \ref{T3.1}
and \ref{T3.2}.\bigskip
\end{proof}

Corollary \ref{C4.2n} further yields

\begin{theorem}
\label{T8.2}\emph{(i)} If $\mathfrak{r}_{_{R_{0}}}^{\triangleright}\nsubseteq
I\neq A$ then there is an ideal $I\subsetneqq J$ such that $J/I$ is a $R_{0}%
$-algebra$.$\smallskip

\emph{(ii)\ \ }If $I\subseteq\mathfrak{r}_{_{R_{0}}}^{\triangleright}$ then
there is an ascending transfinite $\ll_{_{R_{0}}}$-series of ideals from $I$
to $\mathfrak{r}_{_{R_{0}}}^{\triangleright}.\smallskip$

\emph{(iii)} If $\{0\}\neq I\nsubseteq\mathfrak{p}_{_{R_{0}}}^{\triangleleft}$
then there is an ideal $J\subsetneqq I$ such that $I/J$ is $R_{0}%
$-algebra$.\smallskip$

\emph{(iv)\ }If $\mathfrak{p}_{_{R_{0}}}^{\triangleleft}\subseteq I$ then
there is a descending transfinite $\ll_{_{R_{0}}}$-series of ideals from $I$
to $\mathfrak{p}_{_{R_{0}}}^{\triangleleft}.$
\end{theorem}

A C*-algebra $A$ is \textit{approximately finite-dimensional }($AF$-algebra)
if it is the closure of an increasing union of finite-dimensional
*-subalgebras. It is \textit{nuclear }if, for each C*-algebra $B$, the norms
$\left\Vert \cdot\right\Vert _{\max}$ and $\left\Vert \cdot\right\Vert _{\min
}$ on the algebraic tensor product $A\odot B$ coincide, so that there is only
one C*-norm on $A\odot B.$ Furthermore, $A$ is \textit{exact }if, for each
C*-algebra $B$ and every $I\in$ Id$_{B},$%
\[
0\rightarrow A\otimes_{\min}I\rightarrow A\otimes_{\min}B\rightarrow
A\otimes_{\min}B/I\text{ is an exact sequence.}%
\]

Denote by $AF,$ $NU$ and $EX$ the corresponding C*-properties of C*-algebras.
Finite-dimensional and commutative C*-algebras, $C(\mathcal{H}),$ all $AF$ and
all C*-algebras of type I are nuclear. Nuclear algebras are exact. So
\begin{equation}
AF\subseteq NU\subseteq EX. \label{6.2}%
\end{equation}
The next theorem follows from the results in \cite{Da}, \cite{K}, Corollaries
2.5, 9.3 \cite{W}, Corollary XV,3.4 \cite{T}, \cite{Br}.

\begin{theorem}
\label{NE}\emph{(i)} The \emph{C*}-properties $AF$\emph{,} $NU$ and $EX$ are
lower stable\emph{,} upper stable and closed under inductive
limits\emph{;\smallskip}

\emph{(ii)} The \emph{C*}-properties $AF$ and $NU$ are extension stable.
\end{theorem}

By Theorems \ref{T4.1} and \ref{NE}, the C*-properties $AF$, $NU$, $EX$ define
$\mathbf{H}$- and dual $\mathbf{H}$-relations $\ll_{_{AF}},$ $\ll_{_{NU}},$
$\ll_{_{EX}}$ on Id$_{A}$ for each $A\in\mathfrak{A.}$

\begin{corollary}
\label{C8.3}Theorems \emph{\ref{T8.1}} and \emph{\ref{T8.2}} hold with
\emph{"}$R_{0}$\emph{"} being replaced by $AF,$ $NU$ and $EX,$ respectively.
\end{corollary}

For nuclear and AF algebras, Theorem \ref{NE} yields (cf. Theorem 2.22,
Corollary 2.23 \cite{ST}):

\begin{corollary}
\label{C8.4}\emph{(i) }For each $A\in\mathfrak{A,}$ the relations $\ll_{_{AF}%
},$ $\ll_{_{NU}}$ are $\mathbf{R}$-orders in \emph{Id}$_{A},$ so that%
\begin{equation}
\ll_{_{AF}}\text{ }=\text{ }\ll_{_{AF}}^{\triangleright}\text{ and }\ll
_{_{NU}}\text{ }=\text{ }\ll_{_{NU}}^{\triangleright}. \label{6.1}%
\end{equation}
\emph{(ii) }For each $A\in\mathfrak{A,}$ the $\ll_{_{AF}}$-radical
$\mathfrak{r}_{_{AF}}$ is the largest AF-ideal in $A$ and the $\ll_{_{NU}}%
$-radical $\mathfrak{r}_{_{NU}}$ is the largest nuclear ideal in
$A.\smallskip$

\emph{(ii) }$\mathsf{G}(AF)=AF$ and $\mathsf{G}(NU)=NU.$
\end{corollary}

\begin{proof}
(i) By Theorem \ref{NE}, $AF$ and $NU$ are upper stable and extension stable
(Definition \ref{D3.3}) C*-properties; the inductive limit of $AF$-algebras is
an $AF$-algebra and of $NU$-algebras is a $NU$-algebra. Hence, by Corollary
\ref{C3.2}, $\ll_{_{AF}},$ $\ll_{_{NU}}$ are $\mathbf{R}$-order in Id$_{A}$
and (\ref{6.1}) holds.

(ii) From (i) and from Theorem \ref{T3.1} it follows that $\mathfrak{r}%
_{_{AF}}=\mathfrak{r}_{_{AF}}^{\triangleright}$ (the $\ll_{_{AF}%
}^{\triangleright}$-radical) is the largest $AF$ ideal in $A.$ Similarly,
$\mathfrak{r}_{_{NU}}=\mathfrak{r}_{_{NU}}^{\triangleright}$ (the $\ll_{_{NU}%
}^{\triangleright}$-radical) is the largest nuclear ideal in $A.$

(iii) follows from (\ref{6.1}) and Corollary \ref{C3.4}.$\bigskip$
\end{proof}

From (\ref{3.13}), (\ref{6.2}) and Corollary \ref{C8.4} it follows that
$\mathfrak{r}_{_{AF}}\subseteq\mathfrak{r}_{_{NU}}\subseteq\mathfrak{r}%
_{_{EX}}^{\triangleright},$ so that $A/\mathfrak{r}_{_{AF}}$ has no $AF$
ideals, $A/\mathfrak{r}_{_{NU}}$ has no nuclear ideals and $A/\mathfrak{r}%
_{_{EX}}^{\triangleright}$ has no exact ideals by Corollary \ref{C4.2n}.

\section{Topological radicals and relation-valued functions}

\subsection{Topological radicals on C*-algebras}

In this section we study link between topological radicals\ on $\mathfrak{A}$
and C*-properties in $\mathfrak{A}$. In particular, we obtain the following
result that is the consequence of Theorems \ref{Y} and \ref{T3,1}.

\begin{theorem}
\label{T4.3}A map $R$: $A\mapsto R\left(  A\right)  \in$ \emph{Id}$_{A}$ is a
topological radical on a lower and upper stable C*-property $\mathfrak{M}$ if
and only if there is an upper stable C*-property $P_{1}$ and a lower stable
C*-property $P_{2}$ such that $R(A)=\mathfrak{r}_{_{P_{1}}}^{\triangleright
}(A)=\mathfrak{p}_{_{P_{2}}}^{\triangleleft}(A)$ for all $A\in\mathfrak{M}.$
\end{theorem}

Let $\mathfrak{M}$ be a lower and upper stable C*-property. A map $R$:
$A\mapsto R\left(  A\right)  \in\mathrm{Id}_{A}$ for all $A\in\mathfrak{M,}$
is called a \textit{topological} \textit{radical }on $\mathfrak{M}$ if the the
following axioms hold for all $A,B\in\mathfrak{M}$:\medskip

$(1^{\ast}1)$ $\phi\left(  R\left(  A\right)  \right)  =R\left(  B\right)  ,$
if $A,B\in\mathfrak{M}$ and $\phi$: $A\rightarrow B$ is any
isomorphism;\smallskip

$(1^{\ast}2)$ $p\left(  R\left(  A\right)  \right)  \subseteq R\left(
p\left(  A\right)  \right)  $ for $A\in\mathfrak{M,}$ $I\in$ Id$_{A}$ and the
standard epimorphism $p$: $A\rightarrow A/I$;\smallskip

$(2^{\ast})$ $R\left(  A/R\left(  A\right)  \right)  =\{0\}$; $\ \ (3^{\ast})$
$R\left(  R\left(  A\right)  \right)  =R\left(  A\right)  $; $\ \ \ (4^{\ast
})$ $R\left(  I\right)  \subseteq R\left(  A\right)  $ for every $I\in$
Id$_{A}$.\medskip

By isomorphism in $(1^{\ast}1)$ we mean a *-isomorphism, but it was proved in
\cite[Section 2]{ST} that any topological radical is stable under all
topological isomorphisms. In other words, topological radicals on lower and
upper stable C*-properties of C*-algebras satisfy axioms of topological
radicals on classes of normed algebras.

If a map $R$ satisfies axioms $\left(  1^{\ast}1\right)  $ and $(1^{\ast}2),$
it is called a\textit{ topological preradical}. These axioms are equivalent to
the condition: $\phi\left(  R\left(  A\right)  \right)  \subseteq R\left(
B\right)  $ for any epimorphism $\phi$: $A\rightarrow B$.

Subclasses of $\mathfrak{M}$ satisfying also condition (\ref{1.7}) are
C*-properties. If, for example, a map $R$ satisfies ($1^{\ast}1)$ then the
following sets are C*-properties in $\mathfrak{M}$:%
\begin{equation}
\mathbf{Rad}\left(  R\right)  =\left\{  A\in\mathfrak{M}\text{: }R\left(
A\right)  =A\right\}  \text{ and }\mathbf{Sem}\left(  R\right)  =\left\{
A\in\mathfrak{M}\text{: }R\left(  A\right)  =0\right\}  . \label{3,4}%
\end{equation}
The following result is a simple corollary of the definitions above (in the
general case of Banach algebras the proof is given in \cite{KST1})\textbf{.}

\begin{proposition}
\label{hr}Let $\mathfrak{M}$ be a lower and upper stable C*-property and
$R$\emph{:} $A\mapsto R\left(  A\right)  \in$ \emph{Id}$_{A}$ for
$A\in\mathfrak{M},$ be a map.\smallskip

\emph{(i) }If $R$ satisfies $(1^{\ast}1)$ and $(1^{\ast}2)$ then
$P:=\mathbf{Rad}\left(  R\right)  $ is an upper stable C*-property.\smallskip

\emph{(ii) }If $R$ satisfies $(1^{\ast}1)$ and $(4^{\ast})$ then
$\mathbf{Sem}\left(  R\right)  $ is a lower stable C*-property.
\end{proposition}

Let $A\in\mathfrak{A,}$ $I,K\in$ Id$_{A}$ and $p$: $A\rightarrow A/I.$ Then
$p(K)\approx(K+I)/I\overset{(\ref{4.01})}{\approx}K/(I\cap K)$. Hence, for
each ideal $\widehat{J}$ of $p(K),$ there is $J\in$ Id$_{A}$ such that
\begin{equation}
I\cap K\subseteq J\subseteq K\text{ and }\widehat{J}\approx J/(I\cap K)\text{,
so that }p(K)/\widehat{J}\overset{(\ref{4.01})}{\approx}K/J. \label{3.17}%
\end{equation}

For an upper stable C*-property $P$ in $\mathfrak{A}$, $\ll_{_{P}}$ is an
$\mathbf{H}$-relation and $\ll_{_{P}}^{\triangleright}$ is a $\mathbf{R}%
$-order in Id$_{A}$ for each $A\in\mathfrak{M}$ by Theorems \ref{T4.1} and
\ref{inf}. Let $\mathfrak{r}_{_{P}}^{\triangleright}(A)$ be the $\ll_{_{P}%
}^{\triangleright}$-radical in Id$_{A}$.

\begin{theorem}
\label{Y}Let $\mathfrak{M}$ be a lower and upper stable C*-property in
$\mathfrak{A}.$\smallskip

\emph{(i) }For an upper stable C*-property $P,$\emph{ }the map $\mathfrak{r}%
_{_{P}}^{\triangleright}$\emph{: }$A\in\mathfrak{M}\mapsto\mathfrak{r}_{_{P}%
}^{\triangleright}\left(  A\right)  $ is a topological radical on
$\mathfrak{M},$%
\[
\mathbf{Rad}\left(  \mathfrak{r}_{_{P}}^{\triangleright}\right)
=\mathsf{G}P\cap\mathfrak{M}\text{ and }\mathbf{Sem}\left(  \mathfrak{r}%
_{_{P}}^{\triangleright}\right)  =\mathsf{NG}P\cap\mathfrak{M.}%
\]

\emph{(ii) }For a topological radical $R$ on $\mathfrak{M},$ the class
$P=\mathbf{Rad}\left(  R\right)  $ is an upper stable C*-property in
$\mathfrak{A}$\emph{, }%
\[
\mathbf{Rad}\left(  R\right)  =\mathsf{G}P\cap\mathfrak{M}\text{ and
}R(A)=\mathfrak{r}_{_{P}}^{\triangleright}(A)\text{ for all }A\in\mathfrak{M.}%
\]

\end{theorem}

\begin{proof}
(i) To prove that $\mathfrak{r}_{_{P}}^{\triangleright}$ is a topological
radical we have to verify all conditions $(1^{\ast}1)-(4^{\ast}).\smallskip$

$(1^{\ast}1).$ Let $\phi$: $A\rightarrow B$ be an isomorphism. Set
$I=\phi(\mathfrak{r}_{_{P}}^{\triangleright}(A)).$ If $\mathfrak{r}_{_{P}%
}^{\triangleright}(B)\nsubseteq I$ then, by Corollary \ref{C4.2n}(i), there is
$J\in$ Id$_{B}$ such that $I\neq J$ and $I\ll_{_{P}}J.$ Hence $J/I\in P.$ As
$\phi^{-1}$ is an isomorphism, $\phi^{-1}(J)/\phi^{-1}(I)\approx J/I.$ So
$\phi^{-1}(J)/\phi^{-1}(I)\in P.$ Thus $\mathfrak{r}_{_{P}}^{\triangleright
}(A)=\phi^{-1}(I)\ll_{_{P}}\phi^{-1}(J)$ and $\phi^{-1}(I)\neq\phi^{-1}(J)$
which contradicts Corollary \ref{C4.2n}(i). Thus $\mathfrak{r}_{_{P}%
}^{\triangleright}(B)\subseteq\phi(\mathfrak{r}_{_{P}}^{\triangleright}(A)).$

Similarly, $\mathfrak{r}_{_{P}}^{\triangleright}(A)\subseteq\phi
^{-1}(\mathfrak{r}_{_{P}}^{\triangleright}(B)),$ so that $\phi(\mathfrak{r}%
_{_{P}}^{\triangleright}(A))=\mathfrak{r}_{_{P}}^{\triangleright}(B).$ Thus
$(1^{\ast}1)$ holds.\smallskip

$(1^{\ast}2).$ Let $A\in\mathfrak{M,}$ $I\in$ Id$_{A}$ and $p$: $A\rightarrow
A/I.$ Set $K=\mathfrak{r}_{_{P}}^{\triangleright}(A)$ in (\ref{3.17}). Then,
for each ideal $\widehat{J}$ of $p(\mathfrak{r}_{_{P}}^{\triangleright}(A)),$
there is $J\in$ Id$_{A}$ such that $J\subseteq\mathfrak{r}_{_{P}%
}^{\triangleright}(A)$ and $p(\mathfrak{r}_{_{P}}^{\triangleright
}(A))/\widehat{J}\approx\mathfrak{r}_{_{P}}^{\triangleright}(A)/J$.

Let $\mathfrak{r}_{_{P}}^{\triangleright}(A)/J$ $\neq\{0\}.$ As $\mathfrak{r}%
_{_{P}}^{\triangleright}(A)\in$\textsf{G}$P$ by Theorem \ref{T3.1},
$\mathfrak{r}_{_{P}}^{\triangleright}(A)/J$ has a non-zero $P$-ideal. So
$p(\mathfrak{r}_{_{P}}^{\triangleright}(A))/\widehat{J}$ has a non-zero
$P$-ideal. Thus, by Definition \ref{D3.1}, $p(\mathfrak{r}_{_{P}%
}^{\triangleright}(A))$ is a \textsf{G}$P$-ideal of $A/I.$ As $\mathfrak{r}%
_{_{P}}^{\triangleright}(A/I)$ is the largest \textsf{G}$P$-ideal of $A/I$ by
Theorem \ref{T3.1}, $p(\mathfrak{r}_{_{P}}^{\triangleright}(A))\subseteq
\mathfrak{r}_{_{P}}^{\triangleright}(A/I)=\mathfrak{r}_{_{P}}^{\triangleright
}(p(A))$.$\smallskip$

$(2^{\ast}).$ Set $B=A/\mathfrak{r}_{_{P}}^{\triangleright}(A).$ If
$\mathfrak{r}_{_{P}}^{\triangleright}(B)\neq\{0\},$ we have from Corollary
\ref{C4.2n}(i) that there is $\widehat{I}\in$ Id$_{B}$ such that
$\{0\}\ll_{_{P}}\widehat{I},$ i.e., $\widehat{I}\in P.$ Let $I\in$ Id$_{A}$ be
such that $\widehat{I}\approx I/\mathfrak{r}_{_{P}}^{\triangleright}(A).$ Then
$I/\mathfrak{r}_{_{P}}^{\triangleright}(A)\in P.$ So $\mathfrak{r}_{_{P}%
}^{\triangleright}(A)\ll_{_{P}}I$ which contradicts Corollary \ref{C4.2n}(i).
Thus $\mathfrak{r}_{_{P}}^{\triangleright}\left(  A/\mathfrak{r}_{_{P}%
}^{\triangleright}(A)\right)  =\mathfrak{r}_{_{P}}^{\triangleright}(B)=\{0\}.$
So $(2^{\ast})$ holds.\smallskip

$(3^{\ast}).$ By Theorem \ref{T3.1}(i), $\mathfrak{r}_{_{P}}^{\triangleright
}(A)\in\mathsf{G}P.$ Hence, by Theorem \ref{T3.1}(ii), $\mathfrak{r}_{_{P}%
}^{\triangleright}(\mathfrak{r}_{_{P}}^{\triangleright}(A))=\mathfrak{r}%
_{_{P}}^{\triangleright}(A).\smallskip$

$(4^{\ast}).$ Let $I\in$ Id$_{A}.$ By Corollary \ref{C4.2n}(i), $\mathfrak{r}%
_{_{P}}^{\triangleright}(I)$ is a $\ll_{_{P}}^{\triangleright}$-successor of
$\{0\}$ and $\mathfrak{r}_{_{P}}^{\triangleright}(A)$ is the largest
$\ll_{_{P}}^{\triangleright}$-successor of $\{0\}.$ So $\mathfrak{r}_{_{P}%
}^{\triangleright}(I)\subseteq\mathfrak{r}_{_{P}}^{\triangleright
}(A).\smallskip$

Thus the map $\mathfrak{r}_{_{P}}^{\triangleright}$\emph{: }$A\mapsto
\mathfrak{r}_{_{P}}^{\triangleright}\left(  A\right)  $ is a topological
radical on $\mathfrak{M}$. By Theorem \ref{T3.1} and $(\ref{3,4})$,%
\begin{align*}
\mathbf{Rad}\left(  \mathfrak{r}_{_{P}}^{\triangleright}\right)   &  =\left\{
A\in\mathfrak{M}\text{: }\mathfrak{r}_{_{P}}^{\triangleright}(A)=A\right\}
=\mathsf{G}P\cap\mathfrak{M,}\text{ }\\
\mathbf{Sem}\left(  \mathfrak{r}_{_{P}}^{\triangleright}\right)   &  =\left\{
A\in\mathfrak{M}\text{: }\mathfrak{r}_{_{P}}^{\triangleright}%
(A)=\{0\}\right\}  =\mathsf{NG}P\cap\mathfrak{M}.
\end{align*}

(ii) By Proposition \ref{hr}, $P=\mathbf{Rad}\left(  R\right)  $ is an upper
stable C*-property and $P\subseteq\mathfrak{M}$. By Proposition \ref{L3.1},
$P\subseteq$ \textsf{G}$P$. Thus $P\subseteq\mathsf{G}P\cap\mathfrak{M}.$ If
$A\in\mathfrak{M}$ then
\begin{equation}
R(K)=\{0\}\text{ for each ideal }K\text{ of }A/R(A),\text{ as }%
R(K)\overset{(4^{\ast})}{\subseteq}R(A/R(A))\overset{(2^{\ast})}{=}\{0\}.
\label{3,6}%
\end{equation}

Let $A\in\mathsf{G}P\cap\mathfrak{M}.$ If $R(A)\neq A$ then, by Definition
\ref{D3.1}, $A/R(A)$ has a $P$-ideal $K\neq\{0\},$ i.e., $K\in P=\mathbf{Rad}%
\left(  R\right)  .$ Then $R(K)=K$ by (\ref{3,4}) which contradicts
(\ref{3,6}). Thus $R(A)=A.$ So $A\in P.$ Hence \textsf{G}$P\cap\mathfrak{M}%
\subseteq P.$ So \textsf{G}$P\cap\mathfrak{M}=P.$

Let $A\in\mathfrak{M}.$ If $\mathfrak{r}_{_{P}}^{\triangleright}(A)\nsubseteq
R(A)$ then, by Corollary \ref{C4.2n}, $R(A)\ll_{_{P}}J$ for some ideal $J\neq
R(A).$ Thus $K:=J/R(A)\neq\{0\}$ is an ideal in $A/R(A)$ and $K\in
P=\mathbf{Rad}\left(  R\right)  .$ Hence $\{0\}\neq K=R(K)$ by (\ref{3,4}),
which contradicts (\ref{3,6}). Thus $\mathfrak{r}_{_{P}}^{\triangleright
}(A)\subseteq R(A).$

Let $\mathfrak{r}_{_{P}}^{\triangleright}(A)\neq R(A).$ As $R\left(  R\left(
A\right)  \right)  =R\left(  A\right)  $ by ($3^{\ast}),$ we have $R(A)\in
P=\mathbf{Rad}\left(  R\right)  .$ As $P$ is upper stable, $R(A)/\mathfrak{r}%
_{_{P}}^{\triangleright}(A)\in P.$ So $\mathfrak{r}_{_{P}}^{\triangleright
}(A)\ll_{_{P}}R(A)$ -- contradicts Corollary \ref{C4.2n}. Thus $\mathfrak{r}%
_{_{P}}^{\triangleright}(A)=R(A).\bigskip$
\end{proof}

For a lower stable C*-property $P$ in $\mathfrak{A}$, $\ll_{_{P}}$ is a dual
$\mathbf{H}$-relation and $\ll_{_{P}}^{\triangleleft}$ is a dual $\mathbf{R}%
$-order in Id$_{A}$ for each $A\in\mathfrak{M}$ by Theorems \ref{T4.1} and
\ref{inf}. Let $\mathfrak{p}_{_{P}}^{\triangleleft}(A)$ be the dual $\ll
_{_{P}}^{\triangleleft}$-radical in Id$_{A}.$

\begin{theorem}
\label{T3,1}Let $\mathfrak{M}$ be a lower and upper stable C*-property in
$\mathfrak{A.\smallskip}$

\emph{(i) }For a lower stable C*-property $P,$ the map $\mathfrak{p}_{_{P}%
}^{\triangleleft}$\emph{: }$A\in\mathfrak{M}\mapsto\mathfrak{p}_{_{P}%
}^{\triangleleft}\left(  A\right)  $ is a topological radical on
$\mathfrak{M},$%
\[
\mathbf{Rad}\left(  \mathfrak{p}_{_{P}}^{\triangleleft}\right)  =\mathsf{dNG}%
P\cap\mathfrak{M}\text{ and }\mathbf{Sem}\left(  \mathfrak{p}_{_{P}%
}^{\triangleleft}\right)  =\mathsf{dG}P\cap\mathfrak{M}%
\]

\emph{(ii) }For a topological radical $R$ on $\mathfrak{M,}$ $P=\mathbf{Sem}%
\left(  R\right)  \subseteq\mathfrak{M}$ is a lower stable C*-property in
$\mathfrak{A}$\emph{,}%
\[
\mathbf{Sem}\left(  R\right)  =\mathsf{dG}P\cap\mathfrak{M}\text{ and
}R(A)=\mathfrak{p}_{_{P}}^{\triangleleft}(A)\text{ for all }A\in\mathfrak{M.}%
\]

\end{theorem}

\begin{proof}
(i) To prove that $\mathfrak{p}_{_{P}}^{\triangleleft}$ is a topological
radical we have to verify all conditions $(1^{\ast}1)-(4^{\ast}).$

$(1^{\ast}1).$ Let $\phi$: $A\rightarrow B$ be an isomorphism. Set
$I=\phi(\mathfrak{p}_{_{P}}^{\triangleleft}(A)).$ If $I\nsubseteq
\mathfrak{p}_{_{P}}^{\triangleleft}(B)$ then, by Corollary \ref{C4.2n}(ii),
there is $J\in$ Id$_{B}$ such that $I\neq J$ and $J\ll_{_{P}}I.$ Hence $I/J\in
P.$ As $\phi^{-1}$ is an isomorphism, $\phi^{-1}(I)/\phi^{-1}(J)\approx I/J.$
So $\phi^{-1}(I)/\phi^{-1}(J)\in P.$ Hence $\phi^{-1}(J)\ll_{_{P}}\phi
^{-1}(I)=\mathfrak{p}_{_{P}}^{\triangleleft}(A)$ and $\phi^{-1}(I)\neq
\phi^{-1}(J)$ which contradicts Corollary \ref{C4.2n}(ii). Thus $\phi
(\mathfrak{p}_{_{P}}^{\triangleleft}(A))\subseteq\mathfrak{p}_{_{P}%
}^{\triangleleft}(B).$

Similarly, $\phi^{-1}(\mathfrak{p}_{_{P}}^{\triangleleft}(B))\subseteq
\mathfrak{p}_{_{P}}^{\triangleleft}(A),$ so that $\phi(\mathfrak{p}_{_{P}%
}^{\triangleleft}(A))=\mathfrak{p}_{_{P}}^{\triangleleft}(B).$\smallskip

$(1^{\ast}2).$ Let $A\in\mathfrak{M,}$ $I\in$ Id$_{A}$ and $p$: $A\rightarrow
A/I.$ Set $K=\mathfrak{p}_{_{P}}^{\triangleleft}(A)$ in (\ref{3.17}). Then,
for each ideal $\widehat{J}$ of $p(\mathfrak{p}_{_{P}}^{\triangleleft}(A)),$
there is $J\in$ Id$_{A}$ such that $J\subseteq\mathfrak{p}_{_{P}%
}^{\triangleleft}(A)$ and $p(\mathfrak{p}_{_{P}}^{\triangleleft}%
(A))/\widehat{J}\approx\mathfrak{p}_{_{P}}^{\triangleleft}(A)/J$.

Let $\mathfrak{p}_{_{P}}^{\triangleleft}(A)/J$ $\neq\{0\}.$ As $\mathfrak{p}%
_{_{P}}^{\triangleleft}(A)\in\mathsf{dNG}P$ by Theorem \ref{T3.2},
$\mathfrak{p}_{_{P}}^{\triangleleft}(A)/J$ is not a $P$-algebra. So
$p(\mathfrak{p}_{_{P}}^{\triangleleft}(A))/\widehat{J}$ is not a a $P$-algebra
for each ideal $\widehat{J}$ of $p(K).$ Thus, by Definition \ref{D3.2},
$p(\mathfrak{p}_{_{P}}^{\triangleleft}(A))$ is a \textsf{d}$\mathsf{NG}%
P$-ideal of $A/I.$ As $\mathfrak{p}_{_{P}}^{\triangleleft}(A/I)$ is the
largest \textsf{d}$\mathsf{NG}P$-ideal of $A/I$ by Theorem \ref{T3.2},
$p(\mathfrak{p}_{_{P}}^{\triangleleft}(A))\subseteq\mathfrak{p}_{_{P}%
}^{\triangleleft}(A/I)=\mathfrak{p}_{_{P}}^{\triangleleft}(p(A))$.$\smallskip$

$(2^{\ast}).$ By Theorem \ref{T3.2}, $A/\mathfrak{p}_{_{P}}^{\triangleleft
}(A)$ is a $\mathsf{dG}P$-algebra, so that $\mathfrak{p}_{_{P}}^{\triangleleft
}(A/\mathfrak{p}_{_{P}}^{\triangleleft}(A))=\{0\}.\smallskip$

$(3^{\ast}).$ By Theorem \ref{T3.2}, $\mathfrak{p}_{_{P}}^{\triangleleft
}(A)\in\mathsf{dNG}P,$ so that $\mathfrak{p}_{_{P}}^{\triangleleft
}(\mathfrak{p}_{_{P}}^{\triangleleft}(A))=\mathfrak{p}_{_{P}}^{\triangleleft
}(A).\smallskip$

$(4^{\ast}).$ Let $I\in$ Id$_{A}.$ Set $J=\mathfrak{p}_{_{P}}^{\triangleleft
}(I).$ If $J\nsubseteqq\mathfrak{p}_{_{P}}^{\triangleleft}(A)$ then, by
Corollary \ref{C4.2n}(ii), $J$ has a $\ll_{_{P}}$-predecessor, i.e.,
$\mathfrak{p}_{_{P}}^{\triangleleft}(I)$ has a $\ll_{_{P}}$-predecessor. This
contradicts Corollary \ref{C4.2n}(ii). So $\mathfrak{p}_{_{P}}^{\triangleleft
}(I)\subseteq\mathfrak{p}_{_{P}}^{\triangleleft}(A).\smallskip$

Thus $\mathfrak{p}_{_{P}}^{\triangleleft}$ is a topological radical. By
Theorem \ref{T3.2} and $(\ref{3,4})$,%
\begin{align*}
\mathbf{Rad}\left(  \mathfrak{p}_{_{P}}^{\triangleleft}\right)   &  =\left\{
A\in\mathfrak{M}\text{: }\mathfrak{p}_{_{P}}^{\triangleleft}(A)=A\right\}
=\mathsf{dNG}P\cap\mathfrak{M,}\\
\mathbf{Sem}\left(  \mathfrak{p}_{_{P}}^{\triangleleft}\right)   &  =\left\{
A\in\mathfrak{M}\text{: }\mathfrak{p}_{_{P}}^{\triangleleft}(A)=\{0\}\right\}
=\mathsf{dG}P\cap\mathfrak{M}.
\end{align*}

(ii) Set $P=\mathbf{Sem}\left(  R\right)  .$ By Proposition \ref{hr}, $P$ is
lower stable and $P\subseteq\mathfrak{M}$.

Let $A\in\mathfrak{M}.$ By (3*) and (\ref{3,4}), $R(A)\in\mathbf{Rad}(R).$ As
$\mathbf{Rad}(R)$ is a upper stable C*-property by Proposition \ref{hr}, we
have%
\begin{equation}
R(A)/I\in\mathbf{Rad}(R)\text{ for each ideal }I\subseteq R(A). \label{3.18}%
\end{equation}
If $R(A)\nsubseteq\mathfrak{p}_{_{P}}^{\triangleleft}(A)$ then, by Corollary
\ref{C4.2n}(ii), there is $I\in$ Id$_{A}$ such that $I\neq R(A)$ and
$I\ll_{_{P}}R(A).$ Hence $R(A)/I\in P=\mathbf{Sem}\left(  R\right)  .$ On the
other hand, $R(A)/I\in\mathbf{Rad}(R)$ by (\ref{3.18}). Since $\mathbf{Rad}%
(R)\cap\mathbf{Sem}\left(  R\right)  =\{0\},$ we have $R(A)/I=\{0\}.$ So
$I=R(A),$ a contradiction. Thus $R(A)\subseteq\mathfrak{p}_{_{P}%
}^{\triangleleft}(A).$

Suppose that $R(A)\subsetneqq\mathfrak{p}_{_{P}}^{\triangleleft}(A).$ As
$\mathfrak{p}_{_{P}}^{\triangleleft}(A)/R(A)$ is an ideal of $A/R(A)$ and as
$R\left(  A/R\left(  A\right)  \right)  =\{0\}$ by (2*), it follows from (4*)
that $R(\mathfrak{p}_{_{P}}^{\triangleleft}(A)/R(A))=\{0\}.$ Hence
$\mathfrak{p}_{_{P}}^{\triangleleft}(A)/R(A)\in\mathbf{Sem}\left(  R\right)
=P$ (see (\ref{3,4})), so that $R(A)\ll_{_{P}}\mathfrak{p}_{_{P}%
}^{\triangleleft}(A)$ which contradicts Corollary \ref{C4.2n}. Thus
$R(A)=\mathfrak{p}_{_{P}}^{\triangleleft}(A).$

Clearly $P\subseteq$ $\mathsf{dG}P.$ Thus $P\subseteq$ $\mathsf{dG}%
P\cap\mathfrak{M}.$ Conversely, let $A\in$ $\mathsf{dG}P\cap\mathfrak{M}.$
Suppose that $R(A)\neq\{0\}.$ As $A$ is a $\mathsf{dG}P$-algebra, there is
$I\in$ Id$_{A}$ such that $I\subsetneqq R(A)$ and $R(A)/I\in P=\mathbf{Sem}%
\left(  R\right)  $ (see (\ref{3.22})). Since $R(A)/I\in\mathbf{Rad}(R)$ by
(\ref{3.18}) and since $\mathbf{Rad}(R)\cap\mathbf{Sem}\left(  R\right)
=\{0\},$ we have $R(A)/I=\{0\}.$ So $I=R(A),$ a contradiction. Thus
$R(A)=\{0\},$ so that $A\in P=\mathbf{Sem}\left(  R\right)  .$ Thus
$\mathsf{dG}P\cap\mathfrak{M}=P.\bigskip$
\end{proof}

Denote by $\mathcal{P}(\mathfrak{A})$ the class of all C*-properties in
$\mathfrak{A}$ and by $\mathcal{R}(\mathfrak{M})$ the set of all topological
radicals on a lower and upper stable C*-property $\mathfrak{M}.$ Consider the
following subsets of $\mathcal{P(}\mathfrak{A)}$:%
\[
\mathcal{P}_{\text{up}}(\mathfrak{M})=\{P\in\mathcal{P}(\mathfrak{A})\text{:
}P=\mathsf{G}P\cap\mathfrak{M}\}\text{ and }\mathcal{P}_{\text{lo}%
}(\mathfrak{M})=\{P\in\mathcal{P}(\mathfrak{A})\text{: }P=\text{ \textsf{d}%
}\mathsf{G}P\cap\mathfrak{M}\}.
\]
Theorems \ref{Y} and \ref{T3,1} yield

\begin{corollary}
\label{C3.5}\emph{(i) }The map $\mathcal{\theta}$\emph{: }$R\rightarrow
\mathbf{Rad}(R)$ is a one-to-one map from $\mathcal{R}(\mathfrak{M})$ onto
$\mathcal{P}_{\text{\emph{up}}}(\mathfrak{M}).$

The map $\theta^{-1}$\emph{: }$P\rightarrow\mathfrak{r}_{_{P}}^{\triangleright
}$ from $\mathcal{P}_{\text{\emph{up}}}(\mathfrak{M})$ onto $\mathcal{R}%
(\mathfrak{M})$ is its inverse.\smallskip

\emph{(ii) }The map $\varphi$\emph{: }$R\rightarrow\mathbf{Sem}(R)$ is a
one-to-one map from $\mathcal{R}(\mathfrak{M})$ onto $\mathcal{P}%
_{\text{\emph{lo}}}(\mathfrak{M}).$

The map $\varphi^{-1}$\emph{: }$P\rightarrow\mathfrak{p}_{_{P}}^{\triangleleft
}$ from $\mathcal{P}_{\text{\emph{lo}}}(\mathfrak{M})$ onto $\mathcal{R}%
(\mathfrak{M})$ is its inverse.
\end{corollary}

\subsection{Relation-valued functions}

We say that a map $f$ on a lower and upper stable C*-property $\mathfrak{M}$
is a \textit{relation-valued function }if $f(A)=$ $\ll^{A}$ is a reflexive
relation in Id$_{A}$ for each $A\in\mathfrak{M}.$ Each such function $f$
defines the following subclass of C*-algebras in $\mathfrak{M}$:%
\begin{equation}
P_{f}=\{A\in\mathfrak{M}\text{: }\{0\}\ll^{A}A,\text{ where }\ll
^{A}=f(A)\}\text{.} \label{3.15}%
\end{equation}

On the other hand, each C*-property $P$ in $\mathfrak{A}$ generates a
relation-valued function $f_{_{P}}$ on $\mathfrak{M}$:%
\begin{equation}
f_{_{_{P}}}(A)=\text{ }\ll_{_{P}}^{A},\text{ for }A\in\mathfrak{M},\text{
where }I\ll_{_{P}}^{A}J\text{ if }I\subseteq J\text{ in Id}_{A}\text{ and
}J/I\in P. \label{5.2}%
\end{equation}
We consider now the conditions for a relation-valued function to be generated
by a C*-property$.$

\begin{lemma}
\label{L3.5}Let $f$ be a relation-valued function on a lower and upper stable
C*-property $\mathfrak{M}.$ Then $f=f_{_{_{P}}}$ for some C*-property $P,$ if
and only if the following conditions hold\emph{:\smallskip}

$\qquad(C_{1})$ Let $A,B\in\mathfrak{M}$ and $\varphi$: $A\rightarrow B$ be a
$^{\ast}$-isomorphism. Then%
\[
I\ll^{A}J\text{ in }\emph{Id}_{A}\text{ if and only if }\varphi(I)\ll
^{B}\varphi(J)\text{ in }\emph{Id}_{B};
\]

$\qquad(C_{2})$ Let $I\subseteq J\subseteq K$ in \emph{Id}$_{A}.$ Then
$I\ll^{K}J\Longleftrightarrow I\ll^{A}J;\smallskip$

$\qquad(C_{3})$ Let $I\subseteq J\subseteq K$ in \emph{Id}$_{A}.$ Then
$J/I\ll^{A/I}K/I\Longleftrightarrow J\ll^{A}K.\smallskip$

In this case $P=P_{f},$ so that $\ll^{A}=$ $\ll_{_{P_{f}}}^{A}$ for all
$A\in\mathfrak{M}.$
\end{lemma}

\begin{proof}
Let $f=f_{_{_{P}}}$ for some C*-property $P,$ i.e., $f(A)=$ $\ll^{A}=$
$\ll_{_{P}}^{A}$ for all $A\in\mathfrak{A}.$

Let $\varphi$: $A\rightarrow B$ be a *-isomorphism and $I\ll^{A}J$ for
$I,J\in$ Id$_{A}.$ Then $J/I\in P.$ Since $\varphi(J)/\varphi(I)\approx J/I$,
we have $\varphi(J)/\varphi(I)\in P.$ So $\varphi(I)\ll^{B}\varphi(J).$
Similarly, as $\varphi^{-1}$: $B\rightarrow A$ is a *-isomorphism,
$\varphi(I)\ll^{B}\varphi(J)\Rightarrow I\ll^{A}J.$ Thus $(C_{1})$ holds.

Let $I\subseteq J\subseteq K$ in Id$_{A}.$ Then $(C_{2})$ and $(C_{3})$ hold,
since
\begin{align*}
I  &  \ll^{K}J\Longleftrightarrow J/I\in P\Longleftrightarrow I\ll^{A}J\text{
and}\\
J/I  &  \ll^{A/I}K/I\Longleftrightarrow K/J\approx(K/I)/(J/I)\in
P\Longleftrightarrow K/J\in P\Longleftrightarrow J\ll^{A}K.
\end{align*}

Conversely, let $f$ satisfy conditions $(C_{1})$-$(C_{3}).$ Let $A\in P_{f}$
and $\varphi$: $A\rightarrow B$ be a *-isomorphism$.$ By (\ref{3.15}),
$\{0\}_{A}\ll^{A}A$ and, by $(C_{1}),$ $\{0\}_{B}=\varphi(\{0\}_{A})\ll
^{B}\varphi(A)=B.$ So $B\in P_{f}.$ Thus $P_{f}$ is a C*-property in
$\mathfrak{A}.$

Let us show that $f=f_{_{P_{f}}},$ i.e., $\ll^{A}=$ $\ll_{_{P_{f}}}^{A}$ for
$A\in\mathfrak{M}.$ For $I,J\in$ Id$_{A},$%
\[
I\ll_{_{P_{f}}}^{A}J\overset{(\ref{4.1})}{\Leftrightarrow}J/I\in
P_{f}\overset{(\ref{3.15})}{\Leftrightarrow}\{0\}\ll^{J/I}J/I\overset{(C_{3}%
)}{\Leftrightarrow}I\ll^{J}J\overset{(C_{2})}{\Leftrightarrow}I\ll^{A}J.
\]
So $\ll^{A}=$ $\ll_{_{P_{f}}}^{A}.$ Thus $f=f_{_{P_{_{f}}}}.$

If $f=f_{_{P}}$ for some C*-property $P,$ then, by the above argument,
conditions $(C_{1})$-$(C_{3})$ hold, so that $f=f_{P_{_{f}}}.$ Hence $f_{_{P}%
}=f_{_{P_{f}}}.$ Thus $P=P_{f}.\bigskip$
\end{proof}

Let $f$ be a relation-valued function on a lower and upper stable C*-property
$\mathfrak{M}\subseteq\mathfrak{A}$. Suppose that each relation $f(A)=$
$\ll^{A}$ in Id$_{A}$ has a unique $\ll^{A}$-radical $\mathfrak{r}(A)$ (see
(\ref{2.11})):
\begin{equation}
\{0\}\ll^{A}\mathfrak{r}(A)\text{ }\overleftarrow{\ll^{A}}A. \label{4.10}%
\end{equation}
It was shown in \cite{KST1} that the radical $\mathfrak{r}(A)$ may exist even
if the relation $\ll^{A}$ is not an $\mathbf{R}$-order. We will see further in
Remark \ref{R3} that even if each $\ll^{A}$ is an $\mathbf{R}$-order in
Id$_{A},$ the map $\mathfrak{r}$: $A\in\mathfrak{M}\mapsto\mathfrak{r}\left(
A\right)  $ is not necessarily a topological radical. Hence the following
question arises: \textit{Under what conditions this map }$\mathfrak{r}%
$\textit{ is a topological radical}?

\begin{theorem}
\label{T4.5}Let $f$ be a relation-valued function on $\mathfrak{M}$ and each
relation $f(A)=$ $\ll^{A}$ in \emph{Id}$_{A},$ $A\in\mathfrak{M,}$ has a
unique $\ll^{A}$-radical $\mathfrak{r}(A),$ so that \emph{(\ref{4.10}) }holds.

Then the map $\mathfrak{r}$\emph{: }$A\in\mathfrak{M}\mapsto\mathfrak{r}%
\left(  A\right)  $ is a topological radical if and only if $f=f_{_{P_{f}}},$
where $P_{f}$\emph{ }is an upper stable C*-property and $P_{f}=\mathsf{G}%
P_{f}\cap\mathfrak{M}.$ Moreover\emph{, }in this case $P_{f}=\mathbf{Rad}%
(\mathfrak{r).}$
\end{theorem}

\begin{proof}
Let the map $\mathfrak{r}$ be a topological radical on $\mathfrak{M.}$ By
Theorem \ref{Y}(ii) $($see $(\ref{3,4}))$, the class $P=\mathbf{Rad}%
(\mathfrak{r)=\{}A\in\mathfrak{M}$: $\mathfrak{r}\left(  A\right)  =A\}$ is an
upper stable C*-property, \textsf{G}$P\cap\mathfrak{M}=P$ and $\mathfrak{r}%
(A)=\mathfrak{r}_{_{P}}^{\triangleright}(A)$ for each $A\in\mathfrak{M,}$
where $\mathfrak{r}_{_{P}}^{\triangleright}(A)$ is the $\ll_{_{P}%
}^{\triangleright}$-radical in Id$_{A}.$

If $A\in P$ then $\mathfrak{r}(A)=A$ by (\ref{3,4}). Hence, by (\ref{4.10}),
$\{0\}\ll^{A}A.$ So $A\in P_{f}.$ Thus $P\subseteq P_{f}.$ Let now $A\in
P_{f}.$ Then $\{0\}\ll^{A}A$ $\overleftarrow{\ll^{A}}A.$ As $\mathfrak{r}(A)$
is a unique ideal satisfying (\ref{4.10}), $\mathfrak{r}(A)=A.$ So $A\in P.$
Thus $P_{f}\subseteq P.$ So $P=P_{f}.$ Hence $P_{f}=\mathbf{Rad}%
(\mathfrak{r)}$\emph{ }is an upper stable C*-property and $P_{f}%
=\mathsf{G}P_{f}\cap\mathfrak{M}.$

It follows from Lemma \ref{L3.5} that $\ll^{A}=$ $\ll_{_{P_{f}}}^{A}$ for all
$A\in\mathfrak{M},$ so that $f=f_{_{P_{f}}}.$

Conversely, let $P_{f}$\emph{ }be an upper stable C*-property and
$P_{f}=\mathsf{G}P_{f}\cap\mathfrak{M}.$ Set $Q=P_{f}.$ By Corollary
\ref{C3.3}, the relation $\ll_{_{Q}}^{A}=$ $\left(  \ll_{_{Q}}^{A}\right)
^{\triangleright}$ is an $\mathbf{R}$-order in Id$_{A}$ for each
$A\in\mathfrak{M},$ so that $\mathfrak{r}_{_{Q}}(A)=\mathfrak{r}_{_{Q}%
}^{\triangleright}(A)$ is the $\ll_{_{Q}}^{A}$-radical in Id$_{A}.$ By
(\ref{2.11}),%
\begin{equation}
\{0\}\ll_{_{Q}}^{A}\mathfrak{r}_{_{Q}}(A)\text{ }\overleftarrow{\ll_{_{Q}}%
^{A}}A. \label{4.12}%
\end{equation}
As $\mathfrak{r}_{_{Q}}(A)=\mathfrak{r}_{_{Q}}^{\triangleright}(A),$ the map
$\mathfrak{r}_{_{Q}}$: $A\in\mathfrak{M}\rightarrow\mathfrak{r}_{_{Q}}(A)$ is
a topological radical on $\mathfrak{M}$ by Theorem \ref{Y}(i).

We have from Lemma \ref{L3.5} that $\ll_{_{Q}}^{A}=$ $\ll^{A}.$ Hence, by
(\ref{4.12}), $\{0\}\ll^{A}\mathfrak{r}_{_{Q}}(A)$ $\overleftarrow{\ll^{A}}A.$
Comparing this to (\ref{4.10}) and taking into account that $\mathfrak{r}(A)$
is a unique $\ll^{A}$-radical, we get $\mathfrak{r}(A)=\mathfrak{r}_{_{Q}}(A)$
for $A\in\mathfrak{M}.$ As $\mathfrak{r}_{_{Q}}$ is a topological radical, the
map $\mathfrak{r}$\emph{: }$A\in\mathfrak{M}\mapsto\mathfrak{r}\left(
A\right)  $ is a topological radical.
\end{proof}

\section{Small ideals in C*-algebras}

In this section we study small ideals of C*-algebras $A$ and the relation
$\ll_{\text{sm}}=$\ $\ll_{\text{sm}}^{A}$ they generate in Id$_{A}$. Unlike
the relations constructed in the previous sections from various C*-properties
in $\mathfrak{A}$, this relation is defined using the property of ideals
vis-\`{a}-vis the algebra $A.$

We show that it is an $\mathbf{H}$-relation in each Id$_{A},$ so that the
relation $\ll_{\text{sm}}^{\triangleright}$ is an $\mathbf{R}$-order and
Id$_{A}$ has the $\ll_{\text{sm}}^{\triangleright}$-radical $\mathfrak{r}%
_{\text{sm}}^{\triangleright}\left(  A\right)  .$ However, unlike the radicals
$\mathfrak{r}_{_{P}}^{\triangleright}$ generated by C*-properties $P$, the map
$\mathfrak{r}_{\text{sm}}^{\triangleright}$:\emph{ }$A\mapsto\mathfrak{r}%
_{\text{sm}}^{\triangleright}\left(  A\right)  $\emph{ }is not a topological
radical on\emph{ }$\mathfrak{A}.$

If $A$ is unital then $\ll_{\text{sm}}$ $=$ $\ll_{\text{sm}}^{\triangleright}$
is an $\mathbf{R}$-order and the $\ll_{\text{sm}}$-radical $\mathfrak{r}%
_{\text{sm}}\left(  A\right)  $ in Id$_{A}$ coincides with the radical
rad$_{_{K}}$($A$) introduced by Kasch in \cite[p. 214]{Kas}.

\begin{definition}
Let $A\neq\{0\}.$ An $I\in$ \emph{Id}$_{A}$ is \textbf{small} if $I+K\neq A$
for all $K\in$ \emph{Id}$_{A},$ $A\neq K.$ Otherwise\emph{,} $I$ is
\textbf{non-small}. If $A=\{0\}$ then $I=A=\{0\}$ is a \textbf{small} ideal of
$A.$
\end{definition}

\begin{example}
\label{E1}\emph{1) If }$A\neq\{0\}$ \emph{then }$\{0\}$ \emph{is a small ideal
of }$A$\emph{$,$ and} $A$ \emph{is a non-small ideal.\smallskip}

\emph{2) For }$B(H),$\emph{ Id}$_{B(H)}=(\{0\},C(H),B(H)\}$\emph{ and
}$\{0\},C(H)$\emph{ are small ideal.\smallskip}

\emph{3) Let} $A=A_{\infty}$ \emph{in Proposition \ref{P5.5}. Then Id}%
$_{A}=\{A_{n}\}_{n=0}^{\infty},$ $A_{n}\subset A_{n+1}$ \emph{for all }%
$n,$\emph{ and} $A=\overline{\cup A_{n}}.$ \emph{Clearly, each} $A_{n}$
\emph{is a small ideal in} $A.$\emph{\smallskip}

\emph{4) If }$A$\emph{ is dual, then each }$\{0\}\neq I\in$\emph{ Id}$_{A}%
$\emph{ is non-small, as }$I\dotplus$\emph{ an}$(I)=A$\emph{ and an(}$I)\neq
A.$\emph{\smallskip}

\emph{5) Let }$A=C(X)$\emph{ be the C*-algebra of all continuous functions on
a compact }$X$\emph{ and }$\{0\}\neq I\in$\emph{ Id}$_{A}.$\emph{ There is a
compact }$Y\subsetneqq X$\emph{ such that }$I=\{f\in A$\emph{: }$f(y)=0$\emph{
for all }$y\in Y\}.$\emph{ For }$x\in X$\emph{ and }$x\neq Y,$\emph{ let
}$K=\{f\in A$\emph{: }$f(x)=0\}.$\emph{ Then }$I+K=A.$\emph{ Thus all non-zero
ideals of }$A$\emph{ are non-small.}
\end{example}

\begin{theorem}
\label{topol}The ideal $\{0\}$ is the only small ideal of a \emph{C*}-algebra
$A$ if and only if the space $Prim(A)$ of all primitive ideals of $A$ has the
following property\emph{:\smallskip}

$(\Omega)$ Each non-void open subset of $Prim(A)$ contains a non-void closed subset.
\end{theorem}

\begin{proof}
Suppose that $Prim(A)$ has property ($\Omega).$ For $\{0\}\neq K\in$ Id$_{A},$
let $D=hull(K)$ be the set of all primitive ideals containing $K$. Then $D\neq
Prim(A),$ as $K\neq\{0\}$. So $Q=Prim(A)\diagdown D\neq\varnothing$ is an open
set. Let $M\neq\varnothing$ be a closed subset of $Q$. By definition, there is
$J\in$ Id$_{A}$ with $hull(J)=M$. If $J+K\neq A,$ there is $I\in Prim(A)$ that
contains $J+K$. So $K\subset I$ and $J\subset I$. Hence $I\in hull(J)\cap
hull(K)=M\cap D=\varnothing$, a contradiction. Thus $J+K=A$, so that $K$ is
not small.

Conversely, let $A$ have no small ideals apart from $\{0\}.$ For an open
subset $Q$ of $Prim(A),$ let $K=\cap\{I$: $I\in Prim(A)\diagdown Q\}$. Then
$hull(K)=Prim(A)\diagdown Q,$ so that $K\neq\{0\}.$ By our assumption, there
is $J\in$ Id$_{A}$ with $A=J+K$. This is only possible if $hull(J)\cap
hull(K)=\varnothing$. So $hull(J)\subset Q$. Thus $Q$ contains a non-void
closed subset.
\end{proof}

\begin{corollary}
\label{clospoint}If closed points are dense in $Prim(A)$ then $\{0\}$ is the
only small ideal in $A.$
\end{corollary}

\begin{corollary}
\label{CCR}Let a \emph{C*}-property $S$ consist of simple \emph{C*}%
-algebras$.$ In each $\mathsf{R}(S)$-algebra $A$ \emph{(}in particular\emph{,
}in each CCR-algebra\emph{)} \emph{(}Definition \emph{\ref{D4.1}), }$\{0\}$ is
the only small ideal$.$
\end{corollary}

\begin{proof}
Let $I=\ker\pi\in Prim(A),$ $\pi\in\Pi(A).$ Then $A/I\approx\pi(A)\in S$ is a
simple C*-algebra whence $I$ is a maximal ideal of $A.$ So $I$ is a closed
point in $Prim(A).$\bigskip
\end{proof}

Let $I,J,K\in$ Id$_{A}$ with $I\subset J\cap K$. It is easy to see that
\begin{align}
(J+K)/I  &  =J/I+K/I\text{ and}\label{1'}\\
J/I  &  =K/I\text{ implies }J=K. \label{2'}%
\end{align}

Let $A,B\in\mathfrak{A}$ and $\varphi$ be an isomorphism from $A$ onto $B.$
Note that
\begin{equation}
I\in\text{ Id}_{A}\text{ is small in }A\text{ if and only if }\varphi(I)\text{
is small in }B. \label{6.5}%
\end{equation}

\begin{lemma}
\label{L1}For $I\in$ \emph{Id}$_{A},$ let $p_{_{I}}$\emph{: }$A\rightarrow
A/I$ be the standard epimorphism.\smallskip

\emph{(i) \ \ }If $J$ is a small ideal in $A$ then $p(J)$ is small in
$A/I.\smallskip$

\emph{(ii) \ }Let $I\subset J$ in \emph{Id}$_{A}.$Then $J$ is small in $A,$ if
and only if $I$ is small in $A$ and $J/I$ is small in $A/I.\smallskip$

\emph{(iii) }Let $I\subset J$ in \emph{Id}$_{A}.$ If $I$ is small in $J$ then
$I$ is small in $A.\smallskip$

\emph{(iv) }Any finite sum of small ideals is a small ideal.\smallskip
\end{lemma}

\begin{proof}
(i) Let $p(J)$ be not small in $A/I.$ Then $p(J)+\widehat{K}=A/I$ for some
ideal $\widehat{K}\neq A/I$ in $A/I.$ Let $K$ be the preimage of $\widehat{K}$
in $A.$ Then $K\neq A$ and, for each $a\in A,$ there is $j\in J$ such that
$p_{_{I}}(a)-p_{_{I}}(j)\in\widehat{K}.$ Hence $a-j\in K,$ so that $A=J+K.$ So
$J$ is not small in $A$, a contradiction. Thus $p_{_{I}}(J)$ is small in
$A/I.$

(ii) Let $J$ be small in $A.$ By (i), $J/I=p_{_{I}}(J)$ is small in $A/I.$

If $I$ is not small, $I+K=A,$ $K\neq A.$ Then $J+K=A,$ a contradiction. So $I$
is small.

Conversely, let $I$ be small in $A$ and $J/I$ in $A/I.$ If $J$ is non-small,
$A=J+K$ for some $K\neq A.$ As $I$ is small, $K+I\neq A,$ so that $(K+I)/I\neq
A/I$ by (\ref{2'}). By (\ref{1'}), $J/I+(K+I)/I=(J+K+I)/I=A/I.$ Hence $J/I$ is
not small, a contradiction. Thus $J$ is small.\smallskip

(iii) Assume that $I$ is not small in $A.$ Then $I+W=A$, for some $W\in$
Id$_{A}$, $W\neq A$. It follows that $I+W\cap J=J.$ Indeed, for each $a\in J,$
we have $a=m+w$, $m\in I,w\in W$. Hence $w\in W\cap J$, so that $a\in I+W\cap
J$. Thus $J=I+W\cap J$. As $I$ is a small ideal in $J$, we have $W\cap J=J$.
So $J\subset W$ and $A=I+W=W$, a contradiction.\smallskip

(iv) Let $I$ and $J$ be small ideals. If $I+J$ is non-small, $(I+J)+K=A$ for
some $K\neq A$. If $J+K=A$ then $J$ is non-small, otherwise $I$ is non-small
-- a contradiction. Thus $I+J$ is small.\bigskip
\end{proof}

For a C*-algebra $A,$ denote by $\mathcal{S}(A)$ the set of all small ideals
in Id$_{A}.$ Set
\begin{equation}
S_{\Lambda}=\overline{\sum_{I\in\Lambda}I},\text{ for each subset }%
\Lambda\subseteq\mathcal{S}(A),\text{ and }S_{A}=\overline{\sum_{I\in
\mathcal{S}(A)}I}, \label{6.3}%
\end{equation}
where $\sum_{I\in\Lambda}I$ is the set of all finite sums of elements from the
ideals $I\in\Lambda.$

\begin{lemma}
\label{union}If $A$ is unital then $S_{\Lambda}$ is a small ideal for each
$\Lambda\subseteq\mathcal{S}(A).$
\end{lemma}

\begin{proof}
If $S_{\Lambda}$ is not small then $K+S_{\Lambda}=A$ for some $K\neq A$. Hence
$K+\sum_{I\in\Lambda}I$ is a dense ideal in $A$. As $\mathbf{1}_{A}$\textbf{
}is surrounded by a ball of invertible elements, $A$ has no dense ideals. Thus
$K+\sum_{I\in\Lambda}I=A,$ so that $a+b=\mathbf{1}_{A}$ for some $a\in K$ and
$b=\sum_{i=1}^{n}x_{i}$ for some $x_{i}\in I_{i}.$ By Lemma \ref{L1}, the
ideal $\sum_{i=1}^{n}I_{i}$ is small. As $a+b=\mathbf{1}_{A},$ we have
$K+\sum_{i=1}^{n}I_{i}=A$ -- a contradiction.
\end{proof}

\begin{remark}
\label{R1'}\emph{If }$A$ \emph{is} not unital$,$ $S_{A}$\emph{ is not
necessarily small: }$S_{A}=\overline{\sum_{n}A_{n}}=A$ \emph{in Example
\ref{E1} 3.}
\end{remark}

In general, we have the following result.

\begin{proposition}
\label{P7.3}Let a \emph{C*}-property $S$ consist of simple \emph{C*}%
-algebras$.$ If a \emph{C*}-algebra $A$ has a $\mathsf{R}(S)$-quotient
\emph{(}for example\emph{,} $A\in RFD$ or $A\in\mathsf{dGR}(S)$ and\emph{, }in
particular\emph{, }$A$ is a $\mathsf{d}GCR$-algebra$)$ then $S_{A}\neq A.$
\end{proposition}

\begin{proof}
Suppose that $A/I$ is a $\mathsf{R}(S)$-algebra. If $J\in$ Id$_{A}$ is small
then $p_{_{I}}(J)$ is a small in $A/I$ by Lemma \ref{L1}$.$ By Corollary
\ref{CCR}, this means that $p_{_{I}}(J)=0$. So $J\subset I$. Thus all small
ideals of $A$ are contained in $I$. Hence $S_{A}\subseteq I$ and $S_{A}\neq A$.

A $RFD$-algebra $A$ has $\pi\in\Pi(A)$ such that $A/\ker\pi\approx
C(\mathcal{H}_{\pi})\in\mathsf{R}(C)=CCR$ (see (\ref{3.30})). For a
$\mathsf{dGR}(S)$-algebra (in particular, $\mathsf{d}GCR$-algebra) $A$, there
is $I\in$ Id$_{A}$ such that $A/I$ is a $\mathsf{R}(S)$-algebra by Definition
\ref{D3.2}.
\end{proof}

\begin{corollary}
\label{maximal}\emph{(i) }If $S_{A}$ is a small ideal\emph{, }it is the
largest small ideal in $A$ and $A/S_{A}$ has no non-zero small
ideals.\smallskip

\emph{(ii) }If $A$ has a maximal small ideal $M$ then $M=S_{A}$. \smallskip

\emph{(iii) }If $A$ is unital then $S_{A}$ is the largest small ideal in $A$
and $A/S_{A}$ has no small ideals.
\end{corollary}

\begin{proof}
(i) If $J\nsubseteqq S_{A}$ is a small ideal, the ideal $S_{A}+J$ is small by
Lemma \ref{union} and larger than $S_{A}.$ This contradiction shows that
$S_{A}$ is the largest small ideal in $A.$

If $W\neq\{0\}$ is a small ideal in $A/S_{A}$ then, by Proposition \ref{L1},
its preimage in $A$ is a small ideal larger than $S_{A}$, a contradiction. So
$A/S_{A}$ has no small ideals.

(ii) By (\ref{6.3}), $M\subseteq S_{A}.$ If $J\nsubseteqq M$ is a small ideal,
the ideal $M+J$ is small by Lemma \ref{union} and larger than $M.$ This
contradiction shows that $J\subseteq M.$ By (\ref{6.3}), $M=S_{A}$

(iii) By Lemma \ref{union}, $S_{A}$ is a small ideal. The rest follows from
(i).\bigskip
\end{proof}

Let $I,J,K\in$ Id$_{A}$ and $I\cap J=\{0\}$. Then%
\begin{equation}
(I\dotplus J)\cap K\overset{(\ref{distrib})}{=}(I\cap K)+(J\cap K)=(I\cap
K)\dotplus(J\cap K). \label{6.4}%
\end{equation}

\begin{theorem}
\label{T6.1}Let $A\in\mathfrak{A}$ and let the ideal $S_{A}\neq A$ be not
small\emph{: }$A=S_{A}+R$ for some $R\in$ \emph{Id}$_{A}.$\smallskip

\emph{(i) }If $A=S_{A}\dotplus R$ is the direct sum\emph{ }then $A/S_{A}$ has
no non-trivial small ideals.\smallskip

\emph{(ii) }If each non-trivial ideal in $S_{A}$ is small then $A/S_{A}$ has
no non-trivial small ideals.
\end{theorem}

\begin{proof}
(i) Let $A=S_{A}\dotplus R$ and $S_{A}\subsetneqq S\in$ Id$_{A}.$ By
(\ref{6.4}), $S=S_{A}\dotplus(R\cap S),$ so that $S/S_{A}\approx R\cap S$ and
$A/S_{A}\approx R.$ If $S/S_{A}$ is a small ideal in $A/S_{A}$ then $R\cap S$
is small in $R$ by (\ref{6.5}). By Lemma \ref{L1}(iii), $R\cap S$ is a small
ideal in $A.$ As it is not contained in $S_{A},$ we get a contradiction. Thus
$A/S_{A}$ has no non-trivial small ideals.

(ii) Let $J=S_{A}\cap R\neq\{0\}.$ Then $J\neq S_{A},$ $J\neq R,$
$A/J=S_{A}/J\dotplus R/J$. If $S_{A}\subsetneqq S\in$ Id$_{A}$ then%
\begin{equation}
S/J\overset{(\ref{6.4})}{=}S_{A}/J\dotplus((R/J)\cap(S/J)). \label{6.6}%
\end{equation}
Set $T=(R/J)\cap(S/J).$ We have
\[
\{0\}\neq S/S_{A}\overset{(\ref{4.01})}{\approx}(S/J)/(S_{A}%
/J)\overset{(\ref{6.6})}{\approx}(R/J)\cap(S/J)=T\text{ and }A/S_{A}%
\overset{(\ref{4.01})}{\approx}(A/J)/(S_{A}/J)\approx R/J.
\]

Suppose that $S/S_{A}$ is a small ideal in $A/S_{A}.$ By (\ref{6.5}),
$T\approx S/S_{A}$ is small in $R/J\approx A/S_{A}.$ By Lemma \ref{L1}(ii),
$T$ is a small ideal in $A/J$ not contained in $S_{A}/J.$

Let $p_{_{J}}$:\emph{ }$A\rightarrow A/J$ and $K=p_{_{J}}^{-1}(T).$ As
$T\nsubseteqq S_{A}/J,$ the ideal $K\nsubseteqq S_{A}.$ By the condition of
the theorem, $J=S_{A}\cap R\subsetneqq K$ is a small ideal in $A,$ and
$K/J=p_{_{J}}(K)=T$ is a small ideal in $A/J.$ Hence, by Lemma \ref{L1}(i),
the ideal $K$ is small in $A$. As all small ideals of $A$ lie in $S_{A},$ we
get a contradiction. Hence $A/S_{A}$ has no non-zero small ideals.
\end{proof}

\begin{problem}
\label{P6.1}\emph{Does }$A/S_{A}$ \emph{always have no non-zero small ideals?}
\end{problem}

\begin{proposition}
\label{L2}Let $I\subset J$ in \emph{Id}$_{A}$ and $J/I$ be small in $A/I.$
Then the ideal $\widehat{R}=(J+K)/(I+K)$ is small in $A/(I+K)$ for each $K\in$
\emph{Id}$_{A},$ $J+K\neq A.$
\end{proposition}

\begin{proof}
If $\widehat{R}$ is non-small in $A/(I+K),$ there is an ideal $S$ in $A$ such
that%
\begin{equation}
I+K\subset S\neq A\text{ and }\widehat{R}+\widehat{S}=A/(I+K),\text{ where
}\widehat{S}=S/(I+K). \label{3}%
\end{equation}
So
\[
\widehat{R}+\widehat{S}=(J+K)/(I+K)+S/(I+K)\overset{(\ref{3})}{=}A/(I+K)
\]
By (\ref{1'}) and (\ref{2'}), $A=(J+K)+S.$ Hence $A=J+S,$ as $K\subseteq S.$
Thus, by (\ref{1'}), $A/I=J/I+S/I.$ By (\ref{2'}), $A/I\neq S/I,$ as $S\neq
A.$ Hence $J/I$ is non-small -- a contradiction. So $\widehat{R}$ is
small.\bigskip
\end{proof}

Let $A\in\mathfrak{A.}$ Define the relation $\ll_{\text{sm}}^{A}$ in Id$_{A}$
as follows:%
\begin{equation}
I\ll_{\text{sm}}^{A}J,\text{ if }I\subseteq J\text{ and }J/I\text{ is a small
ideal in }A/I. \label{5.3}%
\end{equation}

We write $\ll_{\text{sm}}$ instead of $\ll_{\text{sm}}^{A}$, if it is clear
what algebra $A$ we consider.

\begin{corollary}
\label{C1}The relation $\ll_{\text{\emph{sm}}}$ in \emph{Id}$_{A}$ is a
transitive $\mathbf{H}$-relation.
\end{corollary}

\begin{proof}
The relation $\ll_{\text{sm}}$ is reflexive, as $\{0\}$ is a small ideal in
each $A\in\mathfrak{A}.$ If $I\ll_{\text{sm}}J$ in Id$_{A},$ $J/I$ is small in
$A/I.$ By Proposition \ref{L2}, $(J+K)/(I+K)$ is small in $A/(I+K)$ for each
$K\in$ Id$_{A}.$ Thus $(I+K)=I\vee K\ll_{\text{sm}}(J+K)=J\vee K.$ By
(\ref{1.3}), $\ll_{\text{sm}}$ is an $\mathbf{H}$-relation.

Let $I\ll_{\text{sm}}J\ll_{\text{sm}}K$ in Id$_{A}.$ Set $I^{\prime}=J/I,$
$J^{\prime}=K/I$ and $A^{\prime}=A/I.$ Then $I^{\prime}=J/I$ is small in
$A^{\prime}=A/I$ and $K/J\overset{(\ref{4.01})}{\approx}(K/I)/(J/I)=J^{\prime
}/I^{\prime}$ is small in $A/J\overset{(\ref{4.01})}{\approx}%
(A/I)/(J/I)=A^{\prime}/I^{\prime}.$ Hence, by (\ref{6.5}), $J^{\prime
}/I^{\prime}$ is small in $A^{\prime}/I^{\prime}.$ From Lemma \ref{L1}(i) it
follows that $J^{\prime}=K/I$ is small in $A^{\prime}=A/I.$ Thus
$I\ll_{\text{sm}}K,$ so that $\ll_{\text{sm}}$ is transitive.\bigskip
\end{proof}

As $\ll_{\text{sm}}$ is an $\mathbf{H}$-relation in Id$_{A},$ it follows from
Theorem \ref{T4.1}, Corollary \ref{C4.2n} and (\ref{3.6}) that there is a
unique $\ll_{\text{sm}}^{\triangleright}$-radical $\mathfrak{r}_{\text{sm}%
}^{\triangleright}(A)=\overline{\sum J},$ where $\{0\}\ll_{\text{sm}%
}^{\triangleright}J\in$ Id$_{A}.$ Note that the ideals $J$ are not, generally
speaking, small in $A,$ since $\{0\}\ll_{\text{sm}}^{\triangleright}J$ (not
$\{0\}\ll_{\text{sm}}J).$ Corollary \ref{C4.2n} yields

\begin{corollary}
\label{C2}For $A\in\mathfrak{A,}$ let $\ll_{\text{\emph{sm}}}=$ $\ll
_{\text{\emph{sm}}}^{A}$ and $\mathfrak{r}_{\text{\emph{sm}}}^{\triangleright
}=\mathfrak{r}_{\text{\emph{sm}}}^{\triangleright}(A)$ be the unique
$\ll_{\text{\emph{sm}}}^{\triangleright}$-radical.\smallskip

\emph{(i) }For each $I\subseteq\mathfrak{r}_{\text{\emph{sm}}}^{\triangleright
},$ there is an ascending transfinite $\ll_{\text{\emph{sm}}}$-series $\left(
I_{\lambda}\right)  _{1\leq\lambda\leq\gamma}$ of ideals from $I$\ to
$\mathfrak{r}_{\text{\emph{sm}}}^{\triangleright}$ such that each
$I_{\lambda+1}/I_{\lambda}$ is a small ideal in $A/I_{\lambda}.\smallskip$

\emph{(ii) \ }If $\left(  I_{\lambda}\right)  _{1\leq\lambda\leq\gamma}$ is an
ascending transfinite $\ll_{\text{\emph{sm}}}$-series of ideals from $\{0\}$
to $I$ then $J\subseteq\mathfrak{r}_{\text{\emph{sm}}}^{\triangleright
}.\smallskip$

\emph{(iii)} $\mathfrak{r}_{\text{\emph{sm}}}^{\triangleright}$ contains all
small ideals of $A;$ the algebra $A/\mathfrak{r}_{\text{\emph{sm}}%
}^{\triangleright}$ has no non-zero small ideals$.$\smallskip

\emph{(iv) \ }If $\mathfrak{r}_{\text{\emph{sm}}}^{\triangleright}\nsubseteq
I$ then there is $J\in$ \emph{Id}$_{A}$ such that\emph{ }$I\subsetneqq J$ and
$J/I$ is small in $A/I.$
\end{corollary}

It follows from Remark \ref{R1'} that the ideal $\mathfrak{r}_{\text{sm}%
}^{\triangleright}$ is not necessarily small in $A.$

\begin{corollary}
\label{C6.1}If $A\in\mathfrak{A}$ is unital\emph{ }then $\ll_{\text{\emph{sm}%
}}^{A}$ is an $\mathbf{R}$-order\emph{, }$\ll_{\text{\emph{sm}}}^{A}$ $=$
$(\ll_{\text{\emph{sm}}}^{A})^{\triangleright}$ and $\mathfrak{r}%
_{\text{\emph{sm}}}(A)=\mathfrak{r}_{\text{\emph{sm}}}^{\triangleright}(A).$
\end{corollary}

\begin{proof}
Set $\ll_{\text{sm}}=$ $\ll_{\text{sm}}^{A}$ and $\ll_{\text{sm}%
}^{\triangleright}=$ ($\ll_{\text{sm}}^{A})^{\triangleright}.$ Clearly (see
(\ref{2.2})), $I\ll_{\text{sm}}J$ implies $I\ll_{\text{sm}}^{\triangleright
}J.$

Conversely, let $I\ll_{\text{sm}}^{\triangleright}J.$ By (\ref{2.2}), there is
an ascending transfinite $\ll_{\text{sm}}$-series $\left(  I_{\lambda}\right)
_{1\leq\lambda\leq\gamma}$ of ideals from $I$ to $J$ such that $I_{1}=I,$
$I_{\gamma}=J$ and each $I_{\lambda+1}/I_{\lambda}$ is a small ideal in
$A/I_{\lambda}.$ Let us show that all $I_{\lambda}/I\in\mathcal{S}(A/I)$ for
$1\leq\lambda\leq\gamma.$

Let $I_{\lambda}/I$ be small in $A/I$ for some $\lambda.$ Set $I^{\prime
}=I_{\lambda}/I,$ $J^{\prime}=I_{\lambda+1}/I,$ $A^{\prime}=A/I.$ By
(\ref{4.01}),%
\[
I_{\lambda+1}/I_{\lambda}\approx(I_{\lambda+1}/I)/(I_{\lambda}/I)=J^{\prime
}/I^{\prime}\text{ and }A/I_{\lambda}\approx(A/I)/(I_{\lambda}/I)=A^{\prime
}/I^{\prime}.
\]
As $I_{\lambda+1}/I_{\lambda}$ is small in $A/I_{\lambda},$ it follows from
(\ref{6.5}) that $J^{\prime}/I^{\prime}$ is small in $A^{\prime}/I^{\prime}.$
So, by Lemma \ref{L1}, $J^{\prime}=I_{\lambda+1}/I$ is small in $A^{\prime
}=A/I,$ i.e., $I_{\lambda+1}/I\in\mathcal{S}(A/I).$

Let $\beta$ be a limit ordinal and let all $I_{\lambda}/I\in\mathcal{S}(A/I)$
for $\lambda<\beta$. Then $I_{\beta}=\overline{\cup_{\lambda<\beta}I_{\lambda
}}.$ So $I_{\beta}/I=\overline{\cup_{\lambda<\beta}(I_{\lambda}/I)}%
=\overline{\sum_{\lambda<\beta}(I_{\lambda}/I)}.$ As $A/I$ is unital, it
follows from Lemma \ref{union} that $I_{\beta}/I\in\mathcal{S}(A/I).$ Hence,
by transfinite induction, all $I_{\lambda}/I\in\mathcal{S}(A/I),$
$1\leq\lambda\leq\gamma.$ So $J/I=I_{\gamma}/I\in\mathcal{S}(A/I).$

Thus $I\ll_{\text{sm}}J.$ Hence $\ll_{\text{sm}}=$ $\ll_{\text{sm}%
}^{\triangleright}$ and it is an $\mathbf{R}$-order by Theorem \ref{inf}. So
$\mathfrak{r}_{\text{sm}}(A)=\mathfrak{r}_{\text{sm}}^{\triangleright}(A).$
\end{proof}

\begin{remark}
\label{R3}\emph{The relation-valued function }$f$\emph{ on }$\mathfrak{A}%
$\emph{ defined by }$f(A)=$ $\ll_{\text{\emph{sm}}}^{A}$ \emph{for }%
$A\in\mathfrak{A},$ \emph{illustrates Theorem \ref{T4.5}. Indeed, for each
}$A\in A,$\emph{ there is the }$\left(  \ll_{\text{\emph{sm}}}^{A}\right)
^{\triangleright}$\emph{-radical }$r_{\text{\emph{sm}}}^{\triangleright}(A)$
\emph{in} \emph{Id}$_{A}.$\emph{ However, }the map\emph{ }$r_{\text{\emph{sm}%
}}^{\triangleright}$\emph{: }$A\in\mathfrak{A}\mapsto r_{\text{\emph{sm}}%
}^{\triangleright}\left(  A\right)  $\emph{ }is not a topological
radical\emph{.}

\emph{To show this, note that condition }$(C_{2})$\emph{ in Lemma \ref{L3.5}
does not hold for }$f$\emph{. Indeed, let }$A=B(H),$\emph{ }$I=\{0\}$\emph{
and }$J=K=C(H).$\emph{ Then condition }$(C_{2})$\emph{ gives }$\{0\}\ll
_{\text{\emph{sm}}}^{C(H)}C(H)\Longleftrightarrow\{0\}\ll_{\text{\emph{sm}}%
}^{B(H)}C(H).$

\emph{However (see (\ref{5.3})), the relation }$\{0\}\ll_{\text{\emph{sm}}%
}^{C(H)}C(H)$\emph{ is not true, as }$C(H)$\emph{ is not a small ideal in the
C*-algebra }$C(H);$\emph{ while the relation }$\{0\}\ll_{\text{\emph{sm}}%
}^{B(H)}C(H)$\emph{ is true, as }$C(H)$\emph{ is a small ideal in }%
$B(H).$\emph{ Hence, by Lemma \ref{L3.5}, there does not exist a C*-property
}$P$\emph{ such that }$f=f_{_{P}}.$\emph{ So, by Theorem \ref{T4.5}, the map
}$r_{\text{\emph{sm}}}^{\triangleright}$\emph{ is not a topological radical
on} $\mathfrak{A}$.\ ${\blacksquare\bigskip}$
\end{remark}

In an arbitrary (not necessarily normed) algebra $A$, an ideal $\mathcal{J}$
is \textit{small} (see \cite{Kas}) if $\mathcal{J}+\mathcal{I}\neq A$ for any
ideal $\mathcal{I}\neq A.$ We call them K\textit{-small}. The radical
rad$_{_{K}}$($A)$ is defined in \cite[p. 214]{Kas} as
\begin{equation}
\text{rad}_{_{K}}(A)=\cap\text{ }\{\mathcal{I}\text{: }\mathcal{I}\text{ are
maximal ideals of }A\}=\sum\{\mathcal{J}\text{: }\mathcal{J}\text{ are K-small
ideals of }A\mathcal{\}}. \label{5}%
\end{equation}
The maximal ideals $\mathcal{I}$ of $A$ are not necessarily closed even if $A$
is a Banach algebra$.$

\begin{proposition}
\label{L3}Let $A$ be a unital \emph{C}$^{\ast}$-algebra. Then each small ideal
$I\in$ \emph{Id}$_{A}$ is \emph{K}-small and the closure $\overline
{\mathcal{J}}$ of each \emph{K}-small ideal $\mathcal{J}$ in $A$ is a small ideal.
\end{proposition}

\begin{proof}
If a small ideal $I$ is not K-small, $I+R=A$ for some (not necessarily closed)
ideal $R$ of $A.$ As $\mathbf{1}_{A}$\textbf{ }is surrounded by a ball of
invertible elements, $A$ has no dense ideals. Hence $\overline{R}\neq A$ and
$A=I+R\subseteq I+\overline{R}.$Thus $\overline{R}\in$ Id$_{A}$ and
$I+\overline{R}=A$ -- a contradiction, as $I\ $is small. So $I$ is K-small.

If $\mathcal{J}$ is K-small then $\mathcal{J}+I\neq A$ for all $I\in$
Id$_{A},$ $I\neq A.$ As $\mathcal{J}+I$ is an ideal of $A,$ it is not dense in
$A.$ Hence $\overline{\mathcal{J}}+I\subseteq\overline{\mathcal{J}+I}\neq A,$
so that $\overline{\mathcal{J}}$ is small.\bigskip
\end{proof}

For non-unital C*-algebras the second statement of Proposition \ref{L3} can
fail, i.e., the closure of a K-small ideal is not necessarily a small ideal.
For example, any ideal $I$ of $C(H)$ is dense$.$ So its closure $\overline
{I}=C(H)$ is not a small ideal (Example \ref{E1}). On the other hand, each
countably generated ideal of $C(H)$ is K-small, since (see \cite{BW}) if
$C(H)=I+J$ then neither of summands is countably generated.

\begin{problem}
\label{P3}\emph{Is every small ideal of C*-algebra K-small?}
\end{problem}

\begin{corollary}
\label{C4}\emph{(i)\ }$S_{A}\subseteq\mathfrak{r}_{\text{\emph{sm}}%
}^{\triangleright}(A).$ If $S_{A}$ is small\emph{, }or satisfies the
conditions of Theorem \emph{\ref{T6.1} }then $S_{A}=\mathfrak{r}%
_{\text{\emph{sm}}}^{\triangleright}(A).\smallskip$

\emph{(ii)\ }If $A$ is a unital C$^{\ast}$-algebra then $\emph{rad}_{_{K}%
}(A)=S_{A}=\mathfrak{r}_{\text{\emph{sm}}}(A)$ is a small ideal in
$A.$\emph{\smallskip}
\end{corollary}

\begin{proof}
(i) By Corollary \ref{C2}(iii), $S_{A}\subseteq\mathfrak{r}_{\text{sm}%
}^{\triangleright}(A).$ If $S_{A}\neq\mathfrak{r}_{\text{sm}}^{\triangleright
}(A)$ then, by Corollary \ref{C2}(i), there is $J\in$ Id$_{A}$ such that
$S_{A}\subsetneqq J$ and $J/S_{A}$ is small in $A/S_{A}.$

If $S_{A}$ is small\emph{, }or satisfies the conditions of Theorem \ref{T6.1},
it follows from Corollary \ref{maximal} and Theorem \ref{T6.1} that $A/S_{A}$
has no non-zero small ideals. Hence $S_{A}=\mathfrak{r}_{\text{sm}%
}^{\triangleright}(A).$

(ii) If $A$ is unital, $S_{A}$ is the largest small ideal in $A$ and $A/S_{A}$
has no small ideals by Corollary \ref{maximal}. By (i) and by Corollary
\ref{C6.1}, $S_{A}=\mathfrak{r}_{\text{sm}}^{\triangleright}(A)=\mathfrak{r}%
_{\text{sm}}(A).$ It follows from Proposition \ref{L3} that $S_{A}$ is
K-small. So $S_{A}\subseteq$ rad$_{_{K}}$($A)$ by (\ref{5}). By Proposition
\ref{L3}, for each K-small ideal $\mathcal{J},$ the ideal $\overline
{\mathcal{J}}$ is small. Thus $\mathcal{J}\subseteq S_{A}.$ So rad$_{_{K}%
}(A)\subseteq S_{A}.$ Hence rad$_{_{K}}(A)=S_{A}.$
\end{proof}

\begin{problem}
\label{P1}\emph{What is the link between rad}$_{_{K}}(A),$ $S_{A}$ \emph{and}
$\mathfrak{r}_{\text{\emph{sm}}}^{\triangleright}(A)$ \emph{if }$A$ $\emph{i}%
$\emph{s non-unital?}
\end{problem}

E. Kissin: STORM, London Metropolitan University, 166-220 Holloway Road,
London N7 8DB, Great Britain; e-mail: e.kissin@londonmet.ac.uk\medskip

V. S. Shulman: Department of Mathematics, Vologda State University, Vologda,
Russia; e-mail: shulman.victor80@gmail.com\medskip

Yu. V. Turovskii: Department of Mathematics, Ubuntu Penguin University,
Penguin, Ubuntu; e-mail: yuri.turovskii@gmail.com

\begin{thebibliography}{999999}                                                                                           %


\bibitem[AW]{AW}{\small C. Akemann and S. Wright, Compact actions on Banach
algebras, Glazgow Math.J. 21(1980) 143-149 \vspace{-0.3cm}}

\bibitem[Am]{Am}{\small S. A. Amitsur, A general theory of radicals. I.
Radicals in complete lattices, Amer. J. Math. 74 (1952) 774-786.\vspace
{-0.3cm}}

\bibitem[Am2]{Am2}{\small S. A. Amitsur, A general theory of radicals. II.
Radicals in rings and bicategories, Amer. J. Math. 76 (1954) 100-125.\vspace
{-0.3cm}}

\bibitem[Am3]{Am3}{\small S. A. Amitsur, A general theory of radicals. III.
Applications, Amer. J. Math. 76 (1954) 126-136.\vspace{-0.3cm}}

\bibitem[A]{A}{\small O. Yu. Aristov, Topological radical of a Banach module,
Studia Math., 234 (2016), 149-164.\vspace{-0.3cm}}

\bibitem[BW]{BW}{\small A. Blass and G. Weiss, A characterization and sum
decomposition for operator ideals, Trans. Amer. Math. Soc., 246 (1978),
407-417.\vspace{-0.3cm}}

\bibitem[BD]{BD}{\small F.F. Bonsall and J. Duncan, "Complete normed
algebras", Springer-Verlag, Berlin Heidelberg, New York, 1973.\vspace{-0.3cm}}

\bibitem[Br]{Br}{\small L.G. Brown, Extension of C*-algebras, Operator
Algebras and Applications, Proc. Symp. Pure Math. 38 (1981) 175-176, Amer.
Math. Soc., Providence.\vspace{-0.3cm}}

\bibitem[BP]{BP}{\small L.G. Brown and G.K. Pedersen, C*-Algebras of Real Rank
Zero, J. Funct. Anal., 99 (1991), 131-149.\vspace{-0.3cm}}

\bibitem[Da]{Da}{\small K.R. Davidson, "C*-algebras by Example", AMS,
Providence, Rhode Island, 1991.\vspace{-0.3cm}}

\bibitem[D]{D}{\small J. Dixmier, "Les C*-algebres et leurs representations",
Paris, Gauthier--Villars Editeur, 1969.\vspace{-0.3cm}}

\bibitem[Di]{Di}{\small P.G. Dixon, Topologically irreducible representations
and radicals in Banach algebras, Proc. London Math. Soc. (3), 74 (1997),
174--200.\vspace{-0.3cm}}

\bibitem[DS]{DS}{\small N. Dunford and J.T. Schwartz, "Linear operators",
Interscience Publishers, New York, 1958.\vspace{-0.3cm}}

\bibitem[Gl]{Gl}{\small J. Glimm, Type I C*-algebras, Ann. Math. 73 (1961),
572-612.\vspace{-0.3cm}}

\bibitem[G]{G}{\small G. Gr\"{a}tzer, \textquotedblleft Lattice Theory:
Foundation\textquotedblright, Birkh\"{a}user, 2011.\vspace{-0.3cm}}

\bibitem[Gr]{Gr}{\small M. Gray, "A radical approach to algebra",
Addison-Wesley, 1970.\vspace{-0.3cm}}

\bibitem[H]{H}{\small L.N. Herstein, "Noncommutative rings", John Wiley and
sons, inc., 1971.\vspace{-0.3cm}}

\bibitem[KR]{KR}{\small R. V. Kadison and J. R. Ringrose, Fundamentals of the
theory of operator algebras", v. II, Academic Press, New York, 1986.\vspace
{-0.3cm}}

\bibitem[Ka]{Ka}{\small I. Kaplansky, The structure of certain operator
algebras, Trans. AMS, 70 (1951), 219-255.\vspace{-0.3cm}}

\bibitem[Kas]{Kas}{\small F. Kasch, "Modules and rings", Academic Press,
London, 1982.\vspace{-0.3cm}}

\bibitem[K]{K}{\small E. Kirchberg, On subalgebras of the CAR-algebra, J.
Funct. Anal., 129 (1) (1995), 35-63.\vspace{-0.3cm}}

\bibitem[Ki]{Ki}{\small E. Kissin, Relations and trails in lattices of
projections in W*-algebras, Proc. Edinburgh Math. Soc. (2019), to appear, DOI:
10.1017/S0013091519000245.\vspace{-0.3cm}}

\bibitem[TR]{TR}{\small E. Kissin, V.S. Shulman, Yu. V. Turovskii,.
Topological radicals and Frattini theory of Banach Lie algebras, Integral
Equations and Operator theory, 77 (2012), 51--121.\vspace{-0.3cm}}

\bibitem[KST1]{KST1}{\small E. Kissin, V.S. Shulman and Yu. V. Turovskii, On
theory of topological radicals, Journal of Mathematical Sciences, Contemparary
Mathematics. Fundamental directions, \textbf{64}, No 3 (2018), 490-546.\vspace
{-0.3cm}}

\bibitem[KST2]{KST2}{\small E. Kissin, V.S. Shulman and Yu. V. Turovskii,
Relations and radicals in abstract lattices and in lattices of subspaces of
Banach spaces and ideals of Banach algebras. Amitsur's theory revisited,
Order, Journal on the Theory of Ordered Sets and its Applications, 38 (1)
(2021), 143-201.\vspace{-0.3cm}}

\bibitem[Murphy]{Murphy}{\small G.J. Murphy, C*-algebras and Operator Theory,
Acad. Press, Inc, New York}, {\small 1996.\vspace{-0.3cm}}

\bibitem[Sk]{Sk}{\small L. A. Skornjakov, \textquotedblleft Elements of
Lattice Theory\textquotedblright, Hindustan Publ. Corp., Delhi, 1977.\vspace
{-0.3cm}}

\bibitem[ST]{ST}{\small V.S. Shulman and Yu. V. Turovskii, Topological
radicals, V. From algebra to spectral theory, Operator Theory Advances and
applications 233, Algebraic Methods in Functional analysis (2014),
171--280.\vspace{-0.3cm}}

\bibitem[T]{T}{\small M. Takesaki, "Theory of Operator Algebras I,III",
Springer-Verlag, New-York, 2002.\vspace{-0.3cm}}

\bibitem[To]{To}{\small J. Tomiyama, A characterization of C*-algebras whose
conjugate spaces are separable, Tohoku Math. J., 15 (1963), 96-102.\vspace
{-0.3cm}}

\bibitem[W]{W}{\small S. Wasserman, "Exact C*-algebras and related topics",
Seoul National University, Research Institute of Mathematics, Global Analysis
Research Center, 1994.}
\end{thebibliography}
\end{document}